\documentclass[a4paper,12pt]{amsart}

\usepackage{float}
\usepackage{euscript,eufrak,verbatim}
\usepackage{graphicx}
\usepackage[usenames]{color}
\usepackage[colorlinks,linkcolor=red,anchorcolor=blue,citecolor=blue]{hyperref}
\usepackage{amsmath}
\usepackage{amsthm}
\usepackage[all]{xy}
\usepackage{graphicx} 
\usepackage{amssymb}

\usepackage{mathrsfs}
\usepackage{amscd}

\makeatletter
\newcommand{\svdots}{%
  \vbox{\fontsize{\sf@size}{\sf@size pt}\linespread{0.3}\selectfont
    \kern0.2\baselineskip
    \hbox{.}\hbox{.}\hbox{.}%
    \kern0.1\baselineskip
  }%
}
\makeatother

\theoremstyle{plain}
\newtheorem{main theorem}{Main Theorem}
\newtheorem{theorem}{Theorem}[section]
\newtheorem{lemma}[theorem]{Lemma}

\newtheorem{corollary}[theorem]{Corollary}
\newtheorem{proposition}[theorem]{Proposition}
\newtheorem{claim}[theorem]{Claim}

\newtheorem{lemma-definition}[theorem]{Lemma-Definition}
\theoremstyle{definition}
\newtheorem{definition}[theorem]{Definition}
\newtheorem{remark}[theorem]{Remark}
\newtheorem{example}[theorem]{Example}

\newtheorem{problem}[theorem]{Problem}

\makeatletter 
\@addtoreset{equation}{section}
\numberwithin{equation}{section}

\newcommand{\norm}[1]{\left\lVert#1\right\rVert}

\newcommand{\diam}{\mathrm{Diam}}

\newcommand{\mdim}{\mathrm{mdim}}

\newcommand{\widim}{\mathrm{Widim}}

\newcommand{\umdimm}{\overline{\mathrm{mdim}}_{\mathrm{M}}}
\newcommand{\lmdimm}{\underline{\mathrm{mdim}}_{\mathrm{M}}}
\newcommand{\mdimm}{\mathrm{mdim}_{\mathrm{M}}}

\newcommand{\urdim}{\overline{\mathrm{rdim}}}
\newcommand{\lrdim}{\underline{\mathrm{rdim}}}
\newcommand{\rdim}{\mathrm{rdim}}

\newcommand{\dfs}{d_{\mathrm{FS}}}

\title[Rate distortion dimension of random Brody curves]
{Rate distortion dimension of random Brody curves}

\author{Masaki Tsukamoto}

\address
{Department of Mathematics, Kyoto University, Kitashirakawa Oiwake-cho, Sakyo-ku, Kyoto 606-8502, Japan}

\email{tsukamoto@math.kyoto-u.ac.jp}

\begin{document}



\thanks{The author was supported by JSPS KAKENHI JP21K03227.}

\maketitle

\begin{abstract}
The main purpose of this paper is to propose an ergodic theoretic approach to the study of entire holomorphic curves.
Brody curves are one-Lipschitz holomorphic maps from the complex plane to the complex projective space. 
They naturally form a dynamical system, and 
“random Brody curves” in the title refers to invariant probability measures on it.
We study their geometric and dynamical properties.
Given an invariant probability measure $\mu$ on the space of Brody curves, our first main theorem claims that 
its rate distortion dimension is bounded by the integral of a “potential function” over $\mu$.
This result is analogous to the Ruelle inequality of smooth ergodic theory.
Our second main theorem claims that there exists a rich variety of invariant probability measures attaining 
equality in this “Ruelle inequality for Brody curves”.
The main tools of the proofs are the deformation theory of Brody curves and the variational principle for 
mean dimension with potential.
This approach is motivated by the theory of thermodynamic formalism for Axiom A diffeomorphisms.
\end{abstract}

\section{Backgrounds}  \label{section: Backgrounds}

The classical Nevanlinna theory studies the asymptotic value distribution of meromorphic functions in the complex plane $\mathbb{C}$.
Meromorphic functions are holomorphic maps from $\mathbb{C}$ to the Riemann sphere $\mathbb{C}\cup \{\infty\}$.
Hence it is natural to generalize it to the study of holomorphic curves in complex manifolds.
Over almost a century, several generations of mathematicians have been developing the higher-dimensional Nevanlinna theory.
Substantial progress has been made.
However there still remain many fundamental questions. 
Readers can see a good overview in the book of Noguchi and Winkelmann \cite{Noguchi--Winkelmann}.

The main purpose of this paper is to propose a new approach to the study of holomorphic curves.
We develop an ergodic theoretic approach.
Broadly speaking, the Nevanlinna theory and ergodic theory have a similar spirit;
they both study \textit{recurrence} of functions and orbits. 
Therefore one may imagine some \textit{ergodic theory of holomorphic curves}.
This paper is a serious attempt to such a theory.
We explain backgrounds of the paper in this section. 
Our main results will be given in \S \ref{section: main results}.

There are two main backgrounds for this paper.
We will review them in \S \ref{subsection: ergodic theory of Axiom A diffeomorphisms} and 
\S \ref{subsection: mean dimension of the dynamical system of Brody curves} below.
One is the ergodic theory of Axiom A diffeomorphisms.
This is a milestone of smooth ergodic theory developed by several distinguished mathematicians,
including Anosov, Smale, Adler, Weiss, Sinai, Ruelle and Bowen \cite{Anosov, Smale, Adler--Weiss, Sinai, Ruelle73, Bowen_lecture}.  
Another motivation comes from the study of mean dimension of 
the dynamical system of Brody curves.
Brody curves are one-Lipschitz holomorphic maps from $\mathbb{C}$ to the complex projective space $\mathbb{C}P^N$.
They form a dynamical system and Gromov \cite{Gromov} began to study it.
Roughly speaking, our main theorems show that the dynamical system of Brody curves has an ergodic theoretic 
structure analogous to that of Axiom A diffeomorphisms.

\subsection{Ergodic theory of Axiom A diffeomorphisms} \label{subsection: ergodic theory of Axiom A diffeomorphisms}

In this subsection we briefly review the ergodic theory of Axiom A diffeomorphisms \cite{Bowen_lecture, Katok--Hasselblatt}.
Although this provides a very important motivation for this paper, it is logically independent of 
the main results below.
So readers may skip this subsection if they are unfamiliar to hyperbolic dynamics. 
Our exposition follows the book of Bowen \cite[Chapters 3 and 4]{Bowen_lecture}.

Let $M$ be a compact $C^\infty$ Riemannian manifold and let $T\colon M\to M$ be a diffeomorphism.
We assume that $T$ satisfies \textbf{Axiom A}, namely the non-wandering set is hyperbolic and periodic points are dense in it.
Let $\Omega$ be a \textbf{basic set} of $(M, T)$.
A simple example is that $M = \mathbb{T}^2 := \mathbb{R}^2/\mathbb{Z}^2$ is the two-dimensional torus and 
$T$ is a hyperbolic toral automorphism, e.g. 
$T(x, y) = (x+y, x)$. In this case, the set of periodic points are dense in $\mathbb{T}^2$, 
and the non-wandering set is equal to the whole space $\mathbb{T}^2$.
The basic set $\Omega$ is also equal to $\mathbb{T}^2$.

We denote by $T_x M = E^s_x \oplus E^u_x$ $(x\in \Omega)$ the splitting of the tangent space into 
stable and unstable directions.
We define a function $\phi \colon \Omega \to \mathbb{R}$ as the logarithm of the Jacobian of the map $T$ in the 
unstable direction: 
\begin{equation} \label{eq: unstable Jacobian}
    \phi(x) = \log \left|\det\left(dT_x\colon E^u_x\to E^u_{Tx}\right)\right|. 
\end{equation}    
Let $\mu$ be a $T$-invariant Borel probability measure on $\Omega$.
We denote its Kolmogorov--Sinai entropy by $h_\mu(T)$. Then we have 
the following fundamental inequality \cite[Proposition 4.8 (b)]{Bowen_lecture}:
\begin{equation}  \label{eq: Ruelle inequality}
   h_\mu(T) \leq  \int_\Omega \phi(x) \, d\mu(x). 
\end{equation}
This inequality claims that the entropy is bounded by the average rate of the expansion of the map $T$.
If $\Omega$ is an \textbf{attractor}, then there exists a $T$-invariant Borel probability measure $\mu^+$ on $\Omega$
(called a \textbf{Sinai--Ruelle--Bowen measure}) attaining the equality \cite[Theorem 4.11]{Bowen_lecture}
\begin{equation} \label{eq: Sinai--Ruelle--Bowen}
   h_{\mu^+}(T) = \int_M \phi(x) \, d\mu^+(x). 
\end{equation}    

If $M = \mathbb{T}^2$ and $T(x, y) = (x+y, x)$, then 
$\mathbb{R}^2 = \mathbb{R}\left(\frac{1-\sqrt{5}}{2}, 1\right)\oplus \mathbb{R}\left(\frac{1+\sqrt{5}}{2}, 1\right)$
provides the splitting into stable and unstable directions.
The function $\phi(x)$ is equal to the constant $\log \left(\frac{1+\sqrt{5}}{2}\right)$.
The above inequality (\ref{eq: Ruelle inequality}) gives $h_\mu(T) \leq \log \left(\frac{1+\sqrt{5}}{2}\right)$
for all $T$-invariant Borel probability measures $\mu$ on $\mathbb{T}^2$.
The equality holds if and only if $\mu$ is the Lebesgue measure on $\mathbb{T}^2$.

We notice that the inequality (\ref{eq: Ruelle inequality}) is a special case of the Ruelle inequality
(\cite{Ruelle78}, \cite[Theorem S.2.13]{Katok--Hasselblatt}).
For a diffeomorphism $T$ of a compact Riemannian manifold $M$ (not necessarily satisfying Axiom A),
the Ruelle inequality claims that 
\begin{equation} \label{eq: general Ruelle inequality}
  h_\mu(T) \leq  \int_M \chi_+(x) \, d\mu(x)  
\end{equation}  
for all $T$-invariant  Borel probability measures $\mu$ on $M$.
Here $\chi_+(x)$ is the sum of positive \textbf{Lyapunov exponents} counted with multiplicity.
If the map $T$ satisfies Axiom A and the measure $\mu$ is supported on a basic set $\Omega$, then the integral 
$\int_M \chi_+(x) \, d\mu(x)$ is equal to $\int_\Omega \phi(x) \, d\mu(x)$ in (\ref{eq: Ruelle inequality}) by the ergodic theorem.
So (\ref{eq: Ruelle inequality}) is a special case of (\ref{eq: general Ruelle inequality}).
Our first main result (Theorem \ref{theorem: Ruelle inequality for Brody curves} below) 
is an analogue of the inequality (\ref{eq: Ruelle inequality}) for Brody curves.

\subsection{Mean dimension of the dynamical system of Brody curves} 
\label{subsection: mean dimension of the dynamical system of Brody curves}

We denote by $z = x+y\sqrt{-1}$ the standard coordinate of the complex plane $\mathbb{C}$.
Let $\mathbb{C}P^N$ be the complex projective space equipped with the Fubini--Study metric.
Let $f\colon \mathbb{C}\to \mathbb{C}P^N$ be a holomorphic map. 
We define $|df|(z)$ as the \textbf{local Lipschitz constant} of $f$ at $z$.
Namely, for a tangent vector $v\in T_z\mathbb{C}$, the Fubini--Study length of $df(v)\in T_{f(z)} \mathbb{C}P^N$ is given by 
$|df(v)| = |df|(z)\cdot  |v|$, where $|v|$ is the Euclidean length of $v$.  
More explicitly, for $f(z) = [f_0(z): f_1(z): \dots: f_N(z)]$ in the homogeneous coordinate (each $f_i(z)$ is a holomorphic function),
the local Lipschitz constant is give by 
\begin{equation} \label{eq: local Lipschitz constant}
    |df|^2(z) = \frac{1}{4\pi} \left(\frac{\partial^2}{\partial x^2} + \frac{\partial^2}{\partial y^2}\right) 
     \log \left(|f_0(z)|^2 + |f_1(z)|^2+\cdots +|f_N(z)|^2\right). 
\end{equation}     
The right-hand side is always nonnegative. $|df|(z)$ is its (nonnegative) square root.

We call a holomorphic map $f\colon \mathbb{C}\to \mathbb{C}P^N$ a \textbf{Brody curve} if $|df|(z) \leq 1$ for all $z\in \mathbb{C}$.
An importance of this notion stems from the following fundamental theorem of Brody \cite{Brody}: 
\textit{A compact complex manifold is Kobayashi-hyperbolic 
if and only if it does not admit any non-constant Lipschitz holomorphic map from $\mathbb{C}$}.
See also the papers of Duval \cite{Duval, Duval_around} for very interesting variants of this \textit{Bloch--Brody principle}.
One can find more information about Brody curves in the book \cite[\S 7.6]{Noguchi--Winkelmann} and 
papers \cite{Eremenko, Costa, Costa--Duval, Winkelmann}.

Let $\mathcal{B}^N$ be the set of all Brody curves in $\mathbb{C}P^N$.
This becomes a compact metrizable space under the compact-open topology. (A sequence $\{f_n\}$ of Brody curves converges to 
$f\in \mathcal{B}^N$ if $f_n$ uniformly converges to $f$ over every compact subset of $\mathbb{C}$.
The compactness of $\mathcal{B}^N$ follows from the Arzel\`{a}--Ascoli theorem because Brody curves are one-Lipschitz.)

The group $\mathbb{C}$ continuously acts on $\mathcal{B}^N$ by 
\begin{equation}  \label{eq: group action on Brody curves}
   T\colon  \mathbb{C}\times \mathcal{B}^N \to  \mathcal{B}^N, \quad \left(a, f(z)\right) \mapsto  f(z+a). 
\end{equation}   
Gromov \cite{Gromov} began to study the mean dimension $\mdim(\mathcal{B}^N, T)$ of this action.
This is the number of parameters per unit area of the plane for describing the orbits of $\mathcal{B}^N$.
We review the precise definition of mean dimension in 
\S \ref{subsection: definitions and variational principle}.
Gromov showed that (\cite[p. 396]{Gromov}, \cite{Eremenko}) 
\[  \mdim(\mathcal{B}^N, T) \leq  4N. \]
This is not sharp.
For improving Gromov’s estimate (and, indeed, obtaining the exact formula),
we introduce the \textbf{energy density} $\rho(f)$ of Brody curves $f\in \mathcal{B}^N$ by 
\begin{equation}  \label{eq: energy density}
  \rho(f) = \lim_{R\to \infty}  \left(\sup_{a\in \mathbb{C}} \frac{1}{\pi R^2} \int_{|z-a|<R} |df|^2(z) \, dxdy \right). 
\end{equation}
We define $\rho(\mathbb{C}P^N)$ as the supremum of $\rho(f)$ over $f\in \mathcal{B}^N$.
It is known that \cite{Tsukamoto_packing}
\[  0< \rho(\mathbb{C} P^N) < 1, \quad  \lim_{N\to \infty} \rho(\mathbb{C} P^N) = 1. \]

The papers \cite{Matsuo--Tsukamoto, Tsukamoto_mean_dimension_Brody_curves} proved that
the mean dimension is given by 
\begin{equation} \label{eq: mean dimension of Brody curves}
   \mdim(\mathcal{B}^N, T) = 2(N+1) \rho(\mathbb{C}P^N). 
\end{equation}   
A main purpose of this paper is to get a deeper understanding of this formula in terms of invariant probability measures.

\section{Main results}  \label{section: main results}

The main object of this paper is invariant probability measures on the space of Brody curves.
Let $\mathcal{B}^N$ be the space of Brody curves in $\mathbb{C}P^N$ equipped with 
the continuous action of $\mathbb{C}$ defined by (\ref{eq: group action on Brody curves}).
A Borel probability measure $\mu$ on $\mathcal{B}^N$ is said to be \textbf{$T$-invariant} if 
$\mu(T^a A) = \mu(A)$ for all $a\in \mathbb{C}$ and all Borel subsets $A\subset \mathcal{B}^N$.
We define $\mathscr{M}^T(\mathcal{B}^N)$ as the set of all $T$-invariant Borel probability measures on $\mathcal{B}^N$.

It is not very common to study invariant probability measures in the context of holomorphic curve theory.
So, before going into the details, we present an example for better understanding the subject:

\begin{example} \label{example: random Brody curves}
Let $L$ and $a$ be large positive numbers.
Set 
\[  \Lambda = \mathbb{Z} L + \mathbb{Z} L\sqrt{-1} , \quad K = [0, L]^2, \quad  
     D = \{u \in \mathbb{C}\mid |u-a| \leq  1\}. \]
$\Lambda$ is a lattice in $\mathbb{C}$ and $K$ is a fundamental domain of $\Lambda$. 
Let $D^{\Lambda}$ be the product space of the infinite copies of $D$ indexed by $\Lambda$.
For $u = (u_\lambda)_{\lambda\in \Lambda} \in D^{\Lambda}$ we define a meromorphic function by
\[  f_u(z) = \sum_{\lambda\in \Lambda} \frac{u_\lambda}{(z-\lambda)^3}. \]
It is easy to check that $f_u(z)$ belongs to $\mathcal{B}^1$ as a map from $\mathbb{C}$ to $\mathbb{C}P^1$ if $L$ and $a$
are chosen appropriately large (e.g. $L \gg a^3 \gg 1$).
Then we can consider a continuous map 
\[ \Phi\colon  K\times D^{\Lambda} \to  \mathcal{B}^1, \quad 
    (w, u) \mapsto  f_u(z+w). \]
The image of this map is a $T$-invariant closed subset of $\mathcal{B}^1$.
Let $\mathbf{m}$ be the two-dimensional Lebesgue measure on the plane.
We define probability measures $\mathbf{m}_1$ and $\mathbf{m}_2$ on $K$ and $D$ respectively by 
\[  \mathbf{m}_1 = \frac{\mathbf{m}|_K}{L^2}, \quad \mathbf{m}_2 = \frac{\mathbf{m}|_D}{\pi}, \]
where $\mathbf{m}|_K$ and $\mathbf{m}|_D$ denote the restrictions of the Lebesgue measure $\mathbf{m}$
to $K$ and $D$.
We consider the product measure $\mathbf{m}_1\otimes \mathbf{m}_2^{\otimes \Lambda}$ on $K\times D^{\Lambda}$
and define 
\[  \mu = \Phi_*\left(\mathbf{m}_1\otimes \mathbf{m}_2^{\otimes \Lambda}\right)  \]
as the push-forward measure of $\mathbf{m}_1\otimes \mathbf{m}_2^{\otimes \Lambda}$ under the map $\Phi$.
The measure $\mu$ is $T$-invariant and hence belongs to $\mathscr{M}^T(\mathcal{B}^1)$.
Considering the measure $\mu$ is equivalent to considering a random function
\[    f_u(z+w) = \sum_{\lambda\in \Lambda} \frac{u_\lambda}{(z+w - \lambda)^3}, \]
where $w$ and $u_\lambda$ $(\lambda\in \Lambda)$ are independently chosen according to the uniform distributions on $K$ and $D$.
In general, the study of invariant probability measures on $\mathcal{B}^N$ is equivalent to the study of such
\textit{random Brody curves}.
\end{example}

Our starting point is the formula $\mdim(\mathcal{B}^N, T) = 2(N+1)\rho(\mathbb{C}P^N)$ given in (\ref{eq: mean dimension of Brody curves}).
We will see below that both sides of this formula can be expressed in terms of invariant probability measures.

First we consider the right-hand side $2(N+1)\rho(\mathbb{C}P^N)$.
We define a continuous function $\psi\colon \mathcal{B}^N\to \mathbb{R}$ by 
\begin{equation} \label{eq: potential function for Brody curves}
    \psi(f) = 2(N+1)\,  |df|^2(0). 
\end{equation}   
Here $|df|^2(0)$ is the square of the local Lipschitz constant of $f$ at $z=0$ defined by (\ref{eq: local Lipschitz constant}).
Now we have (see Proposition \ref{prop: energy density and invariant measures} in \S \ref{section: energy integral of Brody curves})
\begin{equation}  \label{eq: energy density and invariant measures}
   2(N+1) \rho(\mathbb{C}P^N) = \sup_{\mu \in \mathscr{M}^T(\mathcal{B}^N)} \int_{\mathcal{B}^N} \psi\, d\mu. 
\end{equation}
Therefore the right-hand side of the formula $\mdim(\mathcal{B}^N, T) = 2(N+1)\rho(\mathbb{C}P^N)$ 
is expressed by using the integral of $\psi$ over invariant 
probability measures.

Next we consider the mean dimension $\mdim(\mathcal{B}^N, T)$. 
For connecting it to invariant measures, we need to consider \textit{rate distortion theory} \cite{Cover--Thomas}.
Rate distortion theory is the basic theory for lossy date compression method.
Here ``lossy’’ means that some data is lost in the process of compression.
For example, given an image (a two-dimensional signal), 
JPEG (Joint Photographic Experts Group) expands it in a wavelet basis and discard small coefficients
and quantizes the remaining ones. Some information is inevitably lost in this process, but the resulting image often 
looks similar to the original one at least for human eyes.
Rate distortion theory describes a fundamental limit of such lossy data compression.

We define a metric $\mathbf{d}$ on $\mathcal{B}^N$ (compatible with the given topology)
by\footnote{This $\mathbf{d}(f, g)$ defines a metric on $\mathcal{B}^N$ because of the 
unique continuation principle; if $\mathbf{d}(f, g) = 0$ then $f=g$ all over the plane. We can also consider other choices of metrics.
For example, 
\[ \mathbf{d}^\prime(f, g) = \max_{|z|\leq 1} d_{\mathrm{FS}}(f(z), g(z))  \quad \text{or} \quad
 \mathbf{d}^{\prime\prime}(f, g) = \sum_{n=1}^\infty 2^{-n} \max_{|z|\leq n}d_{\mathrm{FS}}(f(z), g(z)).  \]
Our main results also hold for $\mathbf{d}^\prime$ and $\mathbf{d}^{\prime\prime}$. But we use $\mathbf{d}$ for simplicity.}
\begin{equation} \label{eq: metric on B^N}
    \mathbf{d}(f, g) = \max_{z\in [0,1]^2} d_{\mathrm{FS}}(f(z), g(z)).
\end{equation}    
Here $d_{\mathrm{FS}}$ is a distance function on $\mathbb{C}P^N$ induced by the Fubini--Study metric, and 
$[0,1]^2 = \{x+y\sqrt{-1}\mid 0\leq x, y \leq 1\} \subset \mathbb{C}$.
Let $\mu \in \mathscr{M}^T(\mathcal{B}^N)$.
For a positive number $\varepsilon$, we denote by $R(\mathbf{d}, \mu, \varepsilon)$ the rate distortion function
with respect to the metric $\mathbf{d}$ and measure $\mu$.
The value of $R(\mathbf{d}, \mu, \varepsilon)$ is a nonnegative real number.
Its precise definition will be given in \S \ref{subsection: rate distortion theory}.
The intuitive meaning of $R(\mathbf{d},\mu,\varepsilon)$ is as follows. 
Suppose that we randomly take a Brody curve $f\in \mathcal{B}^N$ according to the distribution $\mu$.
The rate distortion function $R(\mathbf{d}, \mu, \varepsilon)$ 
measures the bits per unit area of the plane for describing $f$ within 
averaged distortion bounded by $\varepsilon$ with respect to $\mathbf{d}$.

We define \textbf{upper and lower rate distortion dimensions} of $(\mathcal{B}^N, T, \mathbf{d}, \mu)$ by 
\[  \urdim\left(\mathcal{B}^N, T, \mathbf{d}, \mu\right) = \limsup_{\varepsilon\to 0} \frac{R(\mathbf{d}, \mu, \varepsilon)}{\log(1/\varepsilon)},
   \quad   \lrdim\left(\mathcal{B}^N, T, \mathbf{d}, \mu\right) = 
   \liminf_{\varepsilon\to 0} \frac{R(\mathbf{d}, \mu, \varepsilon)}{\log(1/\varepsilon)}.  \]
These are nonnegative real numbers.
If the upper and lower limits coincide, we denote their common value by 
$\rdim \left(\mathcal{B}^N, T, \mathbf{d}, \mu\right)$.
The rate distortion dimension was first introduced by Kawabata--Dembo \cite{Kawabata--Dembo}
in a context of general information processing.
Their original motivation is to study the asymptotic behavior (as $\varepsilon \to 0$) of rate distortion function of
signals taking values in a fractal set.

The mean dimension $\mdim(\mathcal{B}^N, T)$
is given by\footnote{This follows from the argument of \cite[\S 3]{Tsukamoto_mean_dimension_Brody_curves}
and the variational principle given in Theorem \ref{theorem: variational principle for mean dimension}.
See also Remark \ref{remark: ordinary topological and metric mean dimension of Brody curves} below.}
\begin{equation} \label{eq: mdim of B^N and rate distortion dimension}
    \mdim(\mathcal{B}^N, T) = \sup_{\mu\in \mathscr{M}^T(\mathcal{B}^N)} \lrdim\left(\mathcal{B}^N, T, \mathbf{d}, \mu\right)
    = \sup_{\mu \in \mathscr{M}^T(\mathcal{B}^N)} \urdim\left(\mathcal{B}^N, T, \mathbf{d}, \mu\right).
\end{equation}    
From the formula $\mdim(\mathcal{B}^N, T) = 2(N+1)\rho(\mathbb{C}P^N)$ with
 (\ref{eq: energy density and invariant measures}) and (\ref{eq: mdim of B^N and rate distortion dimension}), we conclude 
\begin{equation}  \label{eq: mean dimension formula via invariant measures}
   \sup_{\mu\in \mathscr{M}^T(\mathcal{B}^N)} \lrdim\left(\mathcal{B}^N, T, \mathbf{d}, \mu\right)
    = \sup_{\mu \in \mathscr{M}^T(\mathcal{B}^N)} \urdim\left(\mathcal{B}^N, T, \mathbf{d}, \mu\right)
    = \sup_{\mu \in \mathscr{M}^T(\mathcal{B}^N)} \int_{\mathcal{B}^N} \psi\, d\mu.  
\end{equation}    
This is a reformulation of the equation $\mdim(\mathcal{B}^N, T) = 2(N+1)\rho(\mathbb{C}P^N)$ in terms of 
invariant probability measures.
Now a question naturally arises: 
\begin{problem}[Main problem]
 What is a relation between the integral 
 $\int_{\mathcal{B}^N} \psi\, d\mu$ and the rate distortion dimensions $\lrdim\left(\mathcal{B}^N, T, \mathbf{d}, \mu\right)$
 and $\urdim\left(\mathcal{B}^N, T, \mathbf{d}, \mu\right)$ for each individual invariant 
 measure $\mu$?
\end{problem}
If we take the supremum over $\mu\in \mathscr{M}^T\left(\mathcal{B}^N\right)$ then we have 
the equality (\ref{eq: mean dimension formula via invariant measures}).
But we want to know a more precise result for each measure.

\begin{example} \label{example: naive example}
One might naively expect an equality such as $\urdim\left(\mathcal{B}^N, T, \mathbf{d}, \mu\right) = \int_{\mathcal{B}^N}\psi\, d\mu$.
But this does not hold in general.
Let $L$ be a large positive number and set
$\Lambda = \mathbb{Z}L + \mathbb{Z}L\sqrt{-1}$.
Let $f\colon \mathbb{C}\to \mathbb{C}P^N$ be a $\Lambda$-periodic non-constant Brody curve.
Here “$\Lambda$-periodic” means that $f(z+\lambda) = f(z)$ for all $\lambda \in \Lambda$.
Such $f$ can be constructed by using elliptic functions.
Let $\delta_f$ be the delta measure at $f$, and we define $\mu \in \mathscr{M}^T(\mathcal{B}^N)$ by 
\[  \mu = \frac{1}{L^2}\int_{[0,L]^2} T^u_*\delta_f\, d\mathbf{m}(u), \]
where $\mathbf{m}$ is the two-dimensional Lebesgue measure on $\mathbb{C}$.
Therefore $\mu$ is the uniform measure on the periodic orbit of $f$.
Then we have 
\[  \rdim\left(\mathcal{B}^N, T, \mathbf{d}, \mu\right) = 0, \quad 
     \int_{\mathcal{B}^N} \psi\, d\mu = \frac{2(N+1)}{L^2} \int_{[0,L]^2} |df|^2 dxdy >0. \]
This shows that the equality $\rdim\left(\mathcal{B}^N, T, \mathbf{d}, \mu\right) = \int_{\mathcal{B}^N}\psi\, d\mu$
does not hold in general.
\end{example}

Our main results are the following two theorems.

\begin{theorem} \label{theorem: Ruelle inequality for Brody curves}
For every invariant probability measure $\mu\in \mathscr{M}^T(\mathcal{B}^N)$ we have 
  \begin{equation}   \label{eq: Ruelle inequality for Brody curves}
     \urdim\left(\mathcal{B}^N, T, \mathbf{d}, \mu\right) \leq  \int_{\mathcal{B}^N} \psi \, d\mu. 
  \end{equation}
\end{theorem}

\begin{theorem}  \label{theorem: construction of extremal measure}
For any real number $c$ with $0\leq c < 2(N+1)\rho(\mathbb{C}P^N)$ there exists 
an invariant probability measure $\mu \in \mathscr{M}^T(\mathcal{B}^N)$ satisfying
   \begin{equation}  \label{eq: construction of extremal measure}
     \rdim\left(\mathcal{B}^N, T, \mathbf{d}, \mu\right) = \int_{\mathcal{B}^N} \psi\, d\mu = c. 
   \end{equation}     
Notice that this includes the claim that the upper and lower rate distortion dimensions coincide for this measure $\mu$.   
\end{theorem}

Recall that we have 
\[  2(N+1)\rho(\mathbb{C}P^N) = \sup_{\mu \in \mathscr{M}^T(\mathcal{B}^N)} \int_{\mathcal{B}^N} \psi\, d\mu. \]
Therefore Theorem \ref{theorem: construction of extremal measure} claims that, for every 
$0 \leq c < \sup_{\mu \in \mathscr{M}^T(\mathcal{B}^N)} \int_{\mathcal{B}^N} \psi\, d\mu$, 
there exists $\mu \in \mathscr{M}^T(\mathcal{B}^N)$ satisfying (\ref{eq: construction of extremal measure}).
Now it is immediate to see that Theorems \ref{theorem: Ruelle inequality for Brody curves} and \ref{theorem: construction of extremal measure}
imply the above equation (\ref{eq: mean dimension formula via invariant measures}):
\[ \sup_{\mu\in \mathscr{M}^T(\mathcal{B}^N)} \lrdim\left(\mathcal{B}^N, T, \mathbf{d}, \mu\right)
    = \sup_{\mu \in \mathscr{M}^T(\mathcal{B}^N)} \urdim\left(\mathcal{B}^N, T, \mathbf{d}, \mu\right)
    = \sup_{\mu \in \mathscr{M}^T(\mathcal{B}^N)} \int_{\mathcal{B}^N} \psi\, d\mu.  \]
Thus Theorems \ref{theorem: Ruelle inequality for Brody curves} and \ref{theorem: construction of extremal measure}
provide a deeper understanding of the formula
$\mdim\left(\mathcal{B}^N, T\right) = 2(N+1)\rho(\mathbb{C}P^N)$ in terms of invariant probability measures.

The inequality (\ref{eq: Ruelle inequality for Brody curves}) of
Theorem \ref{theorem: Ruelle inequality for Brody curves} is an analogue of the Ruelle inequality 
$h_\mu(T) \leq \int_\Omega \phi \, d\mu$ of Axiom A diffeomorphisms (\ref{eq: Ruelle inequality}).
Both the inequalities 
\[  h_\mu(T) \leq \int_\Omega \phi \, d\mu  \quad   \text{and} \quad 
     \urdim\left(\mathcal{B}^N, T, \mathbf{d}, \mu\right) \leq  \int_{\mathcal{B}^N} \psi \, d\mu  \]
claim that the information theoretic quantities are bounded by the integral of certain “potential functions”.
We will see later that there is a similarity not only in their statements (\ref{eq: Ruelle inequality}) and (\ref{eq: Ruelle inequality for Brody curves})
but also in their proofs.
So we may say that Theorem \ref{theorem: Ruelle inequality for Brody curves} is a \textit{Ruelle inequality for Brody curves}.

Theorem \ref{theorem: construction of extremal measure} shows that there is also an important difference 
between Axiom A diffeomorphisms and Brody curves.
In the case of Axiom A attractors, the Sinai--Ruelle--Bowen measure (an invariant measure attaining the equality
$h_\mu(T) = \int_\Omega \phi \, d\mu$) is unique \cite{Bowen_unique}.
However Theorem \ref{theorem: construction of extremal measure} shows that there are uncountably many invariant measures 
$\mu$ attaining the equality $\rdim\left(\mathcal{B}^N, T, \mathbf{d}, \mu\right) = \int_{\mathcal{B}^N} \psi\, d\mu$.
It seems that there is no way to select one distinguished invariant measure for Brody curves.

\begin{remark} \label{remark: ergodic theorem and characteristic function}
The integral $\int_{\mathcal{B}^N} \psi \, d\mu$ might look mysterious.
Here we explain that it has a natural connection to the standard Nevanlinna theory.
For a holomorphic map $f\colon \mathbb{C}\to \mathbb{C}P^N$ we define 
its \textbf{Nevanlinna--Shimizu--Ahlfors characteristic function} 
\cite[p. 9, p. 73]{Noguchi--Winkelmann} by
\[  T(R, f) = \int_1^R \left(\int_{|z|<r} |df|^2 \, dxdy\right) \frac{dr}{r}, \quad (R>1). \]
Let $\mu\in \mathscr{M}^T\left(\mathcal{B}^N\right)$.
We assume that $\mu$ is \textit{ergodic}.
This means that if a Borel set $A\subset \mathcal{B}^N$ satisfies 
$\mu\left(T^u A \triangle A\right) = 0$ for all $u\in \mathbb{C}$ then $\mu(A)$ is equal to $0$ or $1$.
Here $\triangle$ denotes the symmetric difference.
Every invariant probability measure can be decomposed into ergodic measures
(Ergodic Decomposition Theorem \cite[Theorem 8.20]{Einsiedler--Ward}).
The pointwise ergodic theorem implies that \cite[Theorem 8.19]{Einsiedler--Ward}
\[  \lim_{r\to \infty} \frac{1}{\pi r^2} \int_{|u|<r} \psi\left(T^u f\right) \, d\mathbf{m}(u)
     = \int_{\mathcal{B}^N} \psi\, d\mu \]
for $\mu$-almost every $f\in \mathcal{B}^N$.
We have
\[  \int_{|u|<r} \psi\left(T^u f\right) \, d\mathbf{m}(u)  
      = 2(N+1) \int_{|z|<r} |df|^2\, dxdy. \]
Therefore 
\[    \lim_{r\to \infty} \frac{2(N+1)}{\pi r^2}\int_{|z|<r} |df|^2\, dxdy  = \int_{\mathcal{B}^N} \psi\, d\mu \]
for $\mu$-almost every $f\in \mathcal{B}^N$.
Then we have 
\[  T(R, f) = \frac{\pi R^2}{4(N+1)}\int_{\mathcal{B}^N} \psi\, d\mu  + o(R^2) 
     \quad   \text{as $R\to \infty$} \]
for $\mu$-almost every $f\in \mathcal{B}^N$.
Thus the integral $\int_{\mathcal{B}^N} \psi \, d\mu$ 
describes the coefficient of the leading term 
of $T(R, f)$ for $\mu$-almost every Brody curve $f$.

For general $\mu\in \mathscr{M}^T\left(\mathcal{B}^N\right)$ (not necessarily ergodic),
$\frac{4(N+1)}{\pi R^2} T(R, f)$ still converges as $R\to \infty$ for $\mu$-almost every $f\in \mathcal{B}^N$,
and its expected value is equal to $\int_{\mathcal{B}^N} \psi \, d\mu$:
\[   \int_{\mathcal{B}^N} \left(\lim_{R\to \infty} \frac{4(N+1)}{\pi R^2} T(R, f)\right) \, d\mu(f)  = \int_{\mathcal{B}^N} \psi \, d\mu. \]
\end{remark}

\begin{remark}
There are several other choices of the “potential function” $\psi$.
For example, for $f\in \mathcal{B}^N$, let
\[  \psi_1(f) = 2(N+1)\int_{[0,1]^2}|df|^2\, dxdy, \quad \psi_2 (f) = \frac{2(N+1)}{\pi}\int_{|z|<1} |df|^2\, dxdy. \]
Then for any $\mu \in \mathscr{M}^T(\mathcal{B}^N)$, we have 
\[ \int_{\mathcal{B}^N} \psi\, d\mu = \int_{\mathcal{B}^N}\psi_1\, d\mu = \int_{\mathcal{B}^N} \psi_2\, d\mu. \]
Therefore we can also formulate our main results by using $\psi_1$ or $\psi_2$.
We decided to use the function $\psi(f) = 2(N+1) |df|^2(0)$ because this seems the simplest choice.
\end{remark}

\begin{example}[Continuation of Example \ref{example: random Brody curves}] \label{example: random Brody curves continued}
Let $\mu\in \mathscr{M}^T(\mathcal{B}^1)$ be the measure constructed in Example \ref{example: random Brody curves}.
Namely $\mu$ is the distribution of a random function
\begin{equation}  \label{eq: random sum}
     \sum_{\lambda\in \Lambda} \frac{u_\lambda}{(z+w - \lambda)^3}, \quad  (\Lambda = \mathbb{Z} L + \mathbb{Z}L\sqrt{-1})  
\end{equation}     
where $w$ and $u_\lambda$ are independently chosen from 
the uniform distributions on $[0, L]^2$ and $\{u \in \mathbb{C} \mid   |u-a|\leq  1\}$ respectively.
Then we have 
\[  \rdim\left(\mathcal{B}^1, T, \mathbf{d}, \mu\right) = \frac{2}{L^2}, \quad 
     \int_{\mathcal{B}^1} \psi \, d\mu = \frac{12}{L^2}. \]
These certainly satisfy the “Ruelle inequality” $\rdim\left(\mathcal{B}^1, T, \mathbf{d}, \mu\right) \leq  \int_{\mathcal{B}^1} \psi \, d\mu$.
\end{example}

\begin{problem} \label{problem: explicit extremal family}
Theorem \ref{theorem: construction of extremal measure} shows that there exist many invariant measures 
$\mu \in \mathscr{M}^T(\mathcal{B}^N)$ attaining the equality 
$\rdim\left(\mathcal{B}^N, T, \mathbf{d}, \mu\right) = \int_{\mathcal{B}^N} \psi\, d\mu$.
However it is not easy to investigate their properties (e.g. ergodicity, mixing property) from the construction given in the proof of 
Theorem \ref{theorem: construction of extremal measure}. 
It is desirable to have a more elementary example.
Can one construct an explicit example\footnote{Of course, we exclude the trivial example $\mu = \delta_f$ with a constant curve $f$.
Such a delta measure satisfies $\rdim\left(\mathcal{B}^N, T, \mathbf{d}, \mu\right) = \int_{\mathcal{B}^N} \psi\, d\mu =0$.
It is desirable to have an example satisfying $\rdim\left(\mathcal{B}^N, T, \mathbf{d}, \mu\right) = \int_{\mathcal{B}^N} \psi\, d\mu>0$}
of $\mu\in \mathscr{M}^T(\mathcal{B}^N)$, like the above random function 
(\ref{eq: random sum}), satisfying $\rdim\left(\mathcal{B}^N, T, \mathbf{d}, \mu\right) = \int_{\mathcal{B}^N} \psi\, d\mu$?
\end{problem}

\noindent
\textbf{Organization of the paper and how to read it.}
We explain main ideas of the proofs of Theorems \ref{theorem: Ruelle inequality for Brody curves} and 
\ref{theorem: construction of extremal measure} in \S \ref{section: idea of the proofs}.
We review basics of mutual information and rate distortion function in \S \ref{section: mutual information and rate distortion function}.
We explain the definitions of topological and metric mean dimensions with potential and a variational principle for them in 
\S \ref{section: topological and metric mean dimensions with potential}.
We also present a convenient method to calculate metric mean dimension with potential in
\S \ref{section: topological and metric mean dimensions with potential}.
We review the deformation theory of Brody curves and prove Theorem \ref{theorem: Ruelle inequality for Brody curves}
in \S \ref{section: proof of the Ruelle inequality for Brody curves}.
We prepare some results on the energy integral $\int |df|^2 dxdy$ of Brody curves $f$ in 
\S \ref{section: energy integral of Brody curves}.
We develop a general method to construct invariant probability measures with a lower bound on rate distortion dimension
in \S \ref{section: general method to bound the rate distortion dimension from below}.
By combining all the methods prepared in \S\S 
\ref{section: proof of the Ruelle inequality for Brody curves}-\ref{section: general method to bound the rate distortion dimension from below},
we prove Theorem \ref{theorem: construction of extremal measure}
in \S \ref{section: construction of extremal measures}.

Unfortunately this paper is rather long.
About half of this paper 
(mainly \S\S \ref{section: energy integral of Brody curves}-\ref{section: construction of extremal measures})
is devoted to the proof of Theorem \ref{theorem: construction of extremal measure}.
It may be reasonable to concentrate on Theorem \ref{theorem: Ruelle inequality for Brody curves}
(a Ruelle inequality for Brody curves)
at the first reading. Hopefully it will not take too much effort for readers to understand the ideas of 
the proof of Theorem \ref{theorem: Ruelle inequality for Brody curves}.

\S \ref{section: general method to bound the rate distortion dimension from below} is probably 
the heaviest part of the proof of Theorem \ref{theorem: construction of extremal measure}.
It uses several technical facts about measure theory and mutual information prepared in 
\S \ref{section: mutual information and rate distortion function}.
We recommend readers to skip technical details of \S \ref{section: mutual information and rate distortion function} and 
\S \ref{section: general method to bound the rate distortion dimension from below} at the first reading
if they look too cumbersome.

\section{Main idea of the proofs: variational principle} \label{section: idea of the proofs}

We explain a main idea of the proofs of 
Theorems \ref{theorem: Ruelle inequality for Brody curves} and \ref{theorem: construction of extremal measure} in this section.
Our proofs are based on a \textit{variational principle} for mean dimension with potential.
This approach is motivated by the theory of Axiom A diffeomorphisms. 
So we first review it.

\subsection{The case of Axiom A diffeomorphisms}  \label{subsection: the case of Axiom A diffeomorphisms}

Here we briefly review the \textit{thermodynamic formalism} for Axiom A diffeomorphisms.
Readers may skip this subsection if they are totally unfamiliar to this theory.
 
Let $M$ be a compact Riemannian manifold and $T\colon M\to M$ an Axiom A diffeomorphism.
Let $\Omega$ be a basic set of $T$.
For any continuous function $\varphi\colon \Omega \to \mathbb{R}$, we have a quantity $P_T(\varphi)$ called 
the \textbf{topological pressure}.
The variational principle connects it to the Kolomogov--Sinai entropy by 
\begin{equation} \label{eq: variational principle for topological pressure} 
    P_T(\varphi) = \sup_{\mu\in \mathscr{M}^T(\Omega)} \left(h_\mu(T) + \int_\Omega \varphi \, d \mu\right), 
\end{equation}    
where $\mathscr{M}^T(\Omega)$ is the set of all invariant Borel probability measures on $\Omega$
\cite{Walters, Walters_book, Bowen_lecture}.
A measure attaining the supremum of the variational principle (\ref{eq: variational principle for topological pressure}) 
is called an \textbf{equilibrium state for $\varphi$}.
Equilibrium states always exist over a basic set $\Omega$ of an Axiom A diffeomorphism $T$ because 
$T$ is expansive on $\Omega$. (The Kolmogorov--Sinai entropy $h_\mu(T)$ is upper semi-continuous with respect to $\mu$
if the given map $T$ is expansive \cite[Proposition 2.19]{Bowen_lecture}.)

Let $\phi\colon \Omega \to \mathbb{R}$ be the logarithm of the Jacobian in the unstable direction 
(\ref{eq: unstable Jacobian}):
\[
    \phi(x)  = \log \left|\det\left(dT_x\colon E^u_x\to E^u_{Tx}\right)\right|. 
\]
We consider the topological pressure for the function $\varphi(x) := -\phi(x)$ (the minus of $\phi$). 
Then we have \cite[Proposition 4.8 (a)]{Bowen_lecture}
\[   P_T(-\phi) \leq  0. \]
From the variational principle (\ref{eq: variational principle for topological pressure}), this implies that
\[  h_\mu(T) \leq  \int_\Omega \phi \, d\mu \quad \text{for every } \mu\in \mathscr{M}^T(\Omega). \]
This is the Ruelle inequality (\ref{eq: Ruelle inequality}) for a basic set of Axiom A diffeomorphisms.
Moreover, if $\Omega$ is an attractor, then we have \cite[Theorem 4.11]{Bowen_lecture}
\begin{equation}  \label{eq: Topological pressure is zero}
     P_T(-\phi) = 0 .
\end{equation}    
Again, by the variational principle (\ref{eq: variational principle for topological pressure}), this implies 
\[    \sup_{\mu\in \mathscr{M}^T(\Omega)} \left(h_\mu(T) - \int_\Omega \phi \, d \mu\right) =0. \]
As we already mentioned, since $T$ is expansive on $\Omega$, 
this supremum is attained by some measure $\mu^+ \in \mathscr{M}^T(\Omega)$ 
and we have 
\[  h_{\mu^+}(T) = \int_\Omega \phi\, d\mu^+. \]
For the function $-\phi(x)$, it is known that the equilibrium state $\mu^+$ is unique \cite{Bowen_unique}.
It is called the Sinai--Ruelle--Bowen measure.

\subsection{The case of Brody curves}

Next we consider the case of Brody curves.
We explain a rough outline of the proofs of Theorems 
\ref{theorem: Ruelle inequality for Brody curves} and \ref{theorem: construction of extremal measure}.
The proofs are based on a variational principle for mean dimension.
This theory was first introduced by Lindenstrauss and the author in the papers
\cite{Lindenstrauss--Tsukamoto_rate_distortion, Lindenstrauss--Tsukamoto_double_VP} and further expanded in
\cite{Tsukamoto_potential, Tsukamoto_R^d}.
The general framework of the theory will be reviewed in \S \ref{subsection: definitions and variational principle}.
Here we present results suitable for the application to Brody curves.

Let $\mathcal{B}^N$ be the space of Brody curves in $\mathbb{C}P^N$ with
the natural action $T\colon \mathbb{C}\times \mathcal{B}^N  \to  \mathcal{B}^N$.
Let $\mathbf{d}$ be the metric on $\mathcal{B}^N$ defined in (\ref{eq: metric on B^N}):
$\mathbf{d}(f, g) = \sup_{z\in [0,1]^2} d_{\mathrm{FS}}\left(f(z), g(z)\right)$.
For any continuous function $\varphi\colon \mathcal{B}^N\to \mathbb{R}$, we have 
two quantities 
\[ \mdim\left(\mathcal{B}^N, T, \varphi\right) \quad  \text{and} \quad
    \umdimm\left(\mathcal{B}^N, T, \mathbf{d}, \varphi\right) \]
called \textit{mean dimension with potential} and \textit{upper metric mean dimension with potential} respectively.
(We also have a quantity $\lmdimm\left(\mathcal{B}^N, T, \mathbf{d}, \varphi\right)$
called \textit{lower metric mean dimension with potential}. But we do not need it here.)
Mean dimension with potential is a “fusion” of topological dimension and topological pressure.
Upper metric mean dimension with potential is a “fusion” of upper Minkowski dimension and topological pressure.
There quantities are related to rate distortion dimension via the following variational principle.
\begin{equation} \label{eq: lower bound in variational principle in the case of Brody curves}
   \mdim\left(\mathcal{B}^N, T, \varphi\right) \leq 
   \sup_{\mu\in \mathscr{M}^T\left(\mathcal{B}^N\right)} 
   \left(\lrdim\left(\mathcal{B}^N, T, \mathbf{d}, \mu\right) + \int_{\mathcal{B}^N} \varphi\, d\mu\right), 
\end{equation}
\begin{equation} \label{eq: upper bound in variational principle in the case of Brody curves}
  \sup_{\mu\in \mathscr{M}^T\left(\mathcal{B}^N\right)} 
  \left(\urdim\left(\mathcal{B}^N, T, \mathbf{d}, \mu\right) + \int_{\mathcal{B}^N} \varphi\, d\mu\right)
  \leq  \umdimm\left(\mathcal{B}^N, T, \mathbf{d}, \varphi\right).
\end{equation}

Let $\psi\colon \mathcal{B}^N\to \mathbb{R}$ be the continuous function defined in 
(\ref{eq: potential function for Brody curves}): 
\[  \psi(f) = 2(N+1) \,  |df|^2(0).    \]
We consider mean dimension with potential and upper metric mean dimension with potential 
with respect to the potential function $\varphi  := - \psi$ (the minus of $\psi$).
The crucial equation is that
\begin{equation}  \label{eq: mean dimension with potential is zero}
   \mdim\left(\mathcal{B}^N, T, -\psi\right) = \umdimm\left(\mathcal{B}^N, T, \mathbf{d}, -\psi\right) = 0. 
\end{equation}
The proof is based on the deformation theory of Brody curves.
The equation (\ref{eq: mean dimension with potential is zero}) is analogous to the equation 
$P_T(-\phi)=0$ for the Axiom A attractors (\ref{eq: Topological pressure is zero}).
Then all the four quantities in (\ref{eq: lower bound in variational principle in the case of Brody curves}) and 
(\ref{eq: upper bound in variational principle in the case of Brody curves}) are equal to zero for $\varphi = -\psi$.
In particular
\[ \sup_{\mu\in \mathscr{M}^T\left(\mathcal{B}^N\right)} 
  \left(\urdim\left(\mathcal{B}^N, T, \mathbf{d}, \mu\right) - \int_{\mathcal{B}^N} \psi \, d\mu\right) =0. \]
This implies Theorem \ref{theorem: Ruelle inequality for Brody curves} (a Ruelle inequality for Brody curves):
\[   \urdim\left(\mathcal{B}^N, T, \mathbf{d}, \mu\right) \leq \int_{\mathcal{B}^N} \psi\, d\mu \quad 
      \text{for all $\mu \in \mathscr{M}^T\left(\mathcal{B}^N\right)$}.  \]  
Therefore the key ingredients of the proof of Theorem \ref{theorem: Ruelle inequality for Brody curves}   
are the variational principle (\ref{eq: upper bound in variational principle in the case of Brody curves}) and 
the equation (\ref{eq: mean dimension with potential is zero}).
      
Theorem \ref{theorem: construction of extremal measure} does not directly follow from
the variational principle itself.
We need to take a closer look at the machinery behind it.
The proof of the variational principle (\ref{eq: lower bound in variational principle in the case of Brody curves})
presents a general method to construct invariant probability measures with a good lower bound on rate distortion dimension.      
By combining this method with the deformation theory of Brody curves, we can construct 
many invariant measures $\mu \in \mathscr{M}^T\left(\mathcal{B}^N\right)$ satisfying 
\[  \rdim\left(\mathcal{B}^N, T, \mathbf{d}, \mu\right) = \int_{\mathcal{B}^N} \psi\, d\mu. \]
This provides Theorem \ref{theorem: construction of extremal measure}.

\section{Mutual information and rate distortion function}   \label{section: mutual information and rate distortion function}

The purpose of this section is to prepare basics of mutual information and rate distortion function.
Here we omit the proofs.
Readers can find more details in \cite[\S 2]{Tsukamoto_R^d}.
Throughout this paper, $\log x$ denotes the logarithm of base 2:
\[   \log x = \log_2 x \quad (x>0). \]

\subsection{Measure theoretic preparations}  \label{subsection: measure theoretic preparations}

Here we prepare several technical notions in measure theory.
They will be used in \S \ref{section: general method to bound the rate distortion dimension from below}.
We do not need them in the proof of Theorem \ref{theorem: Ruelle inequality for Brody curves}.
Readers may skip this subsection at the first reading if it looks cumbersome. 

A pair $(\mathcal{X}, \mathcal{A})$ is called a \textbf{measurable space} if $\mathcal{X}$ is a set and 
$\mathcal{A}$ is its $\sigma$-algebra.
A triplet $(\mathcal{X}, \mathcal{A}, \mathbb{P})$ is called a \textbf{probability space} if 
$(\mathcal{X}, \mathcal{A})$ is a measurable space and $\mathbb{P}$ is a probability measure on it.

Let $(\mathcal{X}, \mathcal{A})$ and $(\mathcal{Y}, \mathcal{B})$ be measurable spaces.
A \textbf{transition probability} on $\mathcal{X}  \times \mathcal{Y}$ is a map 
$\nu\colon \mathcal{X} \times \mathcal{B} \to [0,1]$ such that 
\begin{itemize}
   \item  for every $x\in \mathcal{X}$, the map $\mathcal{B}  \ni  B \mapsto \nu(x, B) \in [0,1]$ is a probability measure on 
             $(\mathcal{Y}, \mathcal{B})$,
   \item  for every $B  \in  \mathcal{B}$, the map $\mathcal{X} \ni  x \mapsto  \nu(x, B)\in [0,1]$ is measurable.
\end{itemize}
We often denote $\nu(x, B)$ by $\nu(B|x)$.

For a topological space $\mathcal{X}$, its \textbf{Borel $\sigma$-algebra} is the minimum $\sigma$-algebra containing 
all open subsets of $\mathcal{X}$.
We always assume that a topological space is equipped with its Borel $\sigma$-algebra.
We also assume that a finite set is always equipped with the discrete topology and
the discrete $\sigma$-algebra (the set of all subsets).
A topological space $\mathcal{X}$ is called a \textbf{Polish space} if it admits a metric $\mathbf{d}$ for which 
$(\mathcal{X}, \mathbf{d})$ is a complete metric space.

A measurable space $(\mathcal{X}, \mathcal{A})$ is called a \textbf{standard Borel space} if it is isomorphic (as a measurable space)
to some Polish space with its Borel $\sigma$-algebra. 
An importance of this notion comes from Theorem \ref{theorem: regular conditional distribution} below.
Readers can find a detailed explanation of standard Borel spaces in the book of Srivastava \cite{Srivastava}

For two measurable spaces $(\mathcal{X}, \mathcal{A})$ and $(\mathcal{Y}, \mathcal{B})$ we denote their product 
by $(\mathcal{X} \times \mathcal{Y}, \mathcal{A} \otimes \mathcal{B})$ where 
$\mathcal{A} \otimes \mathcal{B}$ is the minimum $\sigma$-algebra containing all  
$A \times B$ $(A \in \mathcal{A}, B \in \mathcal{B})$.
For any $E\in \mathcal{A} \otimes \mathcal{B}$, the section $E_x := \{y\in \mathcal{Y}  \mid (x, y) \in E\}$
belongs to $\mathcal{B}$ for every $x\in \mathcal{X}$.
Moreover, if $(\mathcal{Y}, \mathcal{B})$ is a standard Borel space, then for any transition probability 
$\nu$ on $\mathcal{X} \times \mathcal{Y}$ and any $E\in \mathcal{A} \otimes \mathcal{B}$, the map
$\mathcal{X} \ni x\mapsto  \nu(E_x|x)\in [0,1]$ is measurable \cite[Proposition 3.4.24]{Srivastava}. 

Let $(\Omega, \mathcal{F}, \mathbb{P})$ be a probability space and $(\mathcal{X}, \mathcal{A})$ a measurable space.
Let $X\colon \Omega \to \mathcal{X}$ be a measurable map.
We denote the push-forward measure $X_*\mathbb{P}$ by $\mathrm{Law} X$, where 
$X_*\mathbb{P}(A) = \mathbb{P}\left(X^{-1}(A)\right)$.
We call it the \textbf{law of $X$} or the \textbf{distribution of $X$}.

The next theorem guarantees the existence of regular conditional distribution.
For the proof, see \cite[p.15 Theorem 3.3 and its Corollary]{Ikeda--Watanabe} or \cite[p.182 Corollary 6.2]{Gray_entropy}.

\begin{theorem}  \label{theorem: regular conditional distribution}
Let $(\Omega, \mathcal{F}, \mathbb{P})$ be a probability space. Let $(\mathcal{X}, \mathcal{A})$ and 
$(\mathcal{Y}, \mathcal{B})$ be standard Borel spaces, and let $X \colon \Omega \to \mathcal{X}$ and 
$Y \colon \Omega \to \mathcal{Y}$ be measurable maps.
Set $\mu = \mathrm{Law} X$.
Then there exists a transition probability $\nu$ on $\mathcal{X} \times \mathcal{Y}$ such that for any 
$E\in \mathcal{A} \otimes \mathcal{B}$ we have 
\[  \mathbb{P}\left((X, Y)\in E \right) = \int_{\mathcal{X}} \nu (E_x|x) \, d\mu(x). \]
If another transition probability $\nu^\prime$ on $\mathcal{X} \times \mathcal{Y}$ satisfies the same property then 
there exists a $\mu$-null set $N\in \mathcal{A}$ such that 
$\nu(B|x) = \nu^\prime(B|x)$ for all $x\in \mathcal{X} \setminus  N$ and $B \in \mathcal{B}$.
\end{theorem}

The transition probability $\nu$ introduced in this theorem is called the \textbf{regular conditional distribution of $Y$
given $X=x$}. We sometimes denote $\nu(B|x)$ by $\mathbb{P}\left(Y \in B| X=x\right)$.
If $\mathcal{X}$ and $\mathcal{Y}$ are finite sets, then this coincides with the notion of conditional probability in 
the elementary probability theory:
\[  \mathbb{P}\left(Y \in B| X = x\right) = \frac{\mathbb{P}\left(X=x, Y\in  B \right)}{\mathbb{P}(X=x)}, \quad 
     \text{if $\mathbb{P}(X=x) >0$}. \]
In this case we denote $\nu\left(\{y\}|x\right)$ by $\nu(y|x)$ for $x \in \mathcal{X}$ and $y \in \mathcal{Y}$
and call it a \textbf{conditional probability mass function}.

Let $(\Omega, \mathcal{F}, \mathbb{P})$ be a probability space. 
Let $(\mathcal{X}, \mathcal{A}), (\mathcal{Y}, \mathcal{B}), (\mathcal{Z}, \mathcal{C})$
be standard Borel spaces, and let $X, Y, Z$ be random variables defined on $\Omega$ 
and taking values in $\mathcal{X}, \mathcal{Y}, \mathcal{Z}$
respectively.
We say that $X$ and $Y$ are \textbf{conditionally independent given $Z$} if we have 
\[  \mathbb{P}\left((X, Y)\in A \times  B| Z=z \right) = 
    \mathbb{P}\left(X \in A| Z=z\right) \cdot \mathbb{P}\left(Y\in B|Z=z\right) \]
for all $A \in \mathcal{A}, B \in \mathcal{B}$ and for ${Z}_*\mathbb{P}$-almost every $z \in \mathcal{Z}$.
Here the left-hand side is the regular conditional distribution of $(X, Y)\colon \Omega \to \mathcal{X} \times \mathcal{Y}$
given $Z=z$.

\subsection{Mutual information}
The purpose of this subsection is to define mutual information and review its basic properties.
Readers can find more comprehensive introduction to mutual information in the book of Cover--Thomas \cite{Cover--Thomas}.
The book of Gray \cite{Gray_entropy} presents a mathematically sophisticated explanation.

We fix a probability space $(\Omega, \mathcal{F}, \mathbb{P})$.
We assume that all random variables in this subsection are defined on $(\Omega, \mathcal{F},\mathbb{P})$.

Let $(\mathcal{X}, \mathcal{A})$ be a finite set with the discrete $\sigma$-algebra.
For a measurable map $X\colon \Omega \to \mathcal{X}$ we define its \textbf{Shannon entropy} by 
\[  H(X) = - \sum_{x\in \mathcal{X}} \mathbb{P}(X=x) \log \mathbb{P}(X=x). \]
Here we assume $0 \log 0 =0$.

Let $(\mathcal{X}, \mathcal{A})$ and $(\mathcal{Y}, \mathcal{B})$ be measurable spaces. 
Let $X \colon \Omega \to \mathcal{X}$ and $Y\colon  \Omega \to \mathcal{Y}$ be measurable maps.
We define the \textbf{mutual information} $I(X; Y)$ in the following two steps.
Intuitively $I(X; Y)$ is the amount of information shared by $X$ and $Y$.

\begin{itemize}
  \item  \textbf{Step I:} Suppose $\mathcal{X}$ and $\mathcal{Y}$ are both finite sets.
            Then we define 
            \[  I(X;Y) = H(X) + H(Y) - H(X, Y),  \]
            where $H(X, Y)$ is the Shannon entropy of $(X, Y)\colon \Omega \to \mathcal{X} \times \mathcal{Y}$.
            $I(X;Y)$ is always nonnegative.
            The explicit formula is given by 
            \[  I(X;Y) = \sum_{x \in \mathcal{X}, y \in \mathcal{Y}}
            \mathbb{P}\left(X=x, Y=y\right) 
            \log \frac{\mathbb{P}\left(X=x, Y=y\right)}{\mathbb{P}(X=x) \mathbb{P}(Y=y)}. \]
            Here we assume $0 \log \frac{0}{a} = 0$ for any $a\geq 0$.
   \item  \textbf{Step II:} Next we define $I(X;Y)$ in general (i.e. $\mathcal{X}$ and $\mathcal{Y}$ may be infinite sets).
            Let $\mathcal{X}^\prime$ and $\mathcal{Y}^\prime$ be arbitrary finite sets, 
            and let $f \colon \mathcal{X} \to \mathcal{X}^\prime$ and $g\colon \mathcal{Y}\to \mathcal{Y}^\prime$
            be any measurable maps. Now we can consider $I\left(f(X); g(Y)\right)$ by the above Step I.
            We define $I(X;Y)$ as the supremum of $I\left(f(X); g(Y)\right)$ over all finite sets $\mathcal{X}^\prime, \mathcal{Y}^\prime$ and 
            all measurable maps $f \colon \mathcal{X} \to \mathcal{X}^\prime$ and $g \colon \mathcal{Y}\to \mathcal{Y}^\prime$.
            This definition is compatible with Step I.
            The mutual information is nonnegative and symmetric: 
            \[  I(X; Y) = I(Y; X) \geq 0. \]
            If $f \colon \mathcal{X} \to \mathcal{X}^\prime$ and $g\colon \mathcal{Y} \to \mathcal{Y}^\prime$ 
            are measurable maps to some measurable spaces
            $(\mathcal{X}^\prime, \mathcal{A}^\prime)$ and $(\mathcal{Y}^\prime, \mathcal{B}^\prime)$ 
            (not necessarily finite sets) then we have 
            $I\left(f(X); g(Y)\right) \leq  I(X;Y)$.
\end{itemize}

If $(\mathcal{X}, \mathcal{A})$ and $(\mathcal{Y}, \mathcal{B})$ are both standard Borel spaces, then we can consider the 
regular conditional distribution of $Y$ given $X=x$.  Set
\[  \nu(B|x) = \mathbb{P}(Y \in B| X=x), \quad (x\in \mathcal{X}, B \in \mathcal{B}). \]
We also set $\mu = \mathrm{Law} X$.
The distribution of $(X, Y)$ is determined by $\mu$ and $\nu$.
Therefore the mutual information $I(X;Y)$ is also determined by them.
So we sometimes denote $I(X;Y)$ by $I(\mu, \nu)$.
An important point is that $I(\mu, \nu)$ is a concave function in $\mu$ and a convex function in $\nu$.
See Proposition \ref{prop: concavity and convexity of mutual information} below.

In the rest of this subsection we review basic properties of mutual information without proofs.
Readers can find the proofs  in \cite[\S 2.2]{Tsukamoto_R^d}.
We need to use these properties of mutual information in \S \ref{section: general method to bound the rate distortion dimension from below}.
We do not need them in the proof of Theorem \ref{theorem: Ruelle inequality for Brody curves}.

The next lemma immediately follows from Step I of the above definition.

\begin{lemma} \label{lemma: convergence in law and mutual information}   
  Let $X_n$ and $Y_n$ $(n\geq 1)$ be sequences of random variables taking values in finite sets $\mathcal{X}$ and $\mathcal{Y}$
  respectively.  Suppose $(X_n, Y_n)$ converges to $(X, Y)$ in law\footnote{“Convergence in law” means that 
  $\mathbb{P}(X_n=x, Y_n=y) \to \mathbb{P}(X=x, Y=y)$ as $n\to \infty$ for all $x\in \mathcal{X}$ and $y\in \mathcal{Y}$.}.
  Then $I(X_n;Y_n)$ converges to $I(X;Y)$ as $n\to \infty$.
\end{lemma}

\begin{lemma}[Subadditivity of mutual information]  \label{lemma: subadditivity of mutual information} 
Let $X, Y, Z$ be three random variables taking values in standard Borel spaces 
$(\mathcal{X}, \mathcal{A}), (\mathcal{Y}, \mathcal{B}), (\mathcal{Z}, \mathcal{C})$ respectively. 
Suppose that $X$ and $Y$ are conditionally independent given $Z$. Then 
\[   I\left(X, Y; Z\right) \leq I(X; Z) + I(Y;Z). \]
Here $I\left(X, Y; Z\right) = I\left((X, Y); Z\right)$ is the mutual information between $(X, Y)$ and $Z$.
\end{lemma}

\begin{proposition}[$I(\mu, \nu)$ is concave in $\mu$ and convex in $\nu$] \label{prop: concavity and convexity of mutual information}
Let $(\mathcal{X}, \mathcal{A})$ and $(\mathcal{Y}, \mathcal{B})$ be standard Borel spaces, and let 
$(\mathcal{Z}, \mathcal{C}, m)$ be a probability space.
    \begin{enumerate}
      \item  Let $\nu$ be a transition probability on $\mathcal{X} \times \mathcal{Y}$.
                Suppose that we are given a probability measure $\mu_z$ on $\mathcal{X}$ for each $z \in \mathcal{Z}$
                such that the map $\mathcal{Z}\ni z \mapsto \mu_z(A)\in [0,1]$ is measurable for every $A\in \mathcal{A}$.
                We define a probability measure $\mu$ on $(\mathcal{X}, \mathcal{A})$ by 
                \[  \mu(A) = \int_{\mathcal{Z}} \mu_z(A)\, dm(z), \quad (A\in \mathcal{A}). \]
                Then we have 
                \[  I(\mu, \nu) \geq  \int_{\mathcal{Z}} I\left(\mu_z, \nu\right)\, dm(z). \]
      \item  Let $\mu$ be a probability measure on $\mathcal{X}$. Suppose that we are given a transition probability 
                $\nu_z$ on $\mathcal{X}\times \mathcal{Y}$ for each $z\in \mathcal{Z}$ such that the map 
                $\mathcal{X}\times \mathcal{Z}\ni (x, z) \mapsto \nu_z(B|x)\in [0,1]$ is measurable with respect to
                $\mathcal{A}\otimes \mathcal{C}$ for each $B\in \mathcal{B}$.
                We define a transition probability $\nu$ on $\mathcal{X}\times \mathcal{Y}$ by 
                \[  \nu(B|x) = \int_{\mathcal{Z}} \nu_z(B|x) \, dm(z), \quad (x\in \mathcal{X}, B\in \mathcal{B}). \]
                Then we have           
                \[  I(\mu, \nu) \leq  \int_{\mathcal{Z}} I(\mu, \nu_z)\, dm(z). \]
    \end{enumerate}
\end{proposition}

The next proposition is a key result which connects geometry to rate distortion theory.
This is essentially due to Kawabata--Dembo \cite[Proposition 3.2]{Kawabata--Dembo}.
For a subset $E$ of a metric space $(\mathcal{X}, \mathbf{d})$ we denote 
$\diam E = \sup\{\mathbf{d}(x, y)\mid x, y\in E\}$.

\begin{proposition}[Kawabata--Dembo estimate]  \label{prop: Kawabata--Dembo estimate}
Let $\varepsilon$ and $\delta$ be positive numbers with $2\varepsilon \log (1/\varepsilon) \leq \delta$.
Let $s$ be a nonnegative real number.
Let $(\mathcal{X}, \mathbf{d})$ be a separable metric space with a Borel probability measure $\mu$ satisfying 
\[  \mu(E) \leq \left(\diam E\right)^s  \quad \text{for all Borel subsets $E\subset \mathcal{X}$ with 
     $\diam E < \delta$}. \]
Let $X$ and $Y$ be random variables taking values in $\mathcal{X}$ and satisfying $\mathrm{Law} X = \mu$ and 
$\mathbb{E}\mathbf{d}(X, Y) < \varepsilon$. Then 
\[   I(X;Y) \geq  s\log (1/\varepsilon) - K(s+1). \]
Here $K$ is a universal positive constant independent of the given data $\varepsilon, \delta, s, (\mathcal{X}, \mathbf{d}), \mu$.
\end{proposition}

\subsection{Rate distortion theory}  \label{subsection: rate distortion theory}
The purpose of this subsection is to define rate distortion function.
This was first introduced by Shannon \cite{Shannon, Shannon59}.
Readers can find a nice introduction in the book of Cover--Thomas \cite[Chapter 10]{Cover--Thomas}.
Rate distortion theory for continuous-time stochastic processes was investigated by 
Pursley--Gray \cite{Pursley--Gray}.

Let $k$ be a natural number.
We denote the Lebesgue measure on $\mathbb{R}^k$ by $\mathbf{m}$.
Let $(\mathcal{X}, \mathbf{d})$ be a compact metric space.
Let $A$ be a Borel subset of $\mathbb{R}^k$ with $\mathbf{m}(A) < \infty$.
We define $L^1(A, \mathcal{X})$ as the space of all measurable maps 
$f\colon A\to \mathcal{X}$. We identify two maps if they coincide $\mathbf{m}$-almost everywhere.
Define a metric on $L^1(A, \mathcal{X})$ by 
\[  D(f, g) = \int_A\mathbf{d}\left(f(u), g(u)\right) \, d\mathbf{m}(u), \quad 
     (f, g \in L^1(A, \mathcal{X})). \]
$\left(L^1(A, \mathcal{X}), D\right)$ is a complete separable metric space \cite[Lemma 2.14]{Tsukamoto_R^d}.
Hence it is a standard Borel space with respect to the Borel $\sigma$-algebra.

Let $T\colon \mathbb{R}^k\times \mathcal{X}\to \mathcal{X}$ be a continuous action of $\mathbb{R}^k$.
Here $\mathbb{R}^k$ has the standard topology and additive group structure.
(We only need the case of $k=2$ in the application to Brody curves.)
Let $\mu$ be a $T$-invariant Borel probability measure on $\mathcal{X}$.
Here “$T$-invariant” means that $\mu(T^u A) = \mu(A)$ for all $u\in \mathbb{R}^k$ and all Borel subsets $A\subset \mathcal{X}$.

Let $\varepsilon>0$. Let $A$ be a bounded Borel subset of $\mathbb{R}^k$ with $\mathbf{m}(A)>0$.
We define $R(\varepsilon, A)$ as the infimum of the mutual information $I(X;Y)$ where $X$ and $Y$ are random 
variables defined on some probability space $(\Omega, \mathcal{F}, \mathbb{P})$ such that 
 \begin{itemize}
    \item  $X$ takes values in $\mathcal{X}$ and its distribution is given by $\mu$,
    \item  $Y$ takes values in $L^1(A, \mathcal{X})$ and satisfies 
    \[ \mathbb{E}\left(\frac{1}{\mathbf{m}(A)} \int_A \mathbf{d}\left(T^u X, Y_u\right)\, d\mathbf{m}(u) \right) < \varepsilon. \]
     Here $Y_u = Y_u(\omega)$ $(\omega\in \Omega)$ is the value of $Y(\omega) \in L^1(A, \mathcal{X})$ at $u\in A$.
 \end{itemize}

We define the \textbf{rate distortion function} $R(\mathbf{d}, \mu, \varepsilon)$ by 
\[  R(\mathbf{d}, \mu, \varepsilon) = \lim_{L\to \infty} \frac{R\left(\varepsilon, [0, L)^k\right)}{L^k}. \]
This limit exists, and it is equal to the infimum of $\frac{R\left(\varepsilon, [0, L)^k\right)}{L^k}$ over $L>0$
\cite[Lemma 2.17]{Tsukamoto_R^d}.
The value of $R(\mathbf{d}, \mu, \varepsilon)$ is a nonnegative real number.

We define the \textbf{upper and lower rate distortion dimensions} of 
$(\mathcal{X}, T, \mathbf{d}, \mu)$ by 
\[  \urdim\left(\mathcal{X}, T, \mathbf{d}, \mu\right) = \varlimsup_{\varepsilon \to 0} 
     \frac{R(\mathbf{d}, \mu, \varepsilon)}{\log (1/\varepsilon)}, \quad 
     \lrdim\left(\mathcal{X}, T, \mathbf{d}, \mu\right) = \varliminf_{\varepsilon \to 0} 
     \frac{R(\mathbf{d}, \mu, \varepsilon)}{\log (1/\varepsilon)}.  \]
These are nonnegative real numbers (possibly $+\infty$).
When the upper and lower limits coincide, we denote the common value by $\rdim(\mathcal{X}, T, \mathbf{d}, \mu)$.

The next remark will be used in \S \ref{subsection: constructing invariant measures with a lower bound on rate distortion dimension}.
We do not need it for the proof of Theorem \ref{theorem: Ruelle inequality for Brody curves}.

\begin{remark} \label{remark: finite random variable in the definition of rate distortion function}
In the above definition of $R(\varepsilon, A)$, 
we can assume that $Y$ takes only finitely many values in $L^1(A, \mathcal{X})$ without loss of generality.
Indeed let $X$ and $Y$ be random variables satisfying the conditions in the definition of $R(\varepsilon, A)$.
A priori, $Y$ may take infinitely many values in general.
However we can always find another random variable $Z$ taking values in $L^1(A, \mathcal{X})$ such that 
$Z$ takes only finitely many values and satisfies 
\[  \mathbb{E}\left(\frac{1}{\mathbf{m}(A)} \int_A \mathbf{d}\left(T^u X, Z_u\right)\, d\mathbf{m}(u) \right) < \varepsilon, \quad
     I(X;Z) \leq  I(X;Y). \]
Therefore, when we take the infimum, we only need to consider such random variables $Z$.     
The construction is as follows.
Take a positive number $\tau$ with 
\[ \mathbb{E}\left(\frac{1}{\mathbf{m}(A)} \int_A \mathbf{d}\left(T^u X, Y_u\right)\, d\mathbf{m}(u) \right) < \varepsilon -2\tau. \]
The space $L^1(A, \mathcal{X})$ is separable. Hence we can find a countable dense subset 
$\{f_1, f_2, f_3, \dots\}$. 
We define a map $F\colon L^1(A, \mathcal{X}) \to \{f_1, f_2, f_3, \dots\}$ by 
$F(f) = f_n$ where $n$ is the smallest natural number satisfying $D(f, f_n) < \tau \cdot \mathbf{m}(A)$.
Set $Y^\prime = F(Y)$. Then $D(Y, Y^\prime) < \tau \cdot \mathbf{m}(A)$ and hence
\[   \mathbb{E}\left(\frac{1}{\mathbf{m}(A)} \int_A \mathbf{d}\left(T^u X, Y^\prime_u\right)\, d\mathbf{m}(u) \right) < \varepsilon -\tau. \]
We choose a natural number $n_0$ satisfying 
\[  \left(\sum_{n>n_0} \mathbb{P}(Y^\prime = f_n)\right) \cdot \diam \mathcal{X} < \tau. \]
Define $G\colon \{f_1, f_2, f_3, \dots\}\to \{f_1, f_2, \dots, f_{n_0}\}$ by 
\[  G(f)   = \begin{cases}  f & \text{if $f\in \{f_1, \dots, f_{n_0}\}$} \\
                                    f_{n_0} & \text{otherwise} \end{cases}. \] 
Set $Z = G(Y^\prime)$.
Then $Z$ takes only finitely many values $f_1, f_2, \dots, f_{n_0}$ and satisfies 
\[  \mathbb{E}\left(\frac{1}{\mathbf{m}(A)} \int_A \mathbf{d}\left(T^u X, Z_u\right)\, d\mathbf{m}(u) \right) < \varepsilon, \quad
     I(X;Z) \leq I(X; Y^\prime) \leq  I(X;Y). \]
\end{remark}

\section{Topological and metric mean dimensions with potential}  
\label{section: topological and metric mean dimensions with potential}

The purpose of this section is to prepare basics of 
mean dimension with potential and metric mean dimension with potential.

\subsection{Definitions and variational principle}  \label{subsection: definitions and variational principle}

\begin{figure}[h] 
   \centering
   \includegraphics[width=3.0in]{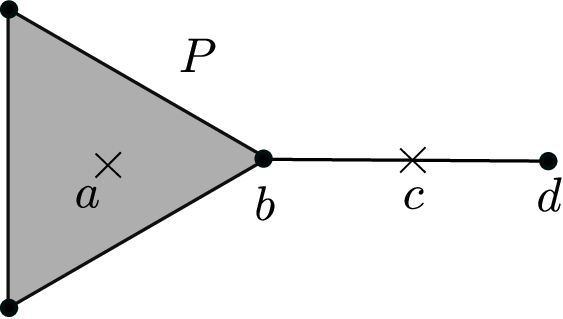}
    \caption{Here $P$ has four vertexes (denoted by dots), four $1$-dimensional simplexes and one $2$-dimensional simplex. 
    The points $b$ and $d$ are vertexes of $P$ whereas $a$ and $c$ are not. 
    We have $\dim_a P = \dim_b P =2$ and $\dim_c P = \dim_d P =1$.}    \label{figure: local dimension}
\end{figure}

Let $P$ be a finite simplicial complex. For a point $a\in P$ we define the \textbf{local dimension} $\dim_a P$ as the 
maximum of $\dim \Delta$ where $\Delta$ is a simplex of $P$ containing $a$.
See Figure \ref{figure: local dimension}. (This is the same as \cite[Fig. 1]{Tsukamoto_potential}.)

Let $(\mathcal{X}, \mathbf{d})$ be a compact metric space.
For a positive number $\varepsilon$, a continuous map $f\colon \mathcal{X}\to \mathcal{Y}$ to some topological space 
$\mathcal{Y}$ is called an \textbf{$\varepsilon$-embedding} if we have $\diam f^{-1}(y) < \varepsilon$ for all $y\in \mathcal{Y}$.
Here $\diam f^{-1}(y) = \sup_{x_1, x_2\in f^{-1}(y)} \mathbf{d}(x_1, x_2)$ is the diameter of the fiber $f^{-1}(y)$.

Let $\varphi\colon \mathcal{X} \to \mathbb{R}$ be a continuous function.
We define \textbf{$\varepsilon$-width dimension with potential} by
\begin{equation*} 
   \begin{split}
     & \widim_\varepsilon(\mathcal{X}, \mathbf{d}, \varphi)  \\ 
     & =  \inf\left\{\max_{x\in \mathcal{X}} \left(\dim_{f(x)} P + \varphi(x)\right) \middle|
    \parbox{3in}{\centering $P$ is a finite simplicial complex and $f\colon \mathcal{X}\to P$ is an $\varepsilon$-embedding}\right\}.  
   \end{split}
\end{equation*}    
We define \textbf{$\varepsilon$-covering number with potential} by
\[ \#\left(\mathcal{X}, \mathbf{d}, \varphi, \varepsilon\right) =  
    \inf\left\{\sum_{i=1}^n (1/\varepsilon)^{\sup_{U_i}\varphi} \middle|\, 
    \parbox{3in}{\centering  $\mathcal{X} = U_1\cup \dots \cup U_n$ is an open cover with $\diam\, U_i < \varepsilon$
     for all $1\leq i\leq n$} \right\}. \]
When $\varphi = 0$, we denote $\#\left(\mathcal{X}, \mathbf{d}, \varphi, \varepsilon\right)$ by 
$\#\left(\mathcal{X}, \mathbf{d}, \varphi\right)$ and call it \textbf{$\varepsilon$-covering number}.
The $\varepsilon$-covering number $\#\left(\mathcal{X}, \mathbf{d}, \varepsilon\right)$
is given by the minimum integer $n\geq 1$ such that $\mathcal{X}$ is covered by $n$ open sets of diameter smaller than $\varepsilon$.

Let $k$ be a natural number. Let $T\colon \mathbb{R}^k \times  \mathcal{X} \to \mathcal{X}$ be a continuous action of 
the group $\mathbb{R}^k$. 
(As in \S \ref{subsection: rate distortion theory}, $\mathbb{R}^k$ has the standard topology and additive group structure.)
For a Borel subset $A$ of $\mathbb{R}^k$ we define a metric $\mathbf{d}_A$ and a function $\varphi_A$ on $\mathcal{X}$ by 
\[  \mathbf{d}_A(x, y) = \sup_{u\in A} \mathbf{d}\left(T^u x, T^u y\right), \quad 
     \varphi_A(x) = \int_A \varphi(T^u x)\, d\mathbf{m}(u). \]
Here $\mathbf{m}$ denotes the $k$-dimensional Lebesgue measure on $\mathbb{R}^k$.
For a positive number $L$, we also denote $\mathbf{d}_L := \mathbf{d}_{[0,L]^k}$ and 
$\varphi_L  := \varphi_{[0, L]^k}$ for simplicity of the notation.

We define \textbf{mean dimension with potential} by 
\[  \mdim\left(\mathcal{X}, T, \varphi\right) 
     = \lim_{\varepsilon\to 0} 
     \left(\lim_{L\to \infty} \frac{\widim_\varepsilon \left(\mathcal{X}, \mathbf{d}_L, \varphi_L \right)}{L^k}\right). \]
This is a topological invariant (i.e. independent of the choice of $\mathbf{d}$).     
When $\varphi=0$, we denote $\mdim\left(\mathcal{X}, T, \varphi\right)$ by 
$\mdim\left(\mathcal{X}, T\right)$ and call it \textbf{mean dimension}.
This coincides with the original definition of mean dimension introduced by 
Gromov \cite{Gromov}.   
     
We define \textbf{upper and lower metric mean dimensions with potential} by 
\[ \umdimm\left(\mathcal{X}, T, \mathbf{d}, \varphi\right) = 
    \limsup_{\varepsilon \to 0} \left(\lim_{L\to \infty} 
    \frac{\log \#\left(\mathcal{X}, \mathbf{d}_L, \varphi_L, \varepsilon\right)}{L^k \log (1/\varepsilon)}\right),   \]
\[  \lmdimm\left(\mathcal{X}, T, \mathbf{d}, \varphi\right) = 
    \liminf_{\varepsilon \to 0} \left(\lim_{L\to \infty} 
    \frac{\log \#\left(\mathcal{X}, \mathbf{d}_L, \varphi_L, \varepsilon\right)}{L^k \log (1/\varepsilon)}\right).   \]
If the upper and lower limits coincide, we denote the common value by $\mdimm\left(\mathcal{X}, T, \mathbf{d}, \varphi\right)$.
When $\varphi=0$, we denote $\umdimm\left(\mathcal{X}, T, \mathbf{d}, \varphi\right)$ and 
$\lmdimm\left(\mathcal{X}, T, \mathbf{d}, \varphi\right)$ by
$\umdimm\left(\mathcal{X}, T, \mathbf{d}\right)$ and $\lmdimm\left(\mathcal{X}, T, \mathbf{d}\right)$
respectively.
They coincide with the original definition of the upper and lower metric mean dimensions 
introduced by Lindenstrauss and Weiss \cite{Lindenstrauss--Weiss, Lindenstrauss}.
We always have \cite[Lemma 3.3, Theorem 3.4]{Tsukamoto_R^d}
\begin{equation} \label{eq: mdim lmdimm umdimm}
    \mdim\left(\mathcal{X}, T, \varphi\right) \leq \lmdimm\left(\mathcal{X}, T, \mathbf{d}, \varphi\right)
     \leq  \umdimm\left(\mathcal{X}, T, \mathbf{d}, \varphi \right). 
\end{equation}

We denote by $\mathscr{M}^T(\mathcal{X})$ the set of all $T$-invariant Borel probability measures on $\mathcal{X}$.
The following variational principle was proved in \cite[Theorem 1.3, Proposition 3.5]{Tsukamoto_R^d}.

\begin{theorem}[Variational principle] \label{theorem: variational principle for mean dimension}
\begin{equation*}
   \begin{split}
     \mdim\left(\mathcal{X}, T, \varphi\right) &\leq \sup_{\mu \in \mathscr{M}^T(\mathcal{X})}
     \left(\lrdim\left(\mathcal{X}, T, \mathbf{d}, \mu\right) + \int_{\mathcal{X}}\varphi\, d\mu\right) \\
    & \leq \sup_{\mu \in \mathscr{M}^T(\mathcal{X})}
     \left(\urdim\left(\mathcal{X}, T, \mathbf{d}, \mu\right) + \int_{\mathcal{X}}\varphi\, d\mu\right)
     \leq  \umdimm\left(\mathcal{X}, T, \mathbf{d}, \varphi\right).
    \end{split} 
\end{equation*}     
\end{theorem}

\begin{remark}  \label{remark: double variational principle conjecture}
The standard variational principle for topological pressure provides the equality 
$P_T(\varphi)  = \sup_{\mu\in \mathscr{M}^T(\mathcal{X})}\left(h_\mu(T) + \int_{\mathcal{X}}\varphi\, d\mu\right)$,
whereas Theorem \ref{theorem: variational principle for mean dimension} presents only inequalities.
This might look unsatisfactory.
Indeed, we conjecture that for any continuous action $T\colon \mathbb{R}^k\times \mathcal{X}\to \mathcal{X}$ on 
a compact metrizable space $\mathcal{X}$ and for any continuous function 
$\varphi\colon \mathcal{X}\to \mathbb{R}$ there exists a metric $\mathbf{d}$ on $\mathcal{X}$ compatible with the given 
topology and satisfying 
\begin{equation*}
   \begin{split}
     \mdim\left(\mathcal{X}, T, \varphi\right) &= \sup_{\mu \in \mathscr{M}^T(\mathcal{X})}
     \left(\lrdim\left(\mathcal{X}, T, \mathbf{d}, \mu\right) + \int_{\mathcal{X}}\varphi\, d\mu\right) \\
    & = \sup_{\mu \in \mathscr{M}^T(\mathcal{X})}
     \left(\urdim\left(\mathcal{X}, T, \mathbf{d}, \mu\right) + \int_{\mathcal{X}}\varphi\, d\mu\right)
     = \mdimm\left(\mathcal{X}, T, \mathbf{d}, \varphi\right).
    \end{split} 
\end{equation*}  
This is more satisfactory (if it is true).   
We will see later that this conjecture holds true for the dynamical system of Brody curves $\mathcal{B}^N$ with respect to
some natural choices of potential functions.
See Remark \ref{remark: ordinary topological and metric mean dimension of Brody curves} below.
\end{remark}

\subsection{Local formula of metric mean dimension with potential} 
\label{subsection: local formula}

Here we prepare a convenient method to calculate metric mean dimension with potential
from certain local quantity.
Let $(\mathcal{X}, \mathbf{d})$ be a compact metric space with a continuous function 
$\varphi\colon \mathcal{X} \to \mathbb{R}$.
For $\varepsilon >0$ and a subset $E \subset \mathcal{X}$ we set 
\[  \#\left(E, \mathbf{d}, \varphi, \varepsilon\right) 
    =  \inf\left\{\sum_{i=1}^n (1/\varepsilon)^{\sup_{U_i}\varphi} \middle|\, 
    \parbox{3in}{\centering  $E\subset U_1\cup \dots \cup U_n$. Each $U_i$ is an open set of $\mathcal{X}$ 
    with $\diam\, U_i < \varepsilon$.} \right\}. \]
We also set $\#\left(E, \mathbf{d}, \varepsilon\right) := \#\left(E, \mathbf{d}, 0, \varepsilon\right)$.
This is the minimum number of $n$ such that there exist open 
sets $U_1, \dots, U_n$ of $\mathcal{X}$ satisfying $E\subset U_1\cup \dots \cup U_n$ and 
$\diam U_i < \varepsilon$ for all $1\leq i \leq n$.

Let $T\colon \mathbb{R}^k\times \mathcal{X}\to \mathcal{X}$ be a continuous action of $\mathbb{R}^k$.
For $\delta>0$ and $x\in \mathcal{X}$ we set 
\[  B_\delta\left(x, \mathbf{d}_{\mathbb{R}^k}\right) := 
    \{y\in \mathcal{X}\mid \mathbf{d}_{\mathbb{R}^k}(x, y) \leq \delta\} 
     = \{y\mid \mathbf{d}\left(T^u x, T^u y\right) \leq \delta \text{ for all $u\in \mathbb{R}^k$}\}. \]
Notice that this is \textit{not} a neighborhood of $x$ (with respect to the original topology of $\mathcal{X}$)
in general.
Nonetheless, the next theorem shows that 
we can calculate metric mean dimension with potential only by studying these sets.
This was proved in \cite[Theorem 7.1]{Tsukamoto_R^d}.

\begin{theorem}  \label{theorem: local formula of metric mean dimension with potential}
For any $\delta>0$ we have 
\[ \umdimm\left(\mathcal{X}, T, \mathbf{d}, \varphi\right) = 
   \varlimsup_{\varepsilon\to 0}\left\{\sup_{x\in \mathcal{X}}\left(\varliminf_{L\to \infty} 
   \frac{\log \#\left(B_\delta\left(x, \mathbf{d}_{\mathbb{R}^k}\right), \mathbf{d}_L, \varphi_L, \varepsilon\right)}{L^k \log (1/\varepsilon)}
   \right)\right\}, \]
\[  \lmdimm\left(\mathcal{X}, T, \mathbf{d}, \varphi\right) = 
   \varliminf_{\varepsilon\to 0}\left\{\sup_{x\in \mathcal{X}}\left(\varliminf_{L\to \infty} 
   \frac{\log \#\left(B_\delta\left(x, \mathbf{d}_{\mathbb{R}^k}\right), \mathbf{d}_L, \varphi_L, \varepsilon\right)}{L^k \log (1/\varepsilon)}
   \right)\right\}.  \]
\end{theorem}

The next proposition is a small modification of this theorem.
This is more suitable for the application to Brody curves.

\begin{proposition} \label{prop: modified local formula}
\[ \umdimm\left(\mathcal{X}, T, \mathbf{d}, \varphi\right) = \lim_{\delta\to 0} \left[\varlimsup_{\varepsilon \to 0} \left\{\sup_{x\in \mathcal{X}}
    \varliminf_{L\to \infty} 
    \left(\frac{\log \#\left(B_\delta\left(x, \mathbf{d}_{\mathbb{R}^k}\right), \mathbf{d}_L, \varepsilon\right)}{L^k \log (1/\varepsilon)}
    + \frac{\varphi_L(x)}{L^k}\right)\right\}\right], \]
\[ \lmdimm\left(\mathcal{X}, T, \mathbf{d}, \varphi\right) = \lim_{\delta\to 0} \left[\varliminf_{\varepsilon \to 0} \left\{\sup_{x\in \mathcal{X}}
    \varliminf_{L\to \infty} 
    \left(\frac{\log \#\left(B_\delta\left(x, \mathbf{d}_{\mathbb{R}^k}\right), \mathbf{d}_L, \varepsilon\right)}{L^k \log (1/\varepsilon)}
    + \frac{\varphi_L(x)}{L^k}\right)\right\}\right], \]
\end{proposition}

This formula may look complicated, but indeed the calculation of the right-hand sides are often much easier than 
the direct calculation of the left-hand sides.

\begin{proof}
We prove the formula for the upper metric mean dimension with potential.
The lower case is similar.
Let $\eta$ be an arbitrary positive number.
We take a positive number $\delta_0$ satisfying 
\[  \mathbf{d}(x, y) < \delta_0 \Longrightarrow  \left|\varphi(x)  - \varphi(y)\right| < \eta. \]
We take arbitrary $\varepsilon \in (0, 1)$ and $\delta\in (0, \delta_0)$.

Let $x\in \mathcal{X}$ and $L>0$.
For every $y\in B_\delta\left(x, \mathbf{d}_{\mathbb{R}^k}\right)$ we have 
$\mathbf{d}_L(x, y) \leq \delta < \delta_0$ and hence $\left|\varphi_L(x) - \varphi_L(y)\right| < \eta L^k$.
Hence
\begin{align*}
  \left(\frac{1}{\varepsilon}\right)^{\varphi_L(x)- \eta L^k}\#\left(B_\delta\left(x, \mathbf{d}_{\mathbb{R}^k}\right), \mathbf{d}_L, \varepsilon\right)
  & \leq  \#\left(B_\delta\left(x, \mathbf{d}_{\mathbb{R}^k}\right), \mathbf{d}_L, \varphi_L, \varepsilon\right) \\
  & \leq  \left(\frac{1}{\varepsilon}\right)^{\varphi_L(x)+ \eta L^k}
  \#\left(B_\delta\left(x, \mathbf{d}_{\mathbb{R}^k}\right), \mathbf{d}_L, \varepsilon\right).
\end{align*}
Taking the logarithm and dividing them by $L^k \log (1/\varepsilon)$,
\begin{align*}
  \frac{\log \#\left(B_\delta\left(x, \mathbf{d}_{\mathbb{R}^k}\right), \mathbf{d}_L, \varepsilon\right)}{L^k \log (1/\varepsilon)}
  + \frac{\varphi_L(x)}{L^k} - \eta  &\leq  
  \frac{\log \#\left(B_\delta\left(x, \mathbf{d}_{\mathbb{R}^k}\right), \mathbf{d}_L, \varphi_L, \varepsilon\right)}{L^k \log (1/\varepsilon)} \\
  & \leq \frac{\log \#\left(B_\delta\left(x, \mathbf{d}_{\mathbb{R}^k}\right), \mathbf{d}_L, \varepsilon\right)}{L^k \log (1/\varepsilon)}
  + \frac{\varphi_L(x)}{L^k} + \eta.
\end{align*}
Taking the lower limit as $L\to \infty$,
\[ \left|\varliminf_{L\to \infty}
\frac{\log \#\left(B_\delta\left(x, \mathbf{d}_{\mathbb{R}^k}\right), \mathbf{d}_L, \varphi_L, \varepsilon\right)}{L^k \log (1/\varepsilon)}
- \varliminf_{L\to \infty} 
  \left(\frac{\log \#\left(B_\delta\left(x, \mathbf{d}_{\mathbb{R}^k}\right), \mathbf{d}_L, \varepsilon\right)}{L^k \log (1/\varepsilon)}
  + \frac{\varphi_L(x)}{L^k}\right)\right| \leq \eta. \]
We take the supremum over $x\in \mathcal{X}$ and let $\varepsilon\to 0$.
By Theorem \ref{theorem: local formula of metric mean dimension with potential} we have 
\[ \left|\umdimm\left(\mathcal{X}, T, \mathbf{d}, \varphi\right) - 
    \varlimsup_{\varepsilon \to 0} \sup_{x\in X} \varliminf_{L\to \infty} 
    \left(\frac{\log \#\left(B_\delta\left(x, \mathbf{d}_{\mathbb{R}^k}\right), \mathbf{d}_L, \varepsilon\right)}{L^k \log (1/\varepsilon)}
  + \frac{\varphi_L(x)}{L^k}\right)\right| \leq \eta. \]
We let $\delta \to 0$:
\[    \left|\umdimm\left(\mathcal{X}, T, \mathbf{d}, \varphi\right) - 
    \lim_{\delta\to 0} \varlimsup_{\varepsilon \to 0} \sup_{x\in X} \varliminf_{L\to \infty} 
    \left(\frac{\log \#\left(B_\delta\left(x, \mathbf{d}_{\mathbb{R}^k}\right), \mathbf{d}_L, \varepsilon\right)}{L^k \log (1/\varepsilon)}
  + \frac{\varphi_L(x)}{L^k}\right)\right| \leq \eta. \]
Here $\eta>0$ is arbitrary. So we let $\eta\to 0$ and get the claim of the proposition. 
\end{proof}

\section{Nondegenerate Brody curves and the proof of Theorem \ref{theorem: Ruelle inequality for Brody curves}} 
\label{section: proof of the Ruelle inequality for Brody curves}

The purpose of this section is to prove Theorem \ref{theorem: Ruelle inequality for Brody curves}
(a Ruelle inequality for Brody curves).

\subsection{Nondegenerate Brody curves}  \label{subsection: nondegenerate Brody curves}

A basic tool to study $\mathcal{B}^N$ is the deformation theory of Brody curves developed in 
\cite{Tsukamoto_deformation, Matsuo--Tsukamoto, Tsukamoto_mean_dimension_Brody_curves}.
This corresponds to the theory of stable and unstable manifolds in the case of Axiom A diffeomorphisms.
The (un)stable manifold theorem requires a hyperbolicity assumption
\cite[3.2. Theorem]{Bowen_lecture}.
Similarly, we need a kind of “hyperbolicity” (transversality) for developing the deformation theory of Brody curves.
The next definition provides such a notion.
This stems from a classical paper of Yosida \cite{Yosida}.

\begin{definition} \label{definition: nondegenerate Brody curves}
Let $R$ be a positive number.
A holomorphic map $f\colon \mathbb{C}\to \mathbb{C}P^N$ is said to be \textbf{$R$-nondegenerate}
if for all points $a\in \mathbb{C}$ we have 
  \begin{equation}  \label{eq: nondegenerate Brody curves}
    \max_{|z-a|\leq R} |df|(z) \geq \frac{1}{R}. 
  \end{equation}  
Namely, $f$ is $R$-nondegenerate if every $R$-disk of the plane contains a point $z$ for which we have $|df|(z) \geq 1/R$.
\end{definition}

Notice that if $f$ is $R$-nondegenerate and $R^\prime \geq R$ then $f$ is $R^\prime$-nondegenerate.
We also say that a holomorphic map $f\colon \mathbb{C}\to \mathbb{C}P^N$ 
is \textbf{nondegenerate} if $f$ is $R$-nondegenerate for some $R>0$.

A typical example of \textit{degenerate} (i.e. not nondegenerate) holomorphic map is a constant map.
The condition (\ref{eq: nondegenerate Brody curves}) implies that $f$ is not close to constant maps 
over every $R$-disk of the plane\footnote{More precisely, $f$ is nondegenerate if and only if the $\mathbb{C}$-orbit of $f$ 
does not accumulate to a constant map. This condition was first introduced by Yosida \cite{Yosida} for meromorphic functions.
The condition (\ref{eq: nondegenerate Brody curves}) is its quantitative version.}.
A rational function is degenerate because it converges to a constant function near the infinity of the plane.
A typical example of nondegenerate holomorphic maps is given by a nonconstant elliptic function.

Definition \ref{definition: nondegenerate Brody curves} plays a role of “hyperbolicity” in the deformation theory of 
Brody curves.
Here we briefly review the main ideas of this theory.
Given a Lipschitz holomorphic map $f\colon \mathbb{C}\to \mathbb{C} P^N$, 
we would like to study a deformation of $f$. 
Let $T\mathbb{C}P^N$ be the tangent bundle of $\mathbb{C}P^N$.
Its Hermitian metric is given by the Fubini--Study metric.
Let $E := f^*T\mathbb{C}P^N$ be the pull-back of $T\mathbb{C}P^N$ by the map $f\colon \mathbb{C}\to \mathbb{C}P^N$.
This is a holomorphic vector bundle over the complex plane.
Its Hermitian metric is given by the pull-back of the Fubini--Study metric.
Let $\bar{\partial}\colon C^\infty(E) \to C^\infty\left(\Lambda^{0,1}(E)\right)$ be the Dolbeault operator.
We define $H_f$ as the space of holomorphic sections $u$ of $E$ satisfying $\norm{u}_{L^\infty(\mathbb{C})} < \infty$.
This is the space of \textit{first-order deformations}\footnote{The condition that a given map is holomorphic is a non-linear equation.
The linearization of this equation is given by the Dolbeault operator $\bar{\partial}$.
So $H_f$ is the space of the solutions of this linearized equation.}.
As is well-known in the standard deformation theory, 
we need a \textit{transversality}\footnote{Transversality is a condition that we need for applying the implicit function theorem.} 
for lifting the first-order deformation to true deformation.
The relevant transversality condition here is the surjectivity of $\bar{\partial}\colon C^\infty(E) \to C^\infty\left(\Lambda^{0,1}(E)\right)$
under a suitable norm restriction.

We have the Bochner formula
\[  \bar{\partial} \bar{\partial}^* \omega 
   = \frac{1}{2}\nabla^* \nabla \omega + \Theta \omega, 
    \quad  (\omega \in C^\infty\left(\Lambda^{0,1}(E)\right)). \]
Here $\bar{\partial}^*\colon C^\infty\left(\Lambda^{0,1}(E)\right) \to C^\infty(E)$ is the formal adjoint of $\bar{\partial}$, and 
$\nabla$ is the canonical connection with respect to the Hermitian metric on $E$.
The operator $\Theta := [\nabla_{\partial/\partial z}, \nabla_{\partial/\partial \bar{z}}]$ is the curvature.
The holomorphic bisectional curvature of the Fubini--Study metric is bounded from below by $2\pi$. Then
\[  \langle \Theta \omega, \omega \rangle \geq \pi |df|^2 \cdot |\omega|^2. \]
Therefore the condition (\ref{eq: nondegenerate Brody curves}) implies that $\Theta$ has a certain positivity property.

Suppose $f$ is nondegenerate. 
Then for any $ \omega \in C^\infty\left(\Lambda^{0,1}(E)\right)$
with $\norm{\omega}_{L^\infty} < \infty$,
by using the Bochner formula and the positivity of $\Theta$, 
we can construct $F(\omega)\in C^\infty(E)$ which depends linearly on $\omega$ and satisfies 
\begin{equation} \label{eq: right inverse of the Dolbeault operator}
    \bar{\partial} \left(F(\omega)\right) = \omega, \quad  \norm{F(\omega)}_{L^\infty} \lesssim \norm{\omega}_{L^\infty}. 
\end{equation}    
This is the transversality condition we need.
From this condition, for any $u\in H_f$ with small $L^\infty$-norm, we can construct $\alpha(u)\in C^\infty(E)$
such that 
\begin{itemize}
   \item   $\norm{\alpha(u)}_{L^\infty} =  O\left(\norm{u}_{L^\infty}^2\right)$,
   \item   the map $g(z) := \exp_{f(z)}\left(u(z)+\alpha(u)(z)\right)$ $(z\in \mathbb{C})$ is holomorphic. Here 
              $\exp_p\colon T_p\mathbb{C}P^N \to \mathbb{C}P^N$ is the exponential map with respect to the Fubini--Study metric. 
\end{itemize}
Moreover, every holomorphic curve sufficiently close to $f$ uniformly over the plane can be written as 
$g(z) = \exp_{f(z)}\left(u(z)+\alpha(u)(z)\right)$.

In summary, if $f\colon \mathbb{C}\to \mathbb{C}P^N$ is a nondegenerate Lipschitz holomorphic curve,
then we can describe a local deformation of $f$ by using the linear space $H_f$.
This is the idea behind Definition \ref{definition: nondegenerate Brody curves}.

Let $\lambda$ be a positive number.
We define $\mathcal{B}^N_\lambda$ as the space of holomorphic maps $f\colon \mathbb{C}\to \mathbb{C}P^N$ satisfying 
$|df|(z) \leq  \lambda$ for all $z\in \mathbb{C}$.
Namely it is the space of $\lambda$-Lipschitz holomorphic curves.
The case $\lambda=1$ coincides with the space of Brody curves $\mathcal{B}^N$.
We are mainly interested in the case of $\lambda=1$, but it is technically convenient to consider 
the case $\lambda>1$ also\footnote{We can transform every result for $\mathcal{B}^N$ to 
that of $\mathcal{B}^N_\lambda$ by a scale change.}.
The group $\mathbb{C}$ acts on $\mathcal{B}^N_\lambda$ by translation.
We also denote it by $T$:
\[  T\colon \mathbb{C}\times \mathcal{B}^N_\lambda \to  \mathcal{B}^N_\lambda, \quad 
     \left(a, f(z)\right)  \mapsto  f(z+a). \]
We define a metric on $\mathcal{B}^N_\lambda$ (also denoted by $\mathbf{d}$) by 
\[   \mathbf{d}\left(f, g\right) = \max_{z\in [0,1]^2} d_{\mathrm{FS}}\left(f(z), g(z)\right). \]
Here recall that $d_{\mathrm{FS}}$ is the distance function on $\mathbb{C}P^N$ defined by the Fubini--Study metric.

The following proposition is a main consequence of the deformation theory described above.
This was proved in \cite[Proposition 3.3]{Tsukamoto_mean_dimension_Brody_curves}.
(Indeed the paper \cite[Proposition 3.3]{Tsukamoto_mean_dimension_Brody_curves} provided a slightly stronger statement than this.)

\begin{proposition} \label{prop: covering number around nondegenerate curves}
For any $R>0$ and $0<\varepsilon<1$ there exist\footnote{Here $\delta_1$ and $C_1$ depend only on 
$R$ and are independent of $\varepsilon$,
while $C_2$ depends only on $\varepsilon$ and is independent of $R$.} positive numbers $\delta_1 =\delta_1(R)$, $C_1 = C_1(R)$
and $C_2=C_2(\varepsilon)$ for which the following statement holds:
For any $R$-nondegenerate $f\in \mathcal{B}^N_2$ and for any square $S\subset \mathbb{C}$ of 
side length $L\geq 1$ we have 
\[  \#\left(\left\{g\in \mathcal{B}^N_2\middle| \, \mathbf{d}_{\mathbb{C}}\left(f, g\right) \leq \delta_1\right\}, 
               \mathbf{d}_S, \varepsilon\right)
     \leq   \left(\frac{C_1}{\varepsilon}\right)^{2(N+1)\int_{S} |df|^2\, dxdy + C_2 L}. \]          
\end{proposition}

Here recall that $\mathbf{d}_{\mathbb{C}}$ and $\mathbf{d}_S$ are metrics on $\mathcal{B}^N_2$ defined by 
\[ \mathbf{d}_{\mathbb{C}}(f, g) = \sup_{u\in \mathbb{C}}  \mathbf{d}\left(T^u f, T^u g\right)
    = \sup_{z\in \mathbb{C}} d_{\mathrm{FS}}\left(f(z), g(z)\right), \quad 
    \mathbf{d}_S(f, g) = \sup_{u \in S} \mathbf{d}\left(T^u f, T^u g\right). \]
In particular, if $S = [0, L]^2 = \{x+y\sqrt{-1} \in \mathbb{C} \mid 0 \leq x, y \leq L\}$ 
then $\mathbf{d}_S(f, g)$ is given by $\sup_{z\in [0, L+1]^2}  d_{\mathrm{FS}}\left(f(z), g(z)\right)$.

\begin{proof}[Sketch of the proof of Proposition \ref{prop: covering number around nondegenerate curves}]
Here is a very brief sketch of the proof.
By the deformation theory, the proof boils down to a problem of estimating the $\varepsilon$-covering number of 
$B_{\delta_1}(H_f) := \{u\in H_f \mid  \norm{u}_{L^\infty(\mathbb{C})} \leq \delta_1\}$ with respect to the
metric $\norm{\cdot}_{L^\infty(S)}$.
We solve this problem by using the Riemann--Roch theorem.
\end{proof}

Proposition \ref{prop: covering number around nondegenerate curves} provides a good control of nondegenerate Brody curves 
and its neighborhoods\footnote{Notice that the metrics $\mathbf{d}$ and  
$\mathbf{d}_{\mathbb{C}}$ define different topology. In particular, a neighborhood with respect to $\mathbf{d}_{\mathbb{C}}$ is not 
a neighborhood with respect to $\mathbf{d}$ in general.}
with respect to the metric $\mathbf{d}_{\mathbb{C}}$.
However some Brody curves are degenerate.
They become “singularity” for our analysis.
The next proposition provides a “ resolution of singularity”.
For $r>0$ and a Brody curve $f\in \mathcal{B}^N$ we denote 
\[ B_r(f, \mathbf{d}_{\mathbb{C}}) = \{g\in \mathcal{B}^N\mid \mathbf{d}_{\mathbb{C}}(f, g) \leq r\}.   \]
Notice that this is a subset of $\mathcal{B}^N$.
We will not use this notation for other $\mathcal{B}^N_\lambda$ $(\lambda \neq 1)$.
So it may not cause a confusion.

\begin{proposition}  \label{prop: resolution of singularity}
Let $\eta$ be a positive number.
There exist\footnote{Here $\delta_2$ and $C_3$ are universal constants independent of $\eta$, while $R_1$ depends on $\eta$.} 
$\delta_2>0$, $C_3>1$ and $R_1 = R_1(\eta)>0$ for which the following statement holds.
For any Brody curve $f\in \mathcal{B}^N$ there exists a map 
$\Phi\colon B_{\delta_2}(f, \mathbf{d}_{\mathbb{C}})\to \mathcal{B}^N_{1+\eta}$
satisfying the following three conditions.
  \begin{enumerate}
    \item  $\Phi(f)$ is $R_1$-nondegenerate.
    \item  For any $g_1, g_2\in B_{\delta_2}(f, \mathbf{d}_{\mathbb{C}})$ and $z\in \mathbb{C}$ we have 
    \[   \dfs\left(\Phi(g_1)(z), \Phi(g_2)(z)\right) \leq  C_3\,  \dfs\left(g_1(z), g_2(z)\right), \]
    \[   \dfs\left(g_1(z), g_2(z)\right)  \leq  C_3 \sup_{|w-z|\leq 3} \dfs\left(\Phi(g_1)(w), \Phi(g_2)(w)\right). \]
             In particular, the image of $\Phi$ is contained in 
             $\{g\in \mathcal{B}^N_{1+\eta}\mid  \mathbf{d}_{\mathbb{C}} \left(\Phi(f), g\right) \leq  C_3 \delta_2\}$.
     \item  For any square $S\subset \mathbb{C}$ of side length $L\geq R_1$
     \[   \left|\int_S |df|^2 \, dxdy - \int_S \left|d\Phi(f)\right|^2\, dxdy \right| < \eta L^2. \]     
     Here $\left|d\Phi(f)\right|$ is the local Lipschitz constant of the map $\Phi(f)\colon \mathbb{C} \to \mathbb{C}P^N$.
  \end{enumerate}
\end{proposition}

This statement may look technical.
Its idea is as follows.
Let $f\in \mathcal{B}^N$ be a degenerate Brody curve.
We want to replace $f$ with a nondegenerate one.
The condition (1) means that the map $\Phi$ provides such a replacement.
The condition (2) means that this replacement procedure does not distort much the metric structure around $f$.
The condition (3) means that $\Phi$ changes the values of the energy integrals only slightly.
 
Proposition \ref{prop: resolution of singularity} is a modification of \cite[Proposition 3.3]{Tsukamoto_mean_dimension_Brody_curves}.
We will prove it in \S \ref{subsection: proof of resolution of singularity}.

\subsection{Proof of Theorem \ref{theorem: Ruelle inequality for Brody curves}} \label{subsection: Proof of Ruelle inequality}

We prove Theorem \ref{theorem: Ruelle inequality for Brody curves} (a Ruelle inequality for Brody curves) here.
Recall some notations: $\mathcal{B}^N$ is the space of Brody curves in $\mathbb{C}P^N$ with the metric 
$\mathbf{d}(f, g) = \max_{z\in [0,1]^2} \dfs\left(f(z), g(z)\right)$.
The function $\psi\colon \mathcal{B}^N\to \mathbb{R}$ is defined by 
\[  \psi(f) = 2(N+1) \, |df|^2(0). \]
For $L>0$ the metric $\mathbf{d}_L$ and function $\psi_L$ are defined by
\[ \mathbf{d}_L(f, g) = \max_{u\in [0,L]^2}\mathbf{d}\left(T^u f, T^u g\right) 
     = \max_{z\in [0,L+1]^2} \dfs\left(f(z), g(z)\right),  \]
\begin{equation} \label{eq: psi_L}
   \psi_L(f) = \int_{[0,L]^2} \psi\left(T^u f\right)\, d\mathbf{m}(u) 
      = 2(N+1)\int_{[0,L]^2} |df|^2(z) \, dxdy. 
\end{equation} 

The next theorem is a key result.
This is a “Brody curve analogue” of the equation $P_T(-\phi) = 0$ of Axiom A attractors (\ref{eq: Topological pressure is zero}).

\begin{theorem} \label{theorem: topological and metric mean dimensions with potential are zero}
The topological and metric mean dimensions of $(\mathcal{B}^N, T, \mathbf{d})$ with respect to the potential function $-\psi$ 
are both zero:
\[  \mdim\left(\mathcal{B}^N, T, -\psi\right) = \mdimm\left(\mathcal{B}^N, T, \mathbf{d},-\psi\right) = 0. \]
\end{theorem}

\begin{proof}
By (\ref{eq: mdim lmdimm umdimm}) we have 
\[  \mdim\left(\mathcal{B}^N, T, -\psi\right) \leq \lmdimm \left(\mathcal{B}^N, T, \mathbf{d},-\psi\right)
     \leq  \umdimm \left(\mathcal{B}^N, T, \mathbf{d},-\psi\right). \]
Let $\mathcal{C}\subset \mathcal{B}^N$ be the space of constant maps.
This is a closed invariant subset of $\mathcal{B}^N$.
The group $\mathbb{C}$ trivially acts on it, and the function $\psi$ is identically zero on it.
Therefore we have 
\[  \mdim\left(\mathcal{B}^N, T, -\psi\right) \geq \mdim\left(\mathcal{C}, T, -\psi\right) = 0. \]
Hence it is enough to prove $\umdimm \left(\mathcal{B}^N, T, \mathbf{d},-\psi\right) \leq 0$.

By Proposition \ref{prop: modified local formula}
\[  \umdimm\left(\mathcal{B}^N, T,\mathbf{d}, -\psi\right) 
    = \lim_{\delta\to 0} \varlimsup_{\varepsilon \to 0} \sup_{f\in \mathcal{B}^N} \varliminf_{L\to \infty}
   \left(\frac{\log \#\left(B_\delta(f, \mathbf{d}_{\mathbb{C}}), \mathbf{d}_L, \varepsilon\right)}{L^2 \log (1/\varepsilon)}
    - \frac{\psi_L(f)}{L^2}\right). \]
By (\ref{eq: psi_L}) this can be also written as
\begin{equation}  \label{eq: upper metric mean dimension with potential of Brody curves}
  \begin{split}
  & \umdimm\left(\mathcal{B}^N, T,\mathbf{d}, -\psi\right)  \\
  &  = \lim_{\delta\to 0} \varlimsup_{\varepsilon \to 0} \sup_{f\in \mathcal{B}^N} \varliminf_{L\to \infty}
   \left(\frac{\log \#\left(B_\delta(f, \mathbf{d}_{\mathbb{C}}), \mathbf{d}_L, \varepsilon\right)}{L^2 \log (1/\varepsilon)}
    - \frac{2(N+1)}{L^2} \int_{[0,L]^2} |df|^2\, dxdy \right). 
  \end{split}  
\end{equation}    
We estimate the right-hand side.

Let $\eta$ be an arbitrary positive number with $\eta <1$.
Let $\delta_2>0, C_3>1$ and $R_1 = R_1(\eta)>0$ be positive numbers introduced in Proposition 
\ref{prop: resolution of singularity} (resolution of singularity).
Take an arbitrary $0<\varepsilon < 1$.
Let $\delta_1 = \delta_1(R_1)$, $C_1 = C_1(R_1)$ and $C_2 = C_2(\varepsilon)$ be positive numbers introduced in 
Proposition \ref{prop: covering number around nondegenerate curves} (covering number estimate around nondegenerate curves)
with respect to these $R_1$ and $\varepsilon$.

Let $0<\delta \leq \min(\delta_2, \delta_1/C_3)$.
Let $f\in \mathcal{B}^N$ be an arbitrary Brody curve.
By Proposition \ref{prop: resolution of singularity}, there exists a map 
$\Phi\colon B_\delta(f, \mathbf{d}_{\mathbb{C}}) \to 
\{g \in \mathcal{B}^N_{1+\eta}\mid \mathbf{d}_{\mathbb{C}}\left(\Phi(f), g\right) \leq C_3\delta\}$
such that 
  \begin{itemize}
   \item  $\Phi(f)$ is $R_1$-nondegenerate,
   \item  for any $g_1, g_2\in B_\delta(f, \mathbf{d}_{\mathbb{C}})$ and $z\in \mathbb{C}$ we have
   \[  \dfs\left(g_1(z), g_2(z)\right) \leq  C_3 \max_{|w-z|\leq 3} \dfs\left(\Phi(g_1)(w), \Phi(g_2)(w)\right), \]
   \item for any $L\geq R_1$
   \begin{equation} \label{eq: energy integral is not changed much}
       \left|\int_{[0,L]^2}|df|^2\, dxdy - \int_{[0,L]^2}\left|d\Phi(f)\right|^2\, dxdy \right| < \eta L^2. 
   \end{equation}   
  \end{itemize}
For $g_1, g_2\in B_\delta(f, \mathbf{d}_{\mathbb{C}})$ we have
\begin{align*}
   \mathbf{d}_L(g_1, g_2) & = \max_{z\in [0,L+1]^2} \dfs\left(g_1(z), g_2(z)\right) \\
   & \leq  C_3 \max_{w \in [-3, L+4]^2} \dfs\left(\Phi(g_1)(w), \Phi(g_2)(w)\right) \\
   &  = C_3 \, \mathbf{d}_{[-3, L+3]^2} \left(\Phi(g_1), \Phi(g_2)\right).
\end{align*}
Therefore 
\begin{align*}
   \#\left(B_\delta(f, \mathbf{d}_{\mathbb{C}}), \mathbf{d}_L, \varepsilon\right) &
   \leq  \#\left(\left\{g \in \mathcal{B}^N_{1+\eta}\mid \mathbf{d}_{\mathbb{C}}\left(\Phi(f), g\right) \leq C_3\delta\right\}, 
   \mathbf{d}_{[-3, L+3]^2}, \frac{\varepsilon}{C_3}\right) \\
   & \leq \#\left(\left\{g \in \mathcal{B}^N_{1+\eta}\mid \mathbf{d}_{\mathbb{C}}\left(\Phi(f), g\right) \leq \delta_1 \right\}, 
   \mathbf{d}_{[-3, L+3]^2}, \frac{\varepsilon}{C_3}\right)  \quad 
   \text{by $\delta\leq \frac{\delta_1}{C_3}$} \\
   & \leq  \left(\frac{C_1 C_3}{\varepsilon}\right)^{2(N+1)\int_{[-3,L+3]^2}\left|d\Phi(f)\right|^2 \, dxdy + C_2 (L+6)}.
\end{align*}
Here we have used the covering number estimate of Proposition \ref{prop: covering number around nondegenerate curves}
in the last line.
We take the logarithm and divide them by $L^2$:
\[ \frac{\log  \#\left(B_\delta(f, \mathbf{d}_{\mathbb{C}}), \mathbf{d}_L, \varepsilon\right)}{L^2}  \leq 
  \left(\frac{2(N+1)}{L^2} \int_{[-3,L+3]^2}\left|d\Phi(f)\right|^2 dxdy + \frac{C_2(L+6)}{L^2}\right)\log \left(\frac{C_1 C_3}{\varepsilon}\right).\]
Now we assume $L\geq R_1$ and use (\ref{eq: energy integral is not changed much}). 
Then the right-hand side is bounded from above by
\[ \left(\frac{2(N+1)}{L^2} \int_{[-3,L+3]^2}\left|df\right|^2\, dxdy + \frac{2(N+1)(L+6)^2 \eta}{L^2}
    + \frac{C_2(L+6)}{L^2}\right)\log \left(\frac{C_1 C_3}{\varepsilon}\right).\]
We have 
\[  \frac{2(N+1)}{L^2} \int_{[-3,L+3]^2}\left|df\right|^2\, dxdy = \frac{2(N+1)}{L^2}\int_{[0,L]^2}|df|^2 \, dxdy
     + O\left(\frac{1}{L}\right). \]
Therefore
\begin{align*}
   &\varliminf_{L\to \infty} \left(\frac{\log  \#\left(B_\delta(f, \mathbf{d}_{\mathbb{C}}), \mathbf{d}_L, \varepsilon\right)}{L^2 \log (1/\varepsilon)}  
   - \frac{2(N+1)}{L^2} \int_{[0, L]^2}|df|^2\, dxdy\right) \\
   & \leq  2(N+1) \frac{\log (C_1 C_3)}{\log (1/\varepsilon)} 
      + 2(N+1) \eta \frac{\log \left(\frac{C_1 C_3}{\varepsilon}\right)}{\log (1/\varepsilon)}.
\end{align*}
Here $C_1 = C_1\left(R_1(\eta)\right)$ and $C_3$ are independent of $\varepsilon$.
We take the supremum over $f\in \mathcal{B}^N$ and let $\varepsilon \to 0$.
Then 
\[  \varlimsup_{\varepsilon \to 0} \sup_{f\in \mathcal{B}^N} \varliminf_{L\to \infty} 
    \left(\frac{\log  \#\left(B_\delta(f, \mathbf{d}_{\mathbb{C}}), \mathbf{d}_L, \varepsilon\right)}{L^2 \log (1/\varepsilon)}  
   - \frac{2(N+1)}{L^2} \int_{[0, L]^2}|df|^2\, dxdy\right) 
    \leq  2(N+1)\eta. \]
We let $\delta \to 0$. 
By (\ref{eq: upper metric mean dimension with potential of Brody curves}) we have
\[ \umdimm\left(\mathcal{B}^N, T,\mathbf{d}, -\psi\right)  \leq  2(N+1)\eta. \]
Recall that $\eta$ is an arbitrary number with $0<\eta<1$.
Therefore we conclude 
\[  \umdimm\left(\mathcal{B}^N, T,\mathbf{d}, -\psi\right) \leq 0.  \]
\end{proof}

Now we can prove Theorem \ref{theorem: Ruelle inequality for Brody curves}.
We write the statement again.

\begin{corollary}[$=$ Theorem \ref{theorem: Ruelle inequality for Brody curves}]
For every invariant probability measure $\mu\in \mathscr{M}^T\left(\mathcal{B}^N\right)$ we have 
\[   \urdim\left(\mathcal{B}^N, T, \mathbf{d}, \mu\right) \leq  \int_{\mathcal{B}^N} \psi\, d\mu. \]
\end{corollary}

\begin{proof}
By the variational principle (Theorem \ref{theorem: variational principle for mean dimension})
\[  \urdim\left(\mathcal{B}^N, T, \mathbf{d}, \mu\right) - \int_{\mathcal{B}^N} \psi \, d\mu 
     \leq  \umdimm\left(\mathcal{B}^N, T, \mathbf{d}, -\psi\right). \]
The right-hand side is zero by Theorem \ref{theorem: topological and metric mean dimensions with potential are zero}.
Thus 
\[  \urdim\left(\mathcal{B}^N, T, \mathbf{d}, \mu\right) - \int_{\mathcal{B}^N} \psi \, d\mu \leq 0. \]
\end{proof}

\begin{remark}  \label{remark: ordinary topological and metric mean dimension of Brody curves}
By Theorem \ref{theorem: topological and metric mean dimensions with potential are zero} and 
the variational principle (Theorem \ref{theorem: variational principle for mean dimension}) we have
\begin{align*}
   \mdim\left(\mathcal{B}^N, T, -\psi\right) 
   & = \sup_{\mu\in \mathscr{M}^T(\mathcal{B}^N)} \left(\lrdim\left(\mathcal{B}^N, T, \mathbf{d}, \mu\right)
    - \int_{\mathcal{B}^N} \psi \, d\mu\right)  \\
   & = \sup_{\mu \in  \mathscr{M}^T(\mathcal{B}^N)} \left(\urdim\left(\mathcal{B}^N, T, \mathbf{d}, \mu\right)
    - \int_{\mathcal{B}^N} \psi \, d\mu\right)  \\
   & = \mdimm\left(\mathcal{B}^N, T, \mathbf{d},-\psi\right) \\
   & = 0.
\end{align*}   
This shows that the conjecture stated in Remark \ref{remark: double variational principle conjecture}
holds for $\left(\mathcal{B}^N, T, -\psi\right)$.
We also have a similar result for the potential function $\varphi := 0$ (the function constantly equal to zero).
Indeed, on the one hand, the paper \cite[Theorem 1.1]{Matsuo--Tsukamoto} showed that 
$\mdim\left(\mathcal{B}^N, T\right) \geq  2(N+1) \rho(\mathbb{C}P^N)$.
On the other hand, another paper \cite[p. 947]{Tsukamoto_mean_dimension_Brody_curves} 
proved\footnote{Strictly speaking, the paper \cite{Tsukamoto_mean_dimension_Brody_curves} stated only 
$\lmdimm\left(\mathcal{B}^N, T, \mathbf{d}\right) \leq 2(N+1) \rho(\mathbb{C}P^N)$.
However the argument of \cite[pp. 947-948]{Tsukamoto_mean_dimension_Brody_curves} can be also applied to 
the upper metric mean dimension without any changes. 
Another paper \cite[\S 4.3]{Tsukamoto_local_nature} explicitly proved 
$\umdimm\left(\mathcal{B}^N, T, \mathbf{d}\right) \leq 2(N+1) \rho(\mathbb{C}P^N)$.}
that $\umdimm\left(\mathcal{B}^N, T, \mathbf{d}\right) \leq 2(N+1) \rho(\mathbb{C}P^N)$.
Since $\mdim\left(\mathcal{B}^N, T\right) \leq \lmdimm\left(\mathcal{B}^N, T, \mathbf{d}\right) 
\leq  \umdimm\left(\mathcal{B}^N, T, \mathbf{d}\right)$, we have 
\[ \mdim\left(\mathcal{B}^N, T\right) = \mdimm\left(\mathcal{B}^N, T, \mathbf{d}\right) = 2(N+1) \rho(\mathbb{C}P^N). \]
By the variational principle (Theorem \ref{theorem: variational principle for mean dimension}), we conclude
\begin{align*}
  \mdim\left(\mathcal{B}^N, T\right)  
  & = \sup_{\mu\in \mathscr{M}^T(\mathcal{B}^N)} \lrdim\left(\mathcal{B}^N, T, \mathbf{d}, \mu\right) \\
  & = \sup_{\mu\in \mathscr{M}^T(\mathcal{B}^N)} \urdim\left(\mathcal{B}^N, T, \mathbf{d}, \mu\right)  \\
  & = \mdimm \left(\mathcal{B}^N, T, \mathbf{d}\right) \\
  & = 2(N+1)\rho(\mathbb{C}P^N).
\end{align*}
\end{remark}

\subsection{Proof of Proposition \ref{prop: resolution of singularity}} \label{subsection: proof of resolution of singularity}

In this section we prove Proposition \ref{prop: resolution of singularity}.
The proof is based on the next lemma. 
This was proved in \cite[Lemma 4.1]{Tsukamoto_mean_dimension_Brody_curves}.
We need to introduce some notations.
We define $\mathcal{H}^N$ as the set of all holomorphic maps $f\colon \mathbb{C}\to \mathbb{C}P^N$.
Let $r>0$, $p\in \mathbb{C}$ and $q\in \mathbb{C}P^N$. 
We denote $D_r(p) = \{z\in \mathbb{C}\mid |z-p|\leq r\}$ and
$B_r(q) = \{w \in \mathbb{C}P^N\mid \dfs(w, q) \leq r\}$.

\begin{lemma}[Gluing lemma] \label{lemma: gluing}
There exist $\delta_3>0, R_2>0$ and $C_4>1$ for which the following statement holds.
Let $p\in \mathbb{C}$, $q\in \mathbb{C}P^N, R\geq R_2$ and define 
$\mathcal{A} = \{f\in \mathcal{H}^N\mid  f\left(D_R(p)\right) \subset B_{\delta_3}(q)\}$.
Then there exists a map $\Psi\colon \mathcal{A} \to \mathcal{H}^N$
satisfying the following conditions.    
  \begin{enumerate}
     \item  For any $f\in \mathcal{A}$ 
              \[  \dfs\left(f(z), \Psi(f)(z)\right) \leq  \frac{C_4}{|z-p|^3}, \quad  (z\neq p).  \]
     \item  We have $\frac{1}{100} < \norm{d\, \Psi(f)}_{L^\infty\left(D_{R/2}(p)\right)} < 1$ for all $f\in \mathcal{A}$.
               Here $\norm{d\, \Psi(f)}_{L^\infty\left(D_{R/2}(p)\right)}$ is the maximum of the local Lipschitz constant of the map $\Psi(f)$
               over the disk $D_{R/2}(p)$.
     \item  For any $f\in \mathcal{A}$ and $z\in \mathbb{C}$ with $|z-p|\geq 1$, we have
              \[  \left||df|(z)- |d\, \Psi(f)|(z)\right| \leq  \frac{C_4}{|z-p|^3} |df|(z) + \frac{C_4}{|z-p|^3}. \]
     \item  Let $f, g \in \mathcal{A}$.
               For any $z\in \mathbb{C}$ with $|z-p|\leq 1$ we have 
               \[  C_4^{-1}\dfs\left(\Psi(f)(z), \Psi(g)(z)\right) \leq \dfs\left(f(z), g(z)\right) \leq 
                    C_4 \sup_{|w-p|\leq 2} \dfs\left(\Psi(f)(w), \Psi(g)(w)\right). \]
               Moreover, for any $z\in \mathbb{C}$ with $|z-p|\geq 1$, we have 
               \begin{align*}
                 \dfs\left(\Psi(f)(z), \Psi(g)(z)\right) & \leq  \left(1+\frac{C_4}{|z-p|^3}\right) \dfs\left(f(z), g(z)\right), \\
                 \dfs\left(f(z), g(z)\right) & \leq  \left(1+\frac{C_4}{|z-p|^3}\right) \dfs\left(\Psi(f)(z), \Psi(g)(z)\right).
               \end{align*}
  \end{enumerate}
\end{lemma}

\begin{proof}[Sketch of the proof]
We explain a brief sketch of the proof for convenience of readers.
The details were given in \cite[Lemma 4.1]{Tsukamoto_mean_dimension_Brody_curves}.
We can assume $p=0$ and $q= [1:0:\dots:0]$ without loss of generality.
Suppose $f\in \mathcal{H}^N$ satisfies $f\left(D_R(0)\right) \subset B_{\delta_2}(q)$.
We assume  that $R$ is sufficiently large and $\delta_2$ is sufficiently small.
Write $f(z) = [1:f_1(z):\dots:f_N(z)]$ in the homogeneous coordinates.
Here $f_n(z)$ are meromorphic functions in $\mathbb{C}$ satisfying $|f_n(z)|\lesssim \delta_2$ on $D_R(0)$.
Therefore $f$ is almost constant over the large disk $D_R(0)$.
We eliminate this “almost constant region” by gluing a rational curve to $f$.

Let $a$ be a positive number and consider a holomorphic map $h_a\colon \mathbb{C}\to \mathbb{C}P^N$ defined by 
$h_a(z) = \left[1:\frac{a}{z^3}:\dots: \frac{a}{z^3}\right]$.
We fix $a>0$ so that $\norm{h_a}_{L^\infty(\mathbb{C})}=\frac{1}{10}$.
Define $\Psi(f)\in \mathcal{H}^N$ by 
\[  \Psi(f)(z)  = \left[1: f_1(z) + \frac{a}{z^3}:\dots: f_N(z)+\frac{a}{z^3}\right]. \]
Geometrically $\Psi(f)$ is a “gluing” of $f$ and $h_a$.
We can check that it satisfies all the requirements by (rather tedious) calculations.
\end{proof}

Now we can start the proof of Proposition \ref{prop: resolution of singularity}.
We write the statement again.

\begin{proposition}[$=$ Proposition \ref{prop: resolution of singularity}]
Let $\eta$ be a positive number.
There exist $\delta_2>0$, $C_3>1$ and $R_1 = R_1(\eta)>0$ for which the following statement holds.
For any Brody curve $f\in \mathcal{B}^N$ there exists a map 
$\Phi\colon B_{\delta_2}(f, \mathbf{d}_{\mathbb{C}})\to \mathcal{B}^N_{1+\eta}$
satisfying the following three conditions.
  \begin{enumerate}
    \item  $\Phi(f)$ is $R_1$-nondegenerate.
    \item  For any $g_1, g_2\in B_{\delta_2}(f, \mathbf{d}_{\mathbb{C}})$ and $z\in \mathbb{C}$ we have 
    \[   \dfs\left(\Phi(g_1)(z), \Phi(g_2)(z)\right) \leq  C_3\,  \dfs\left(g_1(z), g_2(z)\right), \]
    \[   \dfs\left(g_1(z), g_2(z)\right)  \leq  C_3 \sup_{|w-z|\leq 3} \dfs\left(\Phi(g_1)(w), \Phi(g_2)(w)\right). \]
             In particular, the image of $\Phi$ is contained in 
             $\{g\in \mathcal{B}^N_{1+\eta}\mid  \mathbf{d}_{\mathbb{C}} \left(\Phi(f), g\right) \leq  C_3 \delta_2\}$.
     \item  For any square $S\subset \mathbb{C}$ of side length $L\geq R_1$
     \[   \left|\int_S |df|^2 \, dxdy - \int_S \left|d\Phi(f)\right|^2\, dxdy \right| < \eta L^2. \]        
  \end{enumerate}
\end{proposition}

\begin{proof}
The proof is close to \cite[Proposition 3.3]{Tsukamoto_mean_dimension_Brody_curves}.
First we recall elementary facts about infinite products:
Let $a_n$ $(n\geq 1)$ be a sequence of positive numbers.
If the sum $\sum_{n=1}^\infty a_n$ converges, then the product $\prod_{n=1}^\infty (1+a_n)$ also converges.
Moreover, if $\sum_{n=1}^\infty a_n < 1$ then\footnote{Let $b = \sum_{n=1}^\infty a_n$. 
We have 
\[   \log \prod_{n=1}^\infty (1+a_n)  = \sum_{n=1}^\infty \log (1+a_n) 
     = \sum_{n=1}^\infty \left(a_n - \frac{a_n^2}{2} + \frac{a_n^3}{3}- \dots\right) 
     \leq  b + \frac{b^2}{2} + \frac{b^3}{3}+ \dots 
     = \log \frac{1}{1-b}.\]}
\[  1\leq \prod_{n=1}^\infty (1+a_n) \leq \frac{1}{1-\left(\sum_{n=1}^\infty a_n\right)}. \]
In particular, if $\sum_{n=1}^\infty a_n < \frac{1}{2}$ then 
\begin{equation}  \label{eq: infinite sum and product}
   0 \leq  \prod_{n=1}^\infty (1+a_n) -1 \leq 2\sum_{n=1}^\infty a_n. 
\end{equation}   
Namely, if $\sum_{n=1}^\infty a_n$ is close to zero then $\prod_{n=1}^\infty (1+a_n)$ is close to one.

We can assume $0<\eta < 1$. We will show that the statement holds for $\delta_2 := \delta_3/4$
where $\delta_3$ is a positive constant introduced in Lemma \ref{lemma: gluing} (gluing lemma).
Let $R$ be a large positive number which will be fixed in (\ref{eq: choice of R in resolution of singularity}) below.
Set $\Lambda = \mathbb{Z}(2R) + \mathbb{Z}(2R\sqrt{-1}) \subset \mathbb{C}$.
This is a lattice of the plane.
We enumerate its elements as $\Lambda = \{p_1, p_2, p_3, \dots\}$.
We have $\mathbb{C} = \bigcup_{n=1}^\infty D_{2R}(p_n)$.
We fix a small positive number $c \in (0, \frac{1}{2})$ independent\footnote{If we choose $c$ for the case of $R=1$, 
then the same value of $c$ works for all $R>0$ by a scale change.} 
of $R$ such that for any square $S\subset \mathbb{C}$ of side length $L\geq R$ 
we have 
\begin{equation} \label{eq: choice of c in resolution of singularity}
   \mathbf{m}\left(S\cap D_{cR}(\Lambda)\right) < \frac{\eta}{10} L^2,
\end{equation}
where $\mathbf{m}$ denotes the two-dimensional Lebesgue measure and $D_{cR}(\Lambda) = \bigcup_{n=1}^\infty D_{cR}(p_n)$.

Various estimates in this proof are based on the following simple fact:
For any $z\in \mathbb{C}$ we have 
\[  \sum_{|z-p_n|\geq cR} \frac{1}{|z-p_n|^3} \leq \frac{\mathrm{const}}{R^3}. \]
The left-hand side is a sum over all natural numbers $n$ with $|z-p_n|\geq cR$, and 
“$\mathrm{const}$” in the right-hand side is a positive constant independent of $z$ and $R$.

Let $R_2$ and $C_4$ be the positive constants introduced in Lemma \ref{lemma: gluing} (gluing lemma).
We choose $R>2+20\delta_3+R_2$ so large that for any $z\in \mathbb{C}$ we have 
\begin{equation}  \label{eq: choice of R in resolution of singularity}
   \sum_{|z-p_n|\geq cR} \frac{C_4}{|z-p_n|^3} < \min\left(\frac{\eta}{100}, \frac{\delta_3}{100R}\right). 
\end{equation}

A point $p_n\in \Lambda$ is said to be good if we have $\norm{df}_{L^\infty(D_R(p_n))} \geq \frac{\delta_3}{4R}$.
Otherwise, we say that $p_n$ is bad.
If $p_n$ is bad then $f\left(D_R(p_n)\right) \subset B_{\delta_3/4}\left(f(p_n)\right)$.
Hence
\begin{equation}   \label{eq: behavior of bad points}
     \forall g\in B_{\delta_3/4}(f, \mathbf{d}_{\mathbb{C}}): \quad  g\left(D_R(p_n)\right) \subset B_{\delta_3/2}\left(f(p_n)\right), \quad
    \text{if $p_n$ is bad}. 
\end{equation}    
We will apply the gluing lemma (Lemma \ref{lemma: gluing}) at each bad point $p_n$.

We will inductively construct maps $\Phi_n\colon B_{\delta_3/4}(f, \mathbf{d}_{\mathbb{C}}) \to \mathcal{H}^N$ $(n\geq 0)$ such that 
\begin{equation}  \label{eq: induction hypothesis in resolution}
   \dfs\left(g(z), \Phi_n(g)(z)\right) \leq  \sum_{k=1}^n \frac{C_4}{|z-p_k|^3}, \quad 
   \left(z\in \mathbb{C}, g\in B_{\delta_3/4}(f, \mathbf{d}_{\mathbb{C}}) \right).
\end{equation} 
We set $\Phi_0(g) = g$.
Suppose $\Phi_n$ has been constructed. 
If $p_{n+1}$ is good then we set $\Phi_{n+1} = \Phi_n$.
If $p_{n+1}$ is bad then we proceed as follows.
From (\ref{eq: induction hypothesis in resolution}) and (\ref{eq: choice of R in resolution of singularity}), 
for any $g\in B_{\delta_3/4}(f, \mathbf{d}_{\mathbb{C}})$ and $z\in D_R(p_{n+1})$ we have 
\[ \dfs\left(g(z), \Phi_n(g)(z)\right) \leq \sum_{k=1}^n \frac{C_4}{|z-p_k|^3}
     \leq  \sum_{|z-p_k| \geq cR} \frac{C_4}{|z-p_k|^3} < \frac{\delta_3}{2}. \]
Here we have used a fact that $|z-p_k|\geq R > cR$ for all $1\leq k \leq n$.
Since $g\left(D_R(p_{n+1})\right)\subset B_{\delta_3/2}\left(f(p_{n+1})\right)$,
we have $\Phi_n(g)\left(D_R(p_{n+1})\right) \subset B_{\delta_3}\left(f(p_{n+1})\right)$.

Now we use the gluing lemma (Lemma \ref{lemma: gluing}) with the parameters $p=p_{n+1}$ and $q = f(p_{n+1})$:
There exists a map 
$\Psi\colon \{h\in \mathcal{H}^N\mid h\left(D_R(p_{n+1})\right) \subset B_{\delta_3}\left(f(p_{n+1})\right)\}\to \mathcal{H}^N$
satisfying the conditions required in Lemma \ref{lemma: gluing}.
Let $g\in B_{\delta_3/4}(f, \mathbf{d}_{\mathbb{C}})$.
We define $\Phi_{n+1}(g) := \Psi\left(\Phi_n(g)\right)$.
This satisfies the following conditions. 

\begin{enumerate}
   \item[(A)]  We have $\dfs\left(\Phi_n(g)(z), \Phi_{n+1}(g)(z)\right) \leq \frac{C_4}{|z-p_{n+1}|^3}$.
             Hence the induction hypothesis (\ref{eq: induction hypothesis in resolution}) is satisfied by $\Phi_{n+1}$ 
             and the inductions process works.
   \item[(B)]  We have $\frac{1}{100} < \norm{d\Phi_{n+1}(g)}_{L^\infty(D_{R/2}(p_{n+1}))} < 1$.
   \item[(C)]  For $|z-p_{n+1}| \geq 1$ we have
   \[  \left|\left|d\Phi_{n+1}(g)\right|(z) - \left|d\Phi_n(g)\right|(z)\right| \leq  \frac{C_4}{|z-p_{n+1}|^3} \left|d\Phi_n(g)\right|(z) + 
        \frac{C_4}{|z-p_{n+1}|^3}. \]
   \item[(D)]  Let $g_1, g_2\in B_{\delta_3/4}(f, \mathbf{d}_{\mathbb{C}})$.  For $|z-p_{n+1}|\leq 1$ we have 
   \begin{equation*}  
      \begin{split}
        &  \dfs\left(\Phi_{n+1}(g_1)(z), \Phi_{n+1}(g_2)(z)\right)  \leq  C_4\,  \dfs\left(\Phi_n(g_1)(z), \Phi_n(g_2)(z)\right),   \\
        &  \dfs\left(\Phi_n(g_1)(z), \Phi_n(g_2)(z)\right)  \leq  C_4 \sup_{|w-p_{n+1}|\leq 2} \dfs\left(\Phi_{n+1}(g_1)(w), \Phi_{n+1}(g_2)(w)\right). 
      \end{split} 
   \end{equation*}    
   For $|z-p_{n+1}|\geq 1$ we have 
   \begin{equation*} 
      \begin{split}
         \dfs\left(\Phi_{n+1}(g_1)(z), \Phi_{n+1}(g_2)(z)\right)  
         \leq  \left(1+\frac{C_4}{|z-p_{n+1}|^3}\right)  \dfs\left(\Phi_n(g_1)(z), \Phi_n(g_2)(z)\right), \\
         \dfs\left(\Phi_n(g_1)(z), \Phi_n(g_2)(z)\right)  \leq 
         \left(1+\frac{C_4}{|z-p_{n+1}|^3}\right)  \dfs\left(\Phi_{n+1}(g_1)(z), \Phi_{n+1}(g_2)(z)\right).
      \end{split}  
   \end{equation*}
\end{enumerate}
Notice that the conditions (A), (C) and (D) trivially holds if $p_{n+1}$ is good.
Recall that $\Phi_{n+1}=\Phi_n$ if $p_{n+1}$ is good.

From the condition (A), for any $g\in B_{\delta_3/4}(f, \mathbf{d}_{\mathbb{C}})$ and any compact 
subset $K\subset \mathbb{C}$ we have 
\[  \sum_{n=0}^\infty \sup_{z\in K} \dfs\left(\Phi_n(g)(z), \Phi_{n+1}(g)(z)\right) < \infty. \]
Hence the sequence $\Phi_n(g)$ $(n\geq 1)$ converges uniformly over every compact subset of the plane $\mathbb{C}$.
So we can define $\Phi(g): = \lim_{n\to \infty} \Phi_n(g)$ as a holomorphic map from $\mathbb{C}$ to $\mathbb{C} P^N$.

From the conditions (B) and (C), we have that for any $z\in \mathbb{C}$
\[  \left|d\Phi(g)\right|(z) \leq  \left(1+\sum_{|z-p_n| >R/2} \frac{C_4}{|z-p_n|^3}\right)
                                            \prod_{|z-p_n|>R/2} \left(1+\frac{C_4}{|z-p_n|^3}\right). \]
From (\ref{eq: infinite sum and product}) and (\ref{eq: choice of R in resolution of singularity}) with $R/2>cR$ (recall $c\in (0,1/2)$),
\[  \left|d\Phi(g)\right|(z)  \leq   \left(1+\frac{\eta}{100}\right)\left(1+\frac{\eta}{50}\right) < 1 + \eta. \]
Therefore $\Phi$ defines a map from $B_{\delta_3/4}(f, \mathbf{d}_{\mathbb{C}})$ to $\mathcal{B}^N_{1+\eta}$.
Similarly, we have $\Phi_n(g) \in \mathcal{B}^N_{1+\eta}$ for all $n\geq 0$ and $g\in B_{\delta_3/4}(f, \mathbf{d}_{\mathbb{C}})$.
Recall that $0<\eta<1$ and hence $\mathcal{B}^N_{1+\eta}\subset \mathcal{B}^N_2$.
Then, from the condition (C), for $|z-p_{n+1}|\geq 1$ and $g\in B_{\delta_3/4}(f, \mathbf{d}_{\mathbb{C}})$
\begin{equation}  \label{eq: distortion of spherical derivative}
    \left|\left|d\Phi_{n+1}(g)\right|(z) - \left|d\Phi_n(g)\right|(z)\right|   \leq   \frac{3C_4}{|z-p_{n+1}|^3}. 
\end{equation}    

We would like to check the condition (1) in the statement.
For each $p_n$ we can find a point $z\in D_R(p_n)$ such that 
\[  \max \left(|df|(z), \left|d\Phi_n(f)\right|(z)\right) \geq \frac{\delta_3}{4R}. \]
Indeed, if $p_n$ is good then we have $\norm{df}_{L^\infty\left(D_R(p_n)\right)} \geq \frac{\delta_3}{4R}$.
If $p_n$ is bad, then we have $\left|d\Phi_n(f)\right|(z) \geq \frac{1}{100} > \frac{\delta_3}{4R}$ for some 
$z\in D_{R/2}(p_n)$ by the condition (B).
(Here recall that $R>20\delta_3$ and hence $\frac{1}{100} > \frac{\delta_3}{4R}$.)
By (\ref{eq: distortion of spherical derivative}) and (\ref{eq: choice of R in resolution of singularity})
\[  \left|d\Phi(f)\right|(z) \geq   \frac{\delta_3}{4R} - \sum_{|z-p_k|\geq R}\frac{3C_4}{|z-p_k|^3}
     \geq   \frac{\delta_3}{4R}   -  \frac{3\delta_3}{100R} > \frac{\delta_3}{10R}. \]
Since $\mathbb{C} = \bigcup_{n=1}^\infty D_{2R}(p_n)$, the holomorphic map $\Phi(f)$ is $R_1$-nondegenerate 
for $R_1:= 3R+ \frac{10R}{\delta_3}$.
This has shown the condition (1).

Next we consider the condition (2).
Let $g_1, g_2 \in B_{\delta_3/4}(f, \mathbf{d}_{\mathbb{C}})$.
For any $z\in \mathbb{C}$, from the condition (D)
\begin{align*}
  \dfs\left(\Phi(g_1)(z), \Phi(g_2)(z)\right)  
  & \leq   C_4 \prod_{|z-p_n|\geq 1} \left(1+\frac{C_4}{|z-p_n|^3}\right) \dfs\left(g_1(z), g_2(z)\right) \\
  & \leq   2C_4 (1+C_4)\dfs\left(g_1(z), g_2(z)\right).
\end{align*}
\begin{align*}
  &  \dfs\left(g_1(z), g_2(z)\right)  \\ &  \leq \prod_{|z-p_n|\geq 1}\left(1+\frac{C_4}{|z-p_n|^3}\right)
  \sup_{|w-z|\leq 3} \left\{C_4 \prod_{|w-p_n|\geq 1} \left(1+\frac{C_4}{|w-p_n|^3}\right) 
  \dfs\left(\Phi(g_1)(w), \Phi(g_2)(w)\right)\right\}  \\
  & \leq  4C_4(1+C_4)^2 \sup_{|w-z|\leq 3}  \dfs\left(\Phi(g_1)(w), \Phi(g_2)(w)\right).
\end{align*}
Hence $\Phi$ satisfies the condition (2) in the statement with the constant 
$C_3 := 4C_4 (1+C_4)^2$.

Finally we consider the condition (3).
For $|z-p_{n+1}|\geq 1$ we have 
\begin{align*}
  & \left|\left|d\Phi_n(f)\right|^2(z) - \left|d\Phi_{n+1}(f)\right|^2(z)\right|  \\
  &  =
  \left(\left|d\Phi_n(f)\right|(z) + \left|d\Phi_{n+1}(f)\right|(z)\right) \left|\left|d\Phi_n(f)\right|(z) - \left|d\Phi_{n+1}(f)\right|(z)\right|  \\
  & \leq  \frac{12C_4}{|z-p_{n+1}|^3}. 
\end{align*}
Here we have used (\ref{eq: distortion of spherical derivative}) and $\left|d\Phi_n(f)\right|, \left|d\Phi_{n+1}(f)\right| \leq  2$.
Hence for $z\in \mathbb{C}\setminus  D_{cR}(\Lambda)$
\[ \left|\left|d\Phi(f)\right|^2(z) - |df|^2(z)\right| \leq  \sum_{n=1}^\infty \frac{12C_4}{|z-p_n|^3} < \frac{\eta}{2}, \]
by (\ref{eq: choice of R in resolution of singularity}).
Let $S\subset \mathbb{C}$ be any square of side length $L\geq R$.
We have 
\begin{align*}
  \left|\int_S \left|d\Phi(f)\right|^2\, dxdy - \int_S |df|^2\, dxdy\right|  \leq & 
   \left|\int_{S\setminus D_{cR}(\Lambda)} \left|d\Phi(f)\right|^2\, dxdy - \int_{S\setminus D_{cR}(\Lambda)} |df|^2\, dxdy\right| \\
   & + \left|\int_{S\cap  D_{cR}(\Lambda)} \left|d\Phi(f)\right|^2\, dxdy - \int_{S\cap  D_{cR}(\Lambda)} |df|^2\, dxdy\right|. 
\end{align*}
The first term in the right-hand side is smaller than $\frac{\eta L^2}{2}$.
The second term is bounded from above by $4\, \mathbf{m}\left(S\cap D_{cR}(\Lambda)\right)$.
This is smaller than $\frac{\eta L^2}{2}$ by the choice of $c$ in (\ref{eq: choice of c in resolution of singularity}).
Therefore we conclude 
\[ \left|\int_S \left|d\Phi(f)\right|^2\, dxdy - \int_S |df|^2\, dxdy\right|  < \eta \, L^2. \]
This shows the condition (3).
\end{proof}

The next corollary was first proved in \cite[Theorem 1.7]{Matsuo--Tsukamoto}.
This will be used in \S \ref{section: construction of extremal measures}.

\begin{corollary} \label{cor: energy density and nondegenerate curve}
For any $\varepsilon>0$ there exists a nondegenerate Brody curve
$g \colon \mathbb{C}\to \mathbb{C}P^N$ satisfying 
$\rho(g) > \rho(\mathbb{C}P^N)-\varepsilon$.
\end{corollary}

\begin{proof}
By the definition of $\rho(\mathbb{C}P^N)$ there exists $f \in \mathcal{B}^N$
with $\rho(f) > \rho(\mathbb{C}P^N) -\varepsilon$.
Let $\eta$ be a small positive number satisfying 
\begin{equation} \label{eq: choice of eta in corollary}
 \frac{\rho(f)}{(1+\eta)^2} - \frac{\eta}{(1+\eta)^2} > \rho(\mathbb{C}P^N) - \varepsilon.
\end{equation}
Applying Proposition \ref{prop: resolution of singularity} to $f$, we can find a nondegenerate 
curve $h \in \mathcal{B}^N_{1+\eta}$ satisfying 
\[   \left|\int_S |dh|^2\, dxdy - \int_S |df|^2\, dxdy\right| < \eta \, \mathbf{m}(S)  \]
for all sufficiently large squares $S\subset \mathbb{C}$.
Define a nondegenerate curve $g \colon \mathbb{C}\to \mathbb{C}P^N$ 
by $g(z) = h\left(\frac{z}{1+\eta}\right)$.
We have 
\[  |dg|(z) = \frac{1}{1+\eta} |dh|\left(\frac{z}{1+\eta}\right) \leq   1. \]
For any $a\in \mathbb{C}$ and $R>0$
\begin{align*}
   \int_{|z-a|\leq R} |dg|^2(z) \, dxdy & = \frac{1}{(1+\eta)^2} \int_{|z-a|\leq R} |dh|^2\left(\frac{z}{1+\eta}\right)\, dxdy \\
    &  = \int_{\left|z-\frac{a}{1+\eta}\right|\leq \frac{R}{1+\eta}} |dh|^2(z)\, dxdy \\
   & \geq  \int_{\left|z-\frac{a}{1+\eta}\right|\leq \frac{R}{1+\eta}} |df|^2(z) \, dxdy
             - \pi\left(\frac{R}{1+\eta}\right)^2 \eta  - O(R).
\end{align*}
Hence 
\[   \rho(g)
      = \lim_{R\to \infty} \frac{1}{\pi R^2} \sup_{a\in \mathbb{C}} \int_{|z-a|\leq R} |dg|^2(z)\, dxdy
                \geq \frac{\rho(f)}{(1+\eta)^2} - \frac{\eta}{(1+\eta)^2} 
                 > \rho(\mathbb{C}P^N) - \varepsilon. \]
Here we have used (\ref{eq: choice of eta in corollary}) in the last step.
\end{proof}

\section{Energy integral of Brody curves}  \label{section: energy integral of Brody curves}

In this section we prepare some miscellaneous results on the energy integral $\int |df|^2\, dxdy$
of Brody curves $f\in \mathcal{B}^N$.
We will use them in \S \ref{subsection: proof of Theorem construction of extremal measure}.

The next lemma tells us that the asymptotic value of the integral $\int |df|^2\, dxdy$
\textit{does not depend so sensitively} on the choice of $f$.

\begin{lemma}  \label{lemma: energy integral is not sensitive}
There exist positive numbers $\delta_4$ and $C_5$ such that if $f, g \in \mathcal{B}^N$ satisfy
$\mathbf{d}_{\mathbb{C}}(f, g) \leq \delta_4$ then for any square $S\subset \mathbb{C}$ of 
side length $L\geq 1$ we have
\[   \left|\int_S |df|^2\, dxdy - \int_{S} |dg|^2\, dxdy\right| < C_5 L. \]
\end{lemma}

\begin{proof}
Let $\omega$ be the Fubini--Study K\"{a}hler form of $\mathbb{C}P^N$.
We have the basic identity $f^*\omega = |df|^2\, dxdy$ for holomorphic maps 
$f\colon \mathbb{C} \to \mathbb{C} P^N$.
Let $S\subset \mathbb{C}$ be a square of side length $L\geq 1$.
We denote by $S^c$ the complement of $S$ in the plane $\mathbb{C}$.
Set 
$d(z, S^c) = \inf\{|z-w|\mid w\in S^c\}$ where $|z-w|$ is the Euclidean length.

We choose a smooth function $\varphi\colon \mathbb{C}\to [0,1]$ such that 
$\varphi(z) = 0$ on $S^c$ and that $\varphi(z) = 1$ if $d(z, S^c) \geq \frac{1}{2}$.
We can assume 
\[  \left|\frac{\partial \varphi}{\partial x}\right| + \left|\frac{\partial \varphi}{\partial y}\right|
     \leq \mathrm{const} \]
where the right-hand side is an universal positive constant independent of $S$ and $L$.

We choose $\delta_4 >0$ smaller than the injectivity radius of $\mathbb{C}P^N$.
For $0\leq t \leq 1$ and two points $p, q \in \mathbb{C}P^N$ with $\dfs(p, q) \leq \delta_4$,
we denote by $\gamma(p, q, t)$ the point which divides the minimum geodesic segment
between $p$ and $q$ in the ratio $t:(1-t)$.
Namely $\gamma(p, q, t)$ is the point on the minimum geodesic segment between $p$ and $q$ satisfying 
$\dfs\left(p, \gamma(p, q, t)\right) = t \dfs(p, q)$.
In particular, $\gamma(p, q, 0) = p$ and $\gamma(p, q, 1) = q$.

For $0\leq t\leq 1$ we define a $C^\infty$ (not necessarily holomorphic) map 
$f_t\colon \mathbb{C}\to \mathbb{C}P^N$ by $f_t(z) = \gamma\left(f(z), g(z), t\varphi(z)\right)$.
For $z\in S^c$ we have $f_t(z) = f(z)$ for every $t$.
Hence $\{f_t\}$ provides a homotopy between $f_0$ and $f_1$ relative to $S^c$.
This implies that 
\[  \int_S f_0^*\omega = \int_S f_1^* \omega. \]
We have $f_0(z) = \gamma\left(f(z), g(z), 0\right) = f(z)$ for all $z\in \mathbb{C}$. Hence
\[  \int_S f_0^*\omega = \int_S f^*\omega = \int_S |df|^2\, dxdy. \]
If $d(z, S^c) \geq \frac{1}{2}$ then $f_1(z) = \gamma\left(f(z), g(z), 1\right) = g(z)$.
Hence 
\[  \int_S f_1^*\omega = \int_S g^*\omega + O(L) = \int_S |dg|^2 \, dxdy + O(L). \]
Here the implicit constant in $O(L)$ is universal.
Therefore we have 
\[  \left|\int_S |df|^2\, dxdy - \int_{S} |dg|^2\, dxdy\right|  \leq  \mathrm{Const}\cdot  L. \]
\end{proof}

Next we review the Ornstein--Weiss lemma.
This lemma (also known as Ornstein--Weiss' quasi-tiling machinery) was first developed by Ornstein--Weiss \cite{Ornstein--Weiss} in the context of 
amenable group actions. 
Here we state it only for $\mathbb{C}$-actions.

For $r>0$ and $A\subset \mathbb{C}$, we define $\partial_r A$ as the set of $z\in \mathbb{C}$ 
for which the disk $D_r(z)$ intersects both with $A$ and the complement $A^c$.
A sequence of bounded Borel subsets $A_n \subset \mathbb{C}$ $(n\geq 1)$ with $\mathbf{m}(A_n)>0$ is said to be 
a \textbf{F{\o}lner sequence} if we have 
$\lim_{n\to \infty} \frac{\mathbf{m}(\partial_r A_n)}{\mathbf{m}(A_n)} = 0$ for any $r>0$.
For example, $D_n(0)$ and $[0,n]^2$ $(n\geq 1)$ are both F{\o}lner sequences of $\mathbb{C}$.

The following is the Ornstein--Weiss lemma.
For the proof, see Gromov \cite[p. 336]{Gromov}.

\begin{theorem}[Ornstein--Weiss lemma]  \label{theorem: Ornstein--Weiss lemma}
Let $h$ be a nonnegative function defined on 
the set of bounded Borel subsets of $\mathbb{C}$ such that
\begin{enumerate}
  \item   $h(A\cup B) \leq h(A) + h(B)$ for disjoint $A$ and $B$,
  \item   $h(A) \leq  h(B)$ if $A\subset B$,
  \item   $h(a+A) = h(A)$ for $a\in \mathbb{C}$, where $a+A = \{a+x\mid x\in A\}$.
\end{enumerate}
Then for any F{\o}lner sequence $\{A_n\mid n\geq 1\}$ of $\mathbb{C}$, the limit 
$\lim_{n\to \infty} \frac{h(A_n)}{\mathbf{m}(A_n)}$ exists and its value is independent of the choice of 
the F{\o}lner sequences.
\end{theorem}

Recall that the energy density of Brody curves $f\in \mathcal{B}^N$ was defined by 
\[  \rho(f) = \lim_{R\to \infty} \left(\sup_{a\in \mathbb{C}} \frac{1}{\pi R^2} \int_{|z-a|\leq R} |df|^2(z)\, dxdy\right). \]

\begin{corollary}  \label{corollary: reformulation of energy density}
For any $f\in \mathcal{B}^N$ we have 
\[  \rho(f) = \lim_{L\to \infty} \left(\sup_{a\in \mathbb{C}} \frac{1}{L^2} \int_{a+[0, L]^2}|df|^2(z)\, dxdy\right). \]
Here $a+[0,L]^2 = \{a+x+y\sqrt{-1} \mid 0\leq x, y\leq L\}$.
\end{corollary}

\begin{proof}
Apply the Ornstein--Weiss lemma to $h(A) := \sup_{a\in \mathbb{C}} \int_{a+A} |df|^2(z) \, dxdy$.
\end{proof}

Finally we explain that $\rho(\mathbb{C}P^N) = \sup_{f\in \mathcal{B}^N} \rho(f)$ can be expressed 
in terms of invariant probability measures.
We mentioned this result in (\ref{eq: energy density and invariant measures}).
This will not be used in the proofs of the main theorems.
But it is a basic result which initially motivated the research of the paper.

Recall that the potential function $\psi\colon \mathcal{B}^N\to \mathbb{R}$ was defined by 
$\psi(f) = 2(N+1)|df|^2(0)$.

\begin{proposition}  \label{prop: energy density and invariant measures}
\[ 2(N+1) \rho(\mathbb{C}P^N) = \sup_{\mu \in \mathscr{M}^T\left(\mathcal{B}^N\right)} \int_{\mathcal{B}^N}\psi \, d\mu. \]
\end{proposition}

\begin{proof}
For simplicity of the notation, we set $\varphi(f) = |df|^2(0)$ for $f\in \mathcal{B}^N$.
We will prove $\rho(\mathbb{C}P^N) = \sup_{\mu\in \mathscr{M}^T\left(\mathcal{B}^N\right)} \int_{\mathcal{B}^N} \varphi\, d\mu$.

First we show that $\int_{\mathcal{B}^N} \varphi \, d\mu \leq  \rho(\mathbb{C}P^N)$ for every $\mu \in \mathscr{M}^T\left(\mathcal{B}^N\right)$.
This step is close to the argument of Remark \ref{remark: ergodic theorem and characteristic function}.
Let $\mu \in \mathscr{M}^T\left(\mathcal{B}^N\right)$.
By the pointwise ergodic theorem \cite[Theorem 8.19]{Einsiedler--Ward}, 
the limit 
\[  \lim_{R\to \infty} \frac{1}{\pi R^2} \int_{|u|\leq R} \varphi(T^u f)\, d\mathbf{m}(u)  \]
exists for $\mu$-almost every $f\in \mathcal{B}^N$, and its expected value is equal to $\int_{\mathcal{B}^N}\varphi\, d\mu$.
Namely we have 
\[  \int_{\mathcal{B}^N} \left(\lim_{R\to \infty} \frac{1}{\pi R^2} \int_{|z|\leq R} |df|^2(z)\, dxdy\right) d\mu(f)  = \int_{\mathcal{B}^N}\varphi\, d\mu. \]
We have 
\[    \lim_{R\to \infty} \frac{1}{\pi R^2} \int_{|z|\leq R} |df|^2(z)\, dxdy \leq \rho(f) \leq \rho(\mathbb{C}P^N).  \]
Hence $\int_{\mathcal{B}^N} \varphi \, d\mu \leq  \rho(\mathbb{C}P^N)$.

Let $f\in \mathcal{B}^N$. Next we show that there exists $\mu \in \mathscr{M}^T\left(\mathcal{B}^N\right)$ satisfying 
$\int_{\mathcal{B}^N} \varphi\, d\mu = \rho(f)$.
By the definition, there exist sequences of points $a_n\in \mathbb{C}$ and positive numbers $R_n$ with $R_n\to \infty$
such that 
\[  \rho(f) = \lim_{n\to \infty} \frac{1}{\pi R_n^2} \int_{D_{R_n}(a_n)} |df|^2\, dxdy. \]
We define a (not necessarily invariant) Borel probability measure $\mu_n$ on $\mathcal{B}^N$ by
\[  \mu_n = \frac{1}{\pi R_n^2} \int_{D_{R_n}(a_n)} T^u_* \delta_{f} \, d\mathbf{m}(u). \]
Here $\delta_f$ is the delta measure at $f$.
By choosing a subsequence (also denoted by $\mu_n$), 
we can assume that $\mu_n$ converges to a measure $\mu \in \mathscr{M}^T\left(\mathcal{B}^N\right)$
in the weak$^*$ topology.  
We have
\[  \int_{\mathcal{B}^N} \varphi\, d\mu = \lim_{n\to \infty} \int_{\mathcal{B}^N} \varphi\, d\mu_n 
     = \lim_{n\to \infty} \frac{1}{\pi R_n^2} \int_{D_{R_n}(a_n)} |df|^2\, dxdy = \rho(f). \]
\end{proof}

\section{General method to bound the rate distortion dimension from below}
\label{section: general method to bound the rate distortion dimension from below}

In the proof of Theorem \ref{theorem: construction of extremal measure} we will need a method to construct invariant 
probability measures with a good lower bound on rate distortion dimension.
The paper \cite[\S 6]{Tsukamoto_R^d} developed a general method for that purpose
in the course of proving the variational principle for mean dimension with potential 
(Theorem \ref{theorem: variational principle for mean dimension}).
However the formulation of \cite{Tsukamoto_R^d} is not suitable for the application to Brody curves.
The purpose of this section is to reformulate various techniques in \cite{Tsukamoto_R^d} more suitable for the application.

This section is rather technical.
The main result of this section is Proposition \ref{prop: general method to construct invariant measures in L^infty metric}.
We will use it in \S \ref{section: construction of extremal measures}.
We recommend readers to read only Definition \ref{definition: tame growth of covering numbers} and
the statement (not the proof) of Proposition \ref{prop: general method to construct invariant measures in L^infty metric}
and skip to \S \ref{section: construction of extremal measures} at the first reading.

\subsection{Tame growth of covering numbers and $\mathbb{Z}^k$-actions}
\label{subsection: tame growth of covering numbers and Z^k-actions}

The next definition was first introduced in \cite[Condition 3]{Lindenstrauss--Tsukamoto_rate_distortion}.

\begin{definition}  \label{definition: tame growth of covering numbers}
Let $(\mathcal{X}, \mathbf{d})$ be a compact metric space.
It is said to have the \textbf{tame growth of covering numbers} if for any $\delta>0$ we have
\begin{equation} \label{eq: tame growth of covering numbers}
    \lim_{\varepsilon \to 0} \varepsilon^\delta \log \#\left(\mathcal{X}, \mathbf{d}, \varepsilon\right) = 0. 
\end{equation}    
Here recall that $\#\left(\mathcal{X}, \mathbf{d}, \varepsilon\right)$ is the minimum number $n$ such that 
there exists an open covering $\mathcal{X} = U_1\cup \dots \cup U_n$ with $\diam U_i < \varepsilon$ for all $1\leq i \leq n$.
\end{definition}

The condition (\ref{eq: tame growth of covering numbers}) is rather mild.
Every compact subset of the Euclidean space satisfies it
with respect to the Euclidean metric.
A metric $\mathbf{d}$ on the infinite dimensional cube $[0, 1]^{\mathbb{Z}}$ defined by 
$\mathbf{d}(x, y) = \sum_{n\in \mathbb{Z}} 2^{-|n|} |x_n-y_n|$ has the tame growth of covering numbers.
It is known that every compact metrizable space admits a metric having the tame growth of covering numbers
\cite[Lemma 4]{Lindenstrauss--Tsukamoto_rate_distortion}.
We will see later that the metric $\mathbf{d}(f, g) = \max_{z\in [0,1]^2} \dfs\left(f(z), g(z)\right)$ on the space of Brody curves 
$\mathcal{B}^N$ also has the tame growth of covering numbers
(Proposition \ref{prop: Brody curves and tame growth of covering numbers}).

Throughout this paper, we have been studying continuous actions of $\mathbb{R}^k$ (especially, actions of $\mathbb{C} = \mathbb{R}^2$).
For developing a method to construct appropriate invariant measures, 
we also need to consider actions of the discrete subgroup $\mathbb{Z}^k$.
So we prepare some notations about them.
For a natural number $L$ we denote $[L] = \{0,1,2,\dots, L-1\}^k$.
Let $(\mathcal{X}, \mathbf{d})$ be a compact metric space, and let 
$T\colon \mathbb{Z}^k\times \mathcal{X}\to \mathcal{X}$ be a continuous action of the group $\mathbb{Z}^k$.
For a finite (nonempty) subset $A\subset \mathbb{Z}^k$, 
we define metrics $\mathbf{d}^{\mathbb{Z}}_A$ and $\bar{\mathbf{d}}^{\mathbb{Z}}_A$ on $\mathcal{X}$ by 
\[  \mathbf{d}^{\mathbb{Z}}_A(x, y) = \max_{u\in A} \mathbf{d}\left(T^u x, T^u y\right), \quad 
     \bar{\mathbf{d}}^{\mathbb{Z}}_A(x,y) = \frac{1}{|A|} \sum_{u\in A} \mathbf{d}\left(T^u x, T^u y\right). \]
Here “$\mathbb{Z}$” indicates that these are defined for $\mathbb{Z}^k$-actions. 
For a natural number $L$, we also denote $\mathbf{d}^{\mathbb{Z}}_L = \mathbf{d}^{\mathbb{Z}}_{[L]}$
and $\bar{\mathbf{d}}^{\mathbb{Z}}_L = \bar{\mathbf{d}}^{\mathbb{Z}}_{[L]}$.

\begin{lemma}  \label{lemma: tame growth of covering numbers and Z^k actions}
Let $0<\eta<1$ and let $(\mathcal{X}, \mathbf{d})$ be a compact metric space having the tame growth of covering numbers.
There exists a positive number $\varepsilon_0 = \varepsilon_0(\eta, \mathcal{X}, \mathbf{d})$ for which the following statement 
holds. Let $\delta, s, C$ be positive numbers and $L$ a natural number.
Let $T\colon \mathbb{Z}^k\times \mathcal{X}\to \mathcal{X}$ be a continuous action of $\mathbb{Z}^k$ on $\mathcal{X}$.
Suppose a Borel probability measure $\nu$ on $\mathcal{X}$ satisfies 
\begin{equation}  \label{eq: power law estimate in L^infty metric}
    \nu(A) \leq \left(C\cdot \diam\left(A, \mathbf{d}^{\mathbb{Z}}_L\right)\right)^{s L^k} 
\end{equation}    
for any Borel subset $A \subset \mathcal{X}$ with $\diam\left(A, \mathbf{d}^{\mathbb{Z}}_L\right) < \delta$.
Then we also have 
\begin{equation} \label{eq: power law estimate in L^1 metric}
   \nu(A) \leq \left\{4^{\frac{1}{s(1-\eta)}} (2C)^{\frac{1}{1-\eta}} 
   \diam\left(A, \bar{\mathbf{d}}^{\mathbb{Z}}_L\right)\right\}^{(1-\eta)s L^k} 
\end{equation}   
for any Borel subset $A\subset \mathcal{X}$ with 
$\diam\left(A, \bar{\mathbf{d}}^{\mathbb{Z}}_L\right) < \min\left(\varepsilon_0, (\delta/2)^{\frac{1}{1-\eta}}\right)$.
\end{lemma}

This is a rather technical statement.
Roughly speaking, it says that, under the assumption of tame growth of covering numbers, the two metrics
$\mathbf{d}^{\mathbb{Z}}_L$ and $\bar{\mathbf{d}}^{\mathbb{Z}}_L$ behave similarly in the sense that the estimate 
$\nu(A) \leq \left(C\cdot \diam\left(A, \mathbf{d}^{\mathbb{Z}}_L\right)\right)^{s L^k}$ implies a similar estimate 
$\nu(A) \leq \left(C^\prime \cdot \diam\left(A, \bar{\mathbf{d}}^{\mathbb{Z}}_L\right)\right)^{(1-\eta)sL^k}$.
(The precise form of the constant $C^\prime$ is not relevant to the argument in the sequel.)
Conversely, the estimate $\nu(A) \leq \left(C^\prime \cdot \diam\left(A, \bar{\mathbf{d}}^{\mathbb{Z}}_L\right)\right)^{(1-\eta)sL^k}$
obviously implies $\nu(A) \leq \left(C^\prime \cdot \diam\left(A, \mathbf{d}^{\mathbb{Z}}_L\right)\right)^{(1-\eta)sL^k}$
because we have $\bar{\mathbf{d}}^{\mathbb{Z}}_L \leq  \mathbf{d}^{\mathbb{Z}}_L$.

\begin{proof}
Set $M(\varepsilon) = \log \#\left(\mathcal{X}, \mathbf{d}, \varepsilon\right)$.
From the condition of tame growth of covering numbers, we can find $\varepsilon_0 \in (0, 1)$ such that 
for any $0<\varepsilon \leq \varepsilon_0$ we have 
\begin{equation} \label{eq: choice of varepsilon_0 tame growth of covering numbers}
  M(\varepsilon)^{\varepsilon^\eta} < 2.
\end{equation}

Suppose a Borel subset $A\subset \mathcal{X}$ satisfies 
$\diam\left(A, \bar{\mathbf{d}}^{\mathbb{Z}}_L\right) < \min\left(\varepsilon_0, (\delta/2)^{\frac{1}{1-\eta}}\right)$.
We will show (\ref{eq: power law estimate in L^1 metric}).
Set $\tau = \diam\left(A, \bar{\mathbf{d}}^{\mathbb{Z}}_L\right)$. 
We can assume $\tau>0$.
There is an open covering $\mathcal{X} = U_1\cup \dots \cup U_M$ with $M = M(\tau)$ and 
$\diam\left(U_m, \mathbf{d}\right) < \tau$ for all $1\leq m \leq M$.
From (\ref{eq: choice of varepsilon_0 tame growth of covering numbers}) we have $M^{\tau^\eta} < 2$.

Pick and fix a point $a\in A$.
Let $x\in \mathcal{X}$. We have $\bar{\mathbf{d}}^{\mathbb{Z}}_L(x, a) \leq \tau$ and hence 
\[ \left|\{u\in [L] \mid \mathbf{d}\left(T^u x, T^u a\right) > \tau^{1-\eta}\} \right| < \tau^\eta L^k. \] 
Hence 
\begin{equation} \label{eq: cover by partial Bowen balls}
    A\subset \bigcup_{\substack{E\subset [L] \\  |E| < \tau^\eta L^k}}   
     B_{\tau^{1-\eta}} \left(a, \mathbf{d}^{\mathbb{Z}}_{[L]\setminus E} \right). 
\end{equation}     
Here $E$ runs over all subsets of $[L] = \{0,1,2,\dots, L-1\}^k$ of cardinality smaller than $\tau^\eta L^k$, and 
$B_{\tau^{1-\eta}} \left(a, \mathbf{d}_{[L]\setminus E} \right)$ is the ball of radius $\tau^{1-\eta}$ centered at $a$ 
with respect to $\mathbf{d}^{\mathbb{Z}}_{[L]\setminus E}(x, y) = \max_{u\in [L] \setminus E}\mathbf{d}(T^u x, T^u y)$.

Let $E = \{e_1, e_2, \dots, e_r\}\subset [L]$ with $r < \tau^\eta L^k$.
We have 
\[  \mathcal{X} = \bigcup_{1\leq m_1, \dots, m_r\leq M} \left(T^{-e_1}U_{m_1} \cap T^{-e_2}U_{m_2}\cap \dots \cap T^{-e_r}U_{m_r}\right). \]
Hence
\[  B_{\tau^{1-\eta}} \left(a, \mathbf{d}^{\mathbb{Z}}_{[L]\setminus E} \right) = \bigcup_{1\leq m_1, \dots, m_r\leq M}
     \left(B_{\tau^{1-\eta}} \left(a, \mathbf{d}^{\mathbb{Z}}_{[L]\setminus E} \right)\cap T^{-e_1}U_{m_1}\cap T^{-e_2}U_{m_2}\cap \dots \cap 
     T^{-e_r}U_{m_r}\right). \]
The diameter of 
$B_{\tau^{1-\eta}} \left(a, \mathbf{d}^{\mathbb{Z}}_{[L]\setminus E} \right) \cap T^{-e_1}U_{m_1}\cap T^{-e_2}U_{m_2}\cap \dots \cap 
     T^{-e_r}U_{m_r}$
with respect to the metric $\mathbf{d}^{\mathbb{Z}}_L$ is at most $2\tau^{1-\eta} < \delta$.
By applying the estimate (\ref{eq: power law estimate in L^infty metric}) to the sets
$B_{\tau^{1-\eta}} \left(a, \mathbf{d}^{\mathbb{Z}}_{[L]\setminus E} \right) \cap T^{-e_1}U_{m_1}\cap T^{-e_2}U_{m_2}\cap \dots \cap 
     T^{-e_r}U_{m_r}$,
we get 
\[  \nu\left(B_{\tau^{1-\eta}} \left(a, \mathbf{d}^{\mathbb{Z}}_{[L]\setminus E} \right)\right) \leq M^r \left(2C\tau^{1-\eta}\right)^{sL^k}. \]
Since $r< \tau^\eta L^k$ and $M^{\tau^\eta} < 2$, we have 
\[  \nu\left(B_{\tau^{1-\eta}} \left(a, \mathbf{d}^{\mathbb{Z}}_{[L]\setminus E} \right)\right) 
     \leq    2^{L^k} \left(2C\tau^{1-\eta}\right)^{sL^k}. \]
The number of the choices of $E = \{e_1, \dots, e_r\}\subset [L]$ is at most $2^{L^k}$
Therefore, by (\ref{eq: cover by partial Bowen balls})
\[ \nu(A) \leq  2^{L^k} \cdot 2^{L^k} \left(2C\tau^{1-\eta}\right)^{sL^k} = 
   \left\{4^{\frac{1}{s(1-\eta)}} (2C)^{\frac{1}{1-\eta}} \tau\right\}^{(1-\eta)s L^k}.   \]
\end{proof}

\subsection{Constructing invariant measures with a lower bound on rate distortion dimension}
\label{subsection: constructing invariant measures with a lower bound on rate distortion dimension}

We need the following elementary lemma.
This was proved in \cite[Appendix]{Lindenstrauss--Tsukamoto_rate_distortion}.

\begin{lemma}  \label{lemma: elementary optimal transport}
Let $A$ be a finite set and let $\mu_n$ $(n\geq 1)$ be a sequence of probability measures on $A$.
Suppose that $\mu_n$ converges to a probability measure $\mu$ in the weak$^*$ topology, i.e.
$\mu_n(a) \to \mu(a)$ as $n\to \infty$ for every $a\in A$.
Then there exists a sequence of probability measures $\pi_n$ on $A\times A$ such that 
  \begin{itemize}
   \item   $\pi_n$ is a coupling between $\mu_n$ and $\mu$, 
              namely the first and second marginals of $\pi_n$ are given by $\mu_n$ and $\mu$
              respectively,
   \item   $\pi_n$ converges to $(\mathrm{id}\times \mathrm{id})_*\mu$ as $n\to \infty$ in the weak$^*$ topology, namely
      \[  \pi_n(a, b) \to \begin{cases}  \mu(a) & (\text{if $a=b$}) \\ 0  & (\text{if $a\neq b$}) \end{cases}. \]
  \end{itemize}
\end{lemma}

Let $(\mathcal{X}, \mathbf{d})$ be a compact metric space with a continuous $\mathbb{R}^k$-action 
$T\colon  \mathbb{R}^k\times \mathcal{X}\to \mathcal{X}$.
For a positive number $L$ we define a metric $\bar{\mathbf{d}}_L$ on $\mathcal{X}$ by 
\[  \bar{\mathbf{d}}_L(x, y) = \frac{1}{L^k} \int_{[0,L]^k} \mathbf{d}\left(T^u x, T^u y\right) \, d\mathbf{m}(u). \]
Here $\mathbf{m}$ denotes the Lebesgue measure on $\mathbb{R}^k$.

\begin{proposition}   \label{prop: general method to construct measures}
Let $(\mathcal{X}, \mathbf{d})$ be a compact metric space with a continuous action 
$T\colon \mathbb{R}^k\times \mathcal{X}\to \mathcal{X}$.
Let $s$ be a positive number and $L_n$ $(n\geq 1)$ a sequence of positive numbers with 
$\lim_{n\to \infty} L_n = \infty$.
Let $\nu_n$ $(n\geq 1)$ be a sequence of (not necessarily invariant) 
Borel probability measures on $\mathcal{X}$. 
Define 
\[  \mu_n  = \frac{1}{L_n^k} \int_{[0,L_n]^k} T^u_*\nu_n\, d\mathbf{m}(u).  \]
We assume the following two conditions.
  \begin{enumerate}
      \item  There exist positive numbers $\delta$ and $C$ such that we have
        \begin{equation}  \label{eq: power law condition}
          \nu_n(A) \leq \left(C\cdot \diam\left(A, \bar{\mathbf{d}}_{L_n}\right)\right)^{s L_n^k}.
        \end{equation}
         for any Borel subset $A\subset \mathcal{X}$ with $\diam\left(A, \bar{\mathbf{d}}_{L_n}\right) < \delta$.
      \item $\mu_n$ converges to some $\mu \in \mathscr{M}^T(\mathcal{X})$ as $n\to \infty$ in the weak$^*$ topology.
   \end{enumerate}
Then we have
\[  \lrdim\left(\mathcal{X}, T, \mathbf{d}, \mu\right) \geq  s. \]
\end{proposition}

We note that the condition (2) is a mild assumption.
The space of all Borel probability measures on $\mathcal{X}$ is compact in the weak$^*$ topology.
Hence, even if we do not assume (2), we can always take an appropriate subsequence $\{\mu_{n_i}\}$ 
which converges to some Borel probability measure $\mu$. 
From the definition of $\mu_n$, the limit $\mu$ becomes $T$-invariant

\begin{proof}
The following argument was essentially given in the proof of \cite[Theorem 6.3]{Tsukamoto_R^d}.
But the statement of \cite[Theorem 6.3]{Tsukamoto_R^d} is substantially different from this proposition.
(It studies a relation between rate distortion dimension and mean Hausdorff dimension
with potential.)
So we provide a full proof for completeness.
The proof is rather long. Readers may skip it at the first reading.

First we note that if we define a metric $\mathbf{d}^\prime$ on $\mathcal{X}$ by $\mathbf{d}^\prime(x, y) = C \mathbf{d}(x, y)$
then we have $\nu_n(A) \leq \left(\diam\left(A, \overline{\mathbf{d}^\prime}_{L_n}\right)\right)^{s L_n^k}$ for any Borel subset 
$A\subset \mathcal{X}$ with $\diam\left(A, \overline{\mathbf{d}^\prime}_{L_n}\right) < C \delta$.
We have $\lrdim\left(\mathcal{X}, T, \mathbf{d}^\prime, \mu\right) = \lrdim\left(\mathcal{X}, T, \mathbf{d}, \mu\right)$.
Therefore we can assume $C=1$ in (\ref{eq: power law condition}) without loss of generality.
Namely, from the beginning, we assume 
\begin{equation}  \label{eq: power law condition revised}
     \nu_n(A) \leq \left(\diam\left(A, \bar{\mathbf{d}}_{L_n}\right)\right)^{s L_n^k}  
\end{equation}
for any Borel subset 
$A\subset \mathcal{X}$ with $\diam\left(A, \bar{\mathbf{d}}_{L_n}\right) < \delta$.

Let $R(\mathbf{d}, \mu, \varepsilon)$ $(\varepsilon >0)$ be the rate distortion function with respect to the measure $\mu$.
We will prove 
\begin{equation} \label{eq: desired lower estimate on the rate distortion function}
  R(\mathbf{d}, \mu, \varepsilon) \geq  s \log (1/\varepsilon)  - Ks  \quad 
  \text{for sufficiently small $\varepsilon >0$}.
\end{equation}
Here $K$ is the universal positive constant introduced in the Kawabata--Dembo estimate 
(Proposition \ref{prop: Kawabata--Dembo estimate}).
The estimate (\ref{eq: desired lower estimate on the rate distortion function}) implies 
$\lrdim\left(\mathcal{X}, T, \mathbf{d}, \mu\right) = 
\varliminf_{\varepsilon \to 0} \frac{R(\mathbf{d}, \mu, \varepsilon)}{\log (1/\varepsilon)} \geq s$.

Let $\varepsilon$ be a positive number with $2\varepsilon \log (1/\varepsilon) < \delta$.
Let $M$ be a positive number, and let $X$ and $Y$ be random variables satisfying the following two conditions.
  \begin{itemize}
     \item  $X$ takes values in $\mathcal{X}$ and its distribution is $\mu$.
     \item  $Y$ takes values in $L^1\left([0, M)^k, \mathcal{X}\right)$ with 
     \[ \mathbb{E}\left(\frac{1}{M^k}\int_{[0,M)^k} \mathbf{d}\left(T^v X, Y_v\right) d\mathbf{m}(v)\right) < \varepsilon. \]
  \end{itemize}
The aim of the argument below is to prove the lower bound on the mutual information 
\begin{equation} \label{eq: lower bound on mutual information in dynamical Frostman lemma}
   \frac{1}{M^k} I(X;Y) \geq  s \log(1/\varepsilon)  - Ks. 
\end{equation}    
Once this is proved, we will get $\frac{1}{M^k} R\left(\varepsilon, [0, M)^k\right) \geq s \log(1/\varepsilon) - Ks$ and hence
 (\ref{eq: desired lower estimate on the rate distortion function}).
For this purpose, we can assume that the random variable $Y$ takes only finitely many values in 
$L^1\left([0, M)^k, \mathcal{X}\right)$ by Remark \ref{remark: finite random variable in the definition of rate distortion function}.
Namely, we assume that there is a finite subset $\mathcal{Y}\subset  L^1\left([0,M)^k, \mathcal{X}\right)$ such that $Y$ takes 
values only in $\mathcal{Y}$.

We choose a positive number $\tau$ satisfying 
$\mathbb{E}\left(\frac{1}{M^k}\int_{[0,M)^k} \mathbf{d}\left(T^v X, Y_v\right) d\mathbf{m}(v)\right) < \varepsilon -3\tau$.
We take a measurable partition $\mathcal{X} = P_1\cup \dots \cup P_\alpha$ (disjoint union) such that 
$\diam\left(P_i, \bar{\mathbf{d}}_M \right) < \tau$ and $\mu(\partial P_i) = 0$ for all $1\leq i \leq \alpha$.
We pick a point $x_i\in P_i$ for each $i$ and set $A = \{x_1, \dots, x_\alpha\}$.
We define a map $\mathcal{P}\colon \mathcal{X} \to A$ by $\mathcal{P}(P_i) = \{x_i\}$.
We consider the random variable $\mathcal{P}(X)$. 
This takes values in $A$ and satisfies $\bar{\mathbf{d}}_M\left(\mathcal{P}(X), X\right) < \tau$ almost surely.
Hence
\begin{equation}  \label{eq: averaged distance between P(X) and Y}
  \mathbb{E}\left(\frac{1}{M^k} \int_{[0, M)^k} \mathbf{d}\left(T^v \mathcal{P}(X), Y_v\right) d\mathbf{m}(v)\right) < 
  \varepsilon - 2\tau.
\end{equation}
The distribution of $\mathcal{P}(X)$ is given by the push-forward measure $\mathcal{P}_*\mu$.

Let $\mathcal{P}_*\mu_n$ $(n\geq 1)$ be the push-forward measures of $\mu_n$ by $\mathcal{P}$.
They converges to $\mathcal{P}_*\mu$ in the weak$^*$ topology by $\mu(\partial P_i) = 0$.
By Lemma \ref{lemma: elementary optimal transport}, we can take random variables $X_n$ $(n\geq 1)$ such that 
\begin{itemize}
   \item  $X_n$ takes values in $A$ and its distribution is given by $\mathcal{P}_*\mu_n$,
   \item  for each $x_i, x_j\in A$ we have 
   $\mathbb{P}\left(X_n = x_i, \mathcal{P}(X) = x_j\right) \to \delta_{ij} \mathbb{P}\left(\mathcal{P}(X) = x_j\right)$ as $n\to \infty$,
   \item  $X_n$ and $Y$ are conditionally independent given $\mathcal{P}(X)$.
\end{itemize}

Now, for each $x_i \in A$ and $y\in \mathcal{Y}$
\begin{align*}
  &\mathbb{P}\left(X_n= x_i, Y=y\right)  
   = \sum_{j=1}^\alpha \mathbb{P}\left(X_n = x_i, Y=y \mid \mathcal{P}(X) = x_j\right) \mathbb{P}\left(\mathcal{P}(X) = x_j\right) \\
  &  =  \sum_{j=1}^\alpha  
  \mathbb{P}\left(X_n = x_i\mid \mathcal{P}(X) = x_j\right) \mathbb{P}\left(Y=y\mid \mathcal{P}(X) = x_j\right) 
  \mathbb{P}\left(\mathcal{P}(X) = x_j\right) \\
  & = \sum_{j=1}^\alpha  \mathbb{P}\left(X_n=x_i, \mathcal{P}(X) = x_j\right) \mathbb{P}\left(Y=y\mid \mathcal{P}(X) = x_j\right) \\
  & \to  \mathbb{P}\left(\mathcal{P}(X) = x_i\right) \cdot \mathbb{P}\left(Y=y\mid \mathcal{P}(X) = x_i\right)  \quad  (n\to \infty) \\
  & = \mathbb{P}\left(\mathcal{P}(X) = x_i, Y = y\right).
\end{align*}
Therefore by (\ref{eq: averaged distance between P(X) and Y})
\begin{equation}  \label{eq: averaged distance between X_n and Y}
   \mathbb{E}\left(\frac{1}{M^k} \int_{[0, M)^k} \mathbf{d}\left(T^v X_n, Y_v\right) d\mathbf{m}(v)\right) < 
  \varepsilon - 2\tau   \quad \text{for large $n$}
\end{equation}
because the left-hand side converges to 
$\mathbb{E}\left(\frac{1}{M^k} \int_{[0, M)^k} \mathbf{d}\left(T^v \mathcal{P}(X), Y_v\right) d\mathbf{m}(v)\right)$
as $n\to \infty$.

The pair $(X_n, Y)$ takes values in a finite set $A\times \mathcal{Y}$ and its distribution converges to 
$\mathrm{Law}\left(\mathcal{P}(X), Y\right)$ as $n\to \infty$.
Hence by Lemma \ref{lemma: convergence in law and mutual information}   
\begin{equation}  \label{eq: convergence of I(X_n, Y)}
   I\left(\mathcal{P}(X); Y\right)  =  \lim_{n\to \infty} I\left(X_n; Y\right).
\end{equation}
We will estimate $I\left(X_n; Y\right)$ from below.

For $x\in A$ and $y\in L^1\left([0, M)^k, \mathcal{X}\right)$, we define a conditional probability mass function by 
\[   q_n\left(y\mid x\right) = \mathbb{P}\left(Y=y\mid  X_n = x\right). \]
This is zero if $y\not \in \mathcal{Y}$.
For $x\in A$ with $\mathbb{P}(X_n=x)=0$, the function $q_n(\cdot\mid x)$ may be arbitrary.

We define a finite set $\Lambda_n \subset \mathbb{R}^k$ by 
\[  \Lambda_n
   = \left\{(M m_1, M m_2, \dots, M m_k)\middle|\, m_i\in \mathbb{Z}, \, 0\leq m_i \leq \frac{L_n}{M}-2 \text{ for all $1\leq i \leq k$}\right\}. \]
Let $v\in [0, M)^k$. We have 
\[ \bigcup_{\lambda\in \Lambda_n} \left(v+\lambda + [0, M)^k\right) \subset [0, L_n)^k. \]
Here the left-hand side is a disjoint union.
We set 
\[  E_{n,v} = [0, L_n)^k \setminus \bigcup_{\lambda\in \Lambda_n} \left(v+\lambda + [0, M)^k\right). \]
See Figure \ref{figure: squares}. (This figure was given in \cite[Figure 3]{Tsukamoto_R^d}.)

\begin{figure}[h] 
    \centering
    \includegraphics[width=3.0in]{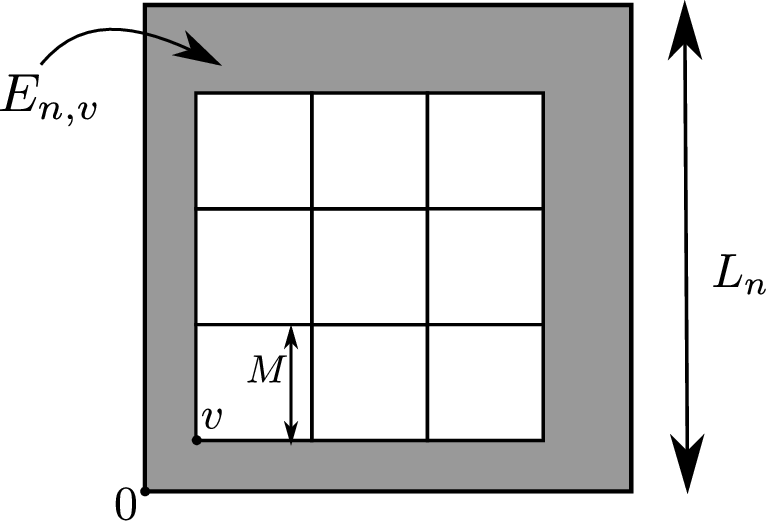}
    \caption{The big square is $[0, L_n)^k$ and small squares are $v+\lambda + [0, M)^k$
    ($\lambda \in \Lambda_n$).  The shadowed region is $E_{n,v}$.}    \label{figure: squares}
\end{figure}

We fix a point $x_0 \in \mathcal{X}$. 
We will denote functions constantly equal to $x_0$ also by $x_0$.
For $x\in \mathcal{X}$ and $f\in L^1\left([0, L_n)^k, \mathcal{X}\right)$ we define a conditional probability mass function 
$\sigma_{n, v}(f\mid x)$ by 
\begin{equation} \label{eq: definition of conditional probability sigma_{n,v}}
     \sigma_{n, v}(f\mid x) = \delta_{x_0}(f|_{E_{n,v}})  \prod_{\lambda\in \Lambda_n} 
     q_n\left(f|_{v+\lambda+[0,M)^k} \mid \mathcal{P}(T^{v+\lambda}x)\right). 
\end{equation}     
For each fixed $x$, this is nonzero only for finitely many $f$.
Here $f|_{E_{n,v}}$ is the restriction of $f$ to $E_{n,v}$, and $\delta_{x_0}$ is the delta measure at 
the constant function $x_0\in L^1\left(E_{n,v}, \mathcal{X}\right)$. 
The function $f|_{v+\lambda+[0,M)^k}$ is the restriction of $f$ to $v+\lambda + [0, M)^k$ and naturally identified with an element of 
$L^1\left([0, M)^k, \mathcal{X}\right)$.

We define a transition probability $\sigma_n$ on $\mathcal{X}\times L^1\left([0, L_n)^k, \mathcal{X}\right)$ by 
\[  \sigma_n(B\mid x) = \frac{1}{M^k}\int_{[0, M)^k} \sigma_{n, v}(B\mid x) d\mathbf{m}(v), \quad 
     \left(x\in \mathcal{X}, B\subset L^1\left([0, L_n)^k, \mathcal{X}\right)\right). \]
Here $\sigma_{n, v}(B\mid x) = \sum_{f\in B} \sigma_{n, v}(f\mid x)$.
We take random variables $Z_n$ and $W_n$ such that 
  \begin{itemize}
    \item  $Z_n$ takes values in $\mathcal{X}$ and its distribution is given by $\nu_n$,
    \item  $W_n$ takes values in $L^1\left([0, L_n)^k, \mathcal{X}\right)$ with 
    \[  \mathbb{P}\left(W_n \in B \mid Z_n = x\right) = \sigma_n(B\mid x), \quad    
     \left(x\in \mathcal{X}, B\subset L^1\left([0, L_n)^k, \mathcal{X}\right)\right). \]
  \end{itemize}

\begin{claim}  \label{claim: distance between Z and W}
  If $n$ is sufficiently large then we have  
\[   \mathbb{E}\left(\frac{1}{L_n^k}\int_{[0, L_n)^k} \mathbf{d}\left(T^u Z_n, (W_n)_u\right) d\mathbf{m}(u)\right) < \varepsilon. \]
\end{claim}

\begin{proof}
For each $v\in [0, M)^k$, we take a random variable $W_n(v)$ such that $W_n(v)$ takes values in $L^1\left([0, L_n)^k, \mathcal{X}\right)$
with $\mathbb{P}\left(W_n(v)=f \mid Z_n=x\right)  =  \sigma_{n, v}(f |x)$ for $x\in \mathcal{X}$ 
and $f\in  L^1\left([0, L_n)^k, \mathcal{X}\right)$.
Then 
\[  \mathbb{E}\left(\int_{[0, L_n)^k} \mathbf{d}\left(T^u Z_n, (W_n)_u\right) d\mathbf{m}(u)\right)
  = \frac{1}{M^k}\int_{[0, M)^k} \mathbb{E}\left(\int_{[0, L_n)^k} \mathbf{d}\left(T^u Z_n, W_n(v)_u\right) d\mathbf{m}(u)\right) d\mathbf{m}(v). \]

Let $v\in [0, M)^k$. 
Since $[0, L_n)^k = E_{n,v} \cup \bigcup_{\lambda\in \Lambda_n} (v+\lambda+[0,M)^k)$ (disjoint union) and 
$\mathbf{m}(E_{n,v}) \leq \mathrm{const} \cdot L_n^{k-1}$,
\begin{align*}
  & \mathbb{E}\left(\int_{[0, L_n)^k} \mathbf{d}\left(T^u Z_n, W_n(v)_u\right) d\mathbf{m}(u)\right) \\
  & = \mathbb{E}\left(\int_{E_{n,v}} \mathbf{d}\left(T^u Z_n, W_n(v)_u\right) d\mathbf{m}(u)\right) +
       \sum_{\lambda\in \Lambda_n} \mathbb{E}\left(\int_{v+\lambda+[0,M)^k} \mathbf{d}\left(T^u Z_n, W_n(v)_u\right) d\mathbf{m}(u)\right) \\
  & \leq C^\prime L_n^{k-1} +  
  \sum_{\lambda\in \Lambda_n} \mathbb{E}\left(\int_{v+\lambda+[0,M)^k} \mathbf{d}\left(T^u Z_n, W_n(v)_u\right) d\mathbf{m}(u)\right).    
\end{align*}
Here $C^\prime$ is a positive constant independent of $v$ and $n$.
Note that $\bar{\mathbf{d}}_M\left(x, \mathcal{P}(x)\right) < \tau$ for all $x\in \mathcal{X}$.
Hence 
\begin{align*}
    \int_{v+\lambda+[0,M)^k} \mathbf{d}\left(T^u Z_n, W_n(v)_u\right) d\mathbf{m}(u) 
   & = \int_{[0,M)^k} \mathbf{d}\left(T^u T^{v+\lambda} Z_n, W_n(v)_{v+\lambda+u}\right) d\mathbf{m}(u) \\
   &\leq  \underbrace{\int_{[0,M)^k} \mathbf{d}\left(T^u T^{v+\lambda} Z_n, 
     T^u \mathcal{P}\left(T^{v+\lambda} Z_n  \right)\right) d\mathbf{m}(u)}_{= M^k 
     \bar{\mathbf{d}}_M\left(T^{v+\lambda}Z_n, \mathcal{P}\left(T^{v+\lambda} Z_n  \right)\right) < M^k \tau} \\
    & + \int_{[0,M)^k} \mathbf{d}\left(T^u \mathcal{P}\left(T^{v+\lambda} Z_n  \right) , W_n(v)_{v+\lambda+u}\right) d\mathbf{m}(u).
\end{align*} 
Therefore 
\begin{equation} \label{eq: averaged distortion between Z and W(v)}
  \begin{split}
  & \mathbb{E}\left(\int_{v+\lambda+[0,M)^k} \mathbf{d}\left(T^u Z_n, W_n(v)_u\right) d\mathbf{m}(u)\right) \\
  & \leq  M^k \tau + 
  \mathbb{E}\left(\int_{[0,M)^k} \mathbf{d}\left(T^u \mathcal{P}\left(T^{v+\lambda} Z_n \right) , W_n(v)_{v+\lambda+u}\right) d\mathbf{m}(u)\right) \\
  & = M^k \tau + \sum_{f\in \mathcal{Y}} \int_{[0,M)^k}\left(\int_{\mathcal{X}}\mathbf{d}\left(T^u x, f_u\right)
        q_n(f|x) d\left(\mathcal{P}_*T^{v+\lambda}_*\nu_n(x)\right)\right) d\mathbf{m}(u).
  \end{split}      
\end{equation} 
In the last equality, we have used that the distribution of 
$\left(\mathcal{P} \left(T^{v+\lambda} Z_n  \right), W_n(v)|_{v+\lambda+[0,M)^k}\right)$
is given by $\mathcal{P}_*T^{v+\lambda}_*\nu_n(x) q_n(f|x)$.
This follows from the definition of $\sigma_{n, v}$ in (\ref{eq: definition of conditional probability sigma_{n,v}}).
 
We sum up the estimate (\ref{eq: averaged distortion between Z and W(v)}) over $\lambda \in \Lambda_n$ and use 
$M^k |\Lambda_n|\leq L_n^k$. Then 
\begin{align*}
  & \mathbb{E}\left(\int_{[0, L_n)^k} \mathbf{d}\left(T^u Z_n, W_n(v)_u\right) d\mathbf{m}(u)\right)  \\
  & \leq  C^\prime L_n^{k-1} + \tau L_n^k  
  +  \sum_{f\in \mathcal{Y}}  
      \int_{[0,M)^k}\left\{\int_{\mathcal{X}}\mathbf{d}\left(T^u x, f_u\right)
        q_n(f|x) d\left(\sum_{\lambda\in \Lambda_n} \mathcal{P}_*T^{v+\lambda}_*\nu_n(x)\right)\right\} d\mathbf{m}(u).
\end{align*}
Integrating this over $v\in [0, M)^k$, we get 
\begin{align*}
  &  \int_{[0, M)^k} \mathbb{E}\left(\int_{[0, L_n)^k} \mathbf{d}\left(T^u Z_n, W_n(v)_u\right) d\mathbf{m}(u)\right) d\mathbf{m}(v) \\
  & \leq  C^\prime L_n^{k-1} M^k + \tau L_n^k M^k \\
   & +  \sum_{f\in \mathcal{Y}}  
      \int_{[0,M)^k}\left\{\int_{\mathcal{X}}\mathbf{d}\left(T^u x, f_u\right)
        q_n(f|x) d\left(\int_{[0, M)^k} \left(\sum_{\lambda\in \Lambda_n} 
        \mathcal{P}_*T^{v+\lambda}_*\nu_n(x)\right) d\mathbf{m}(v)
        \right) \right\} d\mathbf{m}(u).
\end{align*}
We have 
\begin{align*}
    \int_{[0, M)^k} \left(\sum_{\lambda\in \Lambda_n} 
        \mathcal{P}_*T^{v+\lambda}_*\nu_n\right) d\mathbf{m}(v)
    &=  \sum_{\lambda\in \Lambda_n} \int_{\lambda+[0, M)^k} \mathcal{P}_*T^v_*\nu_n d\mathbf{m}(v) \\
    & =  \int_{\bigcup_{\lambda\in \Lambda_n} \left(\lambda+[0, M)^k\right)} \mathcal{P}_*T^v_*\nu_n d\mathbf{m}(v) 
\end{align*}
Since $\bigcup_{\lambda\in \Lambda_n} \left(\lambda+[0, M)^k\right)$ is contained in $[0, L_n)^k$, 
\begin{align*}
     &  \int_{\mathcal{X}}\mathbf{d}\left(T^u x, f_u\right) q_n(f|x) 
       d\left(\int_{[0, M)^k} \left(\sum_{\lambda\in \Lambda} 
        \mathcal{P}_*T^{v+\lambda}_*\nu_n(x)\right) d\mathbf{m}(v)
        \right) \\
      &  \leq  \int_{\mathcal{X}}\mathbf{d}\left(T^u x, f_u\right) q_n(f|x) 
      d\underbrace{\left(\int_{[0, L_n)^k} \mathcal{P}_*T^v_*\nu_n(x)\, d\mathbf{m}(v)\right)}_{L_n^k\, \mathcal{P}_*\mu_n(x)} \\
      & = L_n^k \int_{\mathcal{X}}\mathbf{d}\left(T^u x, f_u\right) q_n(f|x) d\left(\mathcal{P}_*\mu_n(x)\right).
\end{align*}        
Thus 
\begin{align*}
  & \int_{[0, M)^k} \mathbb{E}\left(\int_{[0, L_n)^k} \mathbf{d}\left(T^u Z_n, W_n(v)_u\right) d\mathbf{m}(u)\right) d\mathbf{m}(v) \\
  &  \leq  C^\prime L_n^{k-1} M^k + \tau L_n^k M^k 
    + L_n^k  \sum_{f\in \mathcal{Y}} 
      \int_{[0,M)^k} \left\{\int_{\mathcal{X}}\mathbf{d}\left(T^u x, f_u\right) q_n(f|x) 
      d\left(\mathcal{P}_*\mu_n(x)\right)\right\} d\mathbf{m}(u) \\
  & =      C^\prime L_n^{k-1} M^k + \tau L_n^k M^k 
    + L_n^k \sum_{f\in \mathcal{Y}} \int_{\mathcal{X}} \left(\int_{[0,M)^k} \mathbf{d}\left(T^u x, f_u\right) d\mathbf{m}(u)\right)
        q_n(f|x) \, d\left(\mathcal{P}_*\mu_n(x)\right) \\
  & =     C^\prime L_n^{k-1} M^k + \tau L_n^k M^k 
    + L_n^k \mathbb{E}\left(\int_{[0, M)^k} \mathbf{d}\left(T^u X_n, Y_u\right) d\mathbf{m}(u)\right).  
\end{align*}
In the last step we have used that the distribution of $(X_n, Y)$ is given by $q_n(f|x) \mathcal{P}_*\mu_n(x)$.
 
Therefore 
\begin{align*}
 &  \mathbb{E}\left(\frac{1}{L_n^k}\int_{[0, L_n)^k} \mathbf{d}\left(T^u Z_n, (W_n)_u\right) d\mathbf{m}(u)\right)  \\
  &= \frac{1}{L_n^k M^k} 
  \int_{[0, M)^k} \mathbb{E}\left(\int_{[0, L_n)^k} \mathbf{d}\left(T^u Z_n, W_n(v)_u\right) d\mathbf{m}(u)\right) d\mathbf{m}(v) \\
  & \leq \frac{C^\prime}{L_n} + \tau + \mathbb{E}\left(\frac{1}{M^k} \int_{[0, M)^k} \mathbf{d}\left(T^u X_n, Y_u\right) d\mathbf{m}(u)\right).
\end{align*}
The third term in the last line is smaller than $\varepsilon -2\tau$ for large $n$ by 
(\ref{eq: averaged distance between X_n and Y}).
Thus we conclude that
\[ \mathbb{E}\left(\frac{1}{L_n^k}\int_{[0, L_n)^k} \mathbf{d}\left(T^u Z_n, (W_n)_u\right) d\mathbf{m}(u)\right) < \varepsilon \]
for sufficiently large $n$.   
\end{proof}

\begin{claim}  \label{claim: lower bound on the mutual information between X_n and Y}
For all $n$ we have
\[  \frac{1}{L_n^k} I(Z_n; W_n)  \leq   \frac{1}{M^k} I(X_n;Y). \]
\end{claim}

\begin{proof}
The proof is based on Lemma \ref{lemma: subadditivity of mutual information} (the subadditivity of mutual information)
and Proposition \ref{prop: concavity and convexity of mutual information} (the convexity and concavity of mutual information).

We have $I(Z_n; W_n) = I(\nu_n, \sigma_n)$ and $\sigma_n  = \frac{1}{M^k} \int_{[0, M)^k}\sigma_{n, v}\, d\mathbf{m}(v)$.
By Proposition \ref{prop: concavity and convexity of mutual information} (1) (the convexity of mutual information in 
transition probabilities), 
\begin{equation}  \label{eq: convexity of mutual information in variational principle}
   I(\nu_n, \sigma_n) \leq \frac{1}{M^k} \int_{[0, M)^k} I(\nu_n, \sigma_{n, v})\, d\mathbf{m}(v). 
\end{equation}

Let $W_n(v)$ $(v\in [0, M)^k)$ be the random variable introduced in the proof of Claim \ref{claim: distance between Z and W}.
(Namely, it takes values in $L^1\left([0, L_n)^k, \mathcal{X}\right)$ and its conditional distribution given $Z_n = x$ is 
$\sigma_{n ,v}(f|x)$.)
We have $I(\nu_n, \sigma_{n, v}) = I\left(Z_n; W_n(v)\right)$.
Let $W_n(v)|_{v+\lambda+[0, M)^k}$ $(\lambda \in \Lambda_n)$ and $W_n(v)|_{E_{n,v}}$ be the restrictions of $W_n(v)$ to 
$v+\lambda+[0, M)^k$ and $E_{n,v}$ respectively.
They are conditionally independent given $Z_n$ by the definition of the measure $\sigma_{n, v}$.
By Lemma \ref{lemma: subadditivity of mutual information} (the subadditivity of mutual information),
\[  I\left(Z_n; W_n(v)\right) \leq  I\left(Z_n; W_n(v)|_{E_{n,v}}\right) + \sum_{\lambda\in \Lambda_n} 
     I\left(Z_n; W_n(v)|_{v+\lambda + [0, M)^k}\right). \]
We have 
\[  I\left(Z_n; W_n(v)|_{E_{n,v}}\right) =0  \quad \text{because $W_n(v)|_{E_{n,v}}=x_0$ almost surely}, \]
\[  I\left(Z_n; W_n(v)|_{v+\lambda + [0, M)^k}\right) = I\left(\mathcal{P}(T^{v+\lambda}Z_n); W_n(v)|_{v+\lambda + [0, M)^k}\right)
    = I\left(\mathcal{P}_*T^{v+\lambda}_*\nu_n, q_n\right). \]
Hence $I\left(\nu_n, \sigma_{n, v}\right) \leq  \sum_{\lambda\in \Lambda_n} I\left(\mathcal{P}_*T^{v+\lambda}_*\nu_n, q_n\right)$.
By (\ref{eq: convexity of mutual information in variational principle})
  \begin{align*}
    I(\nu_n, \sigma_n)  & \leq \frac{1}{M^k} \sum_{\lambda\in \Lambda_n} \int_{[0,M)^k} 
    I\left(\mathcal{P}_* T^{v+\lambda}_*\nu_n, q_n\right)\, d\mathbf{m}(v) \\
    & = \frac{1}{M^k} \int_{\bigcup_{\lambda\in \Lambda_n} \left(\lambda+[0, M)^k\right)}  
     I\left(\mathcal{P}_* T^{v}_*\nu_n, q_n\right)\, d\mathbf{m}(v) \\
    & \leq  \frac{1}{M^k} \int_{[0, L_n)^k}   I\left(\mathcal{P}_* T^{v}_*\nu_n, q_n\right)\, d\mathbf{m}(v).
  \end{align*}
We apply Proposition \ref{prop: concavity and convexity of mutual information} (2) (the concavity of mutual information in 
probability measure) and get 
\[   \frac{1}{L_n^k} \int_{[0, L_n)^k}   I\left(\mathcal{P}_* T^{v}_*\nu_n, q_n\right)\, d\mathbf{m}(v)  \leq 
      I\left(\frac{1}{L_n^k}\int_{[0, L_n)^k}\mathcal{P}_*T^v_*\nu_n\, d\mathbf{m}(v), q_n\right)
      = I\left(\mathcal{P}_*\mu_n, q_n\right). \]
Thus 
\[ I(Z_n; W_n) = I(\nu_n, \sigma_n) \leq  \frac{L_n^k}{M^k} I\left(\mathcal{P}_*\mu_n, q_n\right) 
     = \frac{L_n^k}{M^k} I(X_n; Y). \]
\end{proof}

We define a metric $D_n$ on $L^1\left([0, L_n)^k, \mathcal{X}\right)$ by 
\[  D_n(f, g) = \frac{1}{L_n^k} \int_{[0, L_n)^k} \mathbf{d}\left(f(u), g(u)\right)\, d\mathbf{m}(u). \]
We define a map $F_n \colon \mathcal{X}\to L^1\left([0, L_n)^k, \mathcal{X}\right)$ by 
$F_n(x) = \left(T^u x\right)_{u\in [0, L_n)^k}$.
This is an isometric embedding with respect to the metrics $\bar{\mathbf{d}}_{L_n}$ and $D_n$.

Consider the random variable $F_n(Z_n)$. Its distribution is given by the push-forward measure 
$(F_n)_*\nu_n$.
By the estimate (\ref{eq: power law condition revised}), we have 
\[ \mathbb{P}\left(F_n(Z_n)\in A\right) = \nu_n\left(F_n^{-1}(A)\right)
  \leq  \left(\diam\left(F_n^{-1}(A), \bar{\mathbf{d}}_{L_n}\right)\right)^{s L_n^k} \leq  \left(\diam(A, D_n)\right)^{s L_n^k}  \]
for all Borel subsets $A\subset  L^1\left([0, L_n)^k, \mathcal{X}\right)$ with 
$\diam(A, D_n) < \delta$.
By Claim \ref{claim: distance between Z and W} we have $\mathbb{E}\left[D_n\left(F_n(Z_n), W_n\right)\right] < \varepsilon$
for all large $n$.

Noting $2\varepsilon \log(1/\varepsilon) < \delta$, we apply the Kawabata--Dembo estimate 
(Proposition \ref{prop: Kawabata--Dembo estimate}) to the pair $(F_n(Z_n), W_n)$ and get 
\[  I\left(F_n(Z_n); W_n\right) \geq s L_n^k \log(1/\varepsilon) - K(s L_n^k + 1)
    \quad  \text{for large $n$}. \]
Here $K$ is an universal constant.
We have $I\left(F_n(Z_n); W_n\right) = I(Z_n; W_n)$ because $F_n$ is an embedding.
By Claim \ref{claim: lower bound on the mutual information between X_n and Y}
\[  \frac{1}{M^k} I(X_n;Y) \geq  \frac{1}{L_n^k} I(Z_n; W_n)  \geq  s \log(1/\varepsilon) - K\left(s + \frac{1}{L_n^k}\right) 
     \quad  \text{for large $n$}. \]
We let $n\to \infty$.
Recall that $I\left(\mathcal{P}(X); Y\right) = \lim_{n\to \infty} I(X_n;Y)$.
Now we have 
\[  \frac{1}{M^k} I\left(\mathcal{P}(X); Y\right) \geq  s \log(1/\varepsilon) - Ks. \]
Thus we get the estimate (\ref{eq: lower bound on mutual information in dynamical Frostman lemma}):
\[  \frac{1}{M^k} I(X; Y) \geq  \frac{1}{M^k}  I\left(\mathcal{P}(X); Y\right) \geq  s \log(1/\varepsilon) - Ks. \]
This shows that $\frac{1}{M^k}R\left(\varepsilon, [0, M)^k\right) \geq s \log(1/\varepsilon) - Ks$ and hence
$R(\mathbf{d}, \mu, \varepsilon) \geq  s\log(1/\varepsilon) - Ks$.
\end{proof}

Let $(\mathcal{X}, \mathbf{d})$ be a compact metric space with a continuous action 
$T\colon \mathbb{R}^k\times \mathcal{X}\to \mathcal{X}$.
Recall that, for a subset $A\subset \mathbb{R}^k$, we defined a metric $\mathbf{d}_A$ on $\mathcal{X}$ by 
$\mathbf{d}_A(x, y) = \sup_{u\in A} \mathbf{d}\left(T^u x, T^u y\right)$.
In particular, we denoted $\mathbf{d}_L = \mathbf{d}_{[0, L]^k}$ for a positive number $L$.

In the next proposition, we need to consider a $\mathbb{R}^k$-action and its restriction to the discrete subgroup $\mathbb{Z}^k$
simultaneously. 
Recall that, for a natural number $L$, 
we defined metrics $\mathbf{d}^{\mathbb{Z}}_L$ and $\bar{\mathbf{d}}^{\mathbb{Z}}_L$ on $\mathcal{X}$ by 
\[ \mathbf{d}^{\mathbb{Z}}_L(x,y) = \max_{u\in [L]} \mathbf{d}\left(T^u x, T^u y\right), \quad 
    \bar{\mathbf{d}}^{\mathbb{Z}}_L(x, y) = \frac{1}{L^k} \sum_{u\in [L]} \mathbf{d}\left(T^u x, T^u y\right),  \]
where $[L] = \{0,1,2,\dots, L-1\}^k$.

In Proposition \ref{prop: general method to construct measures}, we considered the “$L^1$-metric” $\bar{\mathbf{d}}_L$.
However this is not suitable for the application to Brody curves.
We would like to replace it with the “$L^\infty$-metric” $\mathbf{d}_L$.
The next proposition provides such a result.
This will be used in \S \ref{subsection: proof of Theorem construction of extremal measure}.

\begin{proposition} \label{prop: general method to construct invariant measures in L^infty metric}
Let $(\mathcal{X}, \mathbf{d})$ be a compact metric space with a continuous $\mathbb{R}^k$-action 
$T\colon \mathbb{R}^k\times \mathcal{X}\to \mathcal{X}$. 
Suppose that $(\mathcal{X}, \bar{\mathbf{d}}_1)$ has the tame growth of covering numbers 
and that there exist $a\in \mathbb{R}^k$, $C>0$ and a natural number $\ell$ such that 
     \begin{equation}  \label{eq: L^infty metric and L^infty L^1 metric}
       \mathbf{d}_{a+[0, L]^k}(x, y) \leq C \left(\bar{\mathbf{d}}_1\right)^{\mathbb{Z}}_{L+\ell}(x, y) 
     \end{equation}  
     for all $x, y\in \mathcal{X}$ and all natural numbers $L$.
Let $L_n$ and $M_n$ $(n\geq 1)$ be sequences of positive numbers with $\lim_{n\to \infty} L_n = \infty$.
Let $\nu_n$ be a sequence of Borel probability measures on $\mathcal{X}$. Define 
\[  \mu_n = \frac{1}{L_n^k} \int_{[0,L_n]^k} T^u_*\nu_n\, d\mathbf{m}(u). \]
Assume the following two conditions.
   \begin{enumerate}
      \item There exist positive numbers $\delta$ and $C^\prime$ such that we have
      \[  \nu_n(A) \leq \left(C^\prime\cdot \diam\left(A, \mathbf{d}_{L_n}\right)\right)^{M_n} \]
      for any Borel subset $A\subset \mathcal{X}$ with $\diam\left(A, \mathbf{d}_{L_n}\right) < \delta$.
      \item $\mu_n$ converges to $\mu \in \mathscr{M}^T \left(\mathcal{X}\right)$ as $n\to \infty$ in the weak$^*$ topology.   
   \end{enumerate}
Then we have 
\[   \lrdim\left(\mathcal{X}, T, \mathbf{d}, \mu\right) \geq  \varlimsup_{n\to \infty} \frac{M_n}{L_n^k}.  \]
\end{proposition}

We notice that the metric $\left(\bar{\mathbf{d}}_1\right)^{\mathbb{Z}}_{L+\ell}$ in the right-hand side of 
(\ref{eq: L^infty metric and L^infty L^1 metric}) is given by 
\[ \left(\bar{\mathbf{d}}_1\right)^{\mathbb{Z}}_{L+\ell}(x, y) = \max_{u\in [L+\ell]} \bar{\mathbf{d}}_1\left(T^u x, T^u y\right)
    = \max_{u\in [L+\ell]} \int_{[0,1]^k}\mathbf{d}\left(T^{u+v}x, T^{u+v}y\right) d\mathbf{m}(v). \]

\begin{proof}
By replacing $L_n$ with $\lceil L_n\rceil$, we can assume that $L_n$ are natural numbers.
If $\varlimsup_{n\to \infty} \frac{M_n}{L_n^k}$ is zero then the statement is trivial. So we assume that it is positive.
By taking a subsequence, we can also assume that the limit $\lim_{n\to \infty} \frac{M_n}{L_n^k}$ exists.
Let $s$ and $\eta$ be arbitrary positive numbers with $s< \lim_{n\to \infty} \frac{M_n}{L_n^k}$ and $0<\eta<1$.
We will show that $\lrdim\left(\mathcal{X}, T, \mathbf{d}, \mu\right) \geq (1-\eta)s$.
Once this is proved, we get the result by letting $s\to  \lim_{n\to \infty} \frac{M_n}{L_n^k}$ and $\eta \to 0$.

There is a natural number $n_0$ such that we have $M_n > s (L_n+\ell)^k$ for $n\geq n_0$.
We assume $n\geq n_0$ below.

Set $\nu^\prime_n = T^{-a}_*\nu_n$.
We have 
\[  \frac{1}{(L_n+\ell)^k} \int_{[0, L_n+\ell]^k} T^u_*\nu_n^\prime \, d\mathbf{m}(u)  \to \mu  \quad 
     \text{as $n\to \infty$, in the weak$^*$ topology}. \]

If a Borel subset $A$ of $\mathcal{X}$ satisfies $\diam\left(A, \mathbf{d}_{a+[0,L_n]^k}\right) < \delta$,
then $\diam\left(T^a A, \mathbf{d}_{[0,L_n]^k}\right) < \delta$ and hence 
\[
  \nu_n^\prime(A) = \nu_n\left(T^a A\right)  \leq \left(C^\prime \cdot \diam\left(T^a A, \mathbf{d}_{[0,L_n]^k}\right)\right)^{M_n} 
   = \left(C^\prime \cdot \diam\left(A, \mathbf{d}_{a+[0,L_n]^k}\right)\right)^{M_n}.
\]

Set $D = \bar{\mathbf{d}}_1$. We have $\bar{D}_L^{\mathbb{Z}} = \bar{\mathbf{d}}_L$ for any natural number $L$.
Suppose a Borel subset $A\subset \mathcal{X}$ satisfies 
$\diam \left(A, D^{\mathbb{Z}}_{L_n+\ell}\right) < \min\left(\frac{\delta}{C}, \frac{1}{CC^\prime}\right)$.
By (\ref{eq: L^infty metric and L^infty L^1 metric}) we have
 $\diam\left(A, \mathbf{d}_{a+[0, L_n]^k}\right) < \delta$ and hence 
\begin{align*}
   \nu_n^\prime (A) & \leq \left(C^\prime \cdot \diam\left(A, \mathbf{d}_{a+[0,L_n]^k}\right)\right)^{M_n} \\
   & \leq  \left(C C^\prime \cdot \diam\left(A, D^{\mathbb{Z}}_{L_n+\ell}\right)\right)^{M_n}  \quad 
   \text{by (\ref{eq: L^infty metric and L^infty L^1 metric})} \\
   & \leq   \left(CC^\prime \cdot \diam\left(A, D^{\mathbb{Z}}_{L_n+\ell}\right)\right)^{s(L_n+\ell)^k}.
\end{align*}
In the last step we have used $C C^\prime \cdot \diam \left(A, D^{\mathbb{Z}}_{L_n+\ell}\right)  < 1$ and $M_n > s(L_n+\ell)^k$.

Let $\varepsilon_0 = \varepsilon_0(\eta, \mathcal{X}, D)$ be the positive number introduced in 
Lemma \ref{lemma: tame growth of covering numbers and Z^k actions}.
By Lemma \ref{lemma: tame growth of covering numbers and Z^k actions}, if a Borel subset $A\subset \mathcal{X}$
satisfies $\diam\left(A, \bar{D}^{\mathbb{Z}}_{L_n+\ell}\right) < 
\min\left(\varepsilon_0, \left(\frac{\delta}{2C}\right)^{\frac{1}{1-\eta}}, \left(\frac{1}{2C C^\prime}\right)^{\frac{1}{1-\eta}}\right)$
then we have 
\[  \nu^\prime_n(A) \leq  \left\{4^{\frac{1}{s(1-\eta)}} \left(2CC^\prime \right)^{\frac{1}{1-\eta}}  
     \diam\left(A, \bar{D}^{\mathbb{Z}}_{L_n+\ell}\right)\right\}^{(1-\eta)s(L_n+\ell)^k}.  \]
Since $\bar{D}^{\mathbb{Z}}_{L_n+\ell} = \bar{\mathbf{d}}_{L_n+\ell}$, we also have 
\[ \nu^\prime_n(A) \leq  \left\{4^{\frac{1}{s(1-\eta)}} \left(2CC^\prime \right)^{\frac{1}{1-\eta}}  
     \diam\left(A, \bar{\mathbf{d}}_{L_n+\ell}\right)\right\}^{(1-\eta)s(L_n+\ell)^k}  \]
for any Borel subset $A\subset \mathcal{X}$ with 
$\diam\left(A, \bar{\mathbf{d}}_{L_n+\ell}\right) < 
\min\left(\varepsilon_0, \left(\frac{\delta}{2C}\right)^{\frac{1}{1-\eta}}, \left(\frac{1}{2C C^\prime}\right)^{\frac{1}{1-\eta}}\right)$.
Now we can use Proposition \ref{prop: general method to construct measures} and conclude 
$\lrdim\left(\mathcal{X}, T, \mathbf{d}, \mu\right) \geq (1-\eta)s$.
\end{proof}

\section{Construction of extremal measures: Proof of Theorem \ref{theorem: construction of extremal measure}}
\label{section: construction of extremal measures}

The purpose of this section is to prove Theorem \ref{theorem: construction of extremal measure}.
We prepare technical results on the metric structure of the space of Brody curves in 
\S \ref{subsection: metric structure of the space of Brody curves}.
We prove Theorem \ref{theorem: construction of extremal measure} in 
\S \ref{subsection: proof of Theorem construction of extremal measure}.

\subsection{Metric structure of the space of Brody curves}  \label{subsection: metric structure of the space of Brody curves}

We would like to apply Proposition \ref{prop: general method to construct invariant measures in L^infty metric}
to Brody curves.
Proposition \ref{prop: general method to construct invariant measures in L^infty metric} requires two conditions 
on the metric structure of the given $\mathbb{R}^k$-actions.
Namely, given a compact metric space $(\mathcal{X}, \mathbf{d})$ with a continuous action 
$T\colon \mathbb{R}^k\times \mathcal{X}\to \mathcal{X}$, 
it requires that 
$(\mathcal{X}, \bar{\mathbf{d}}_1)$ has the tame growth of covering numbers and that there
exist $a\in \mathbb{R}^k$, $C>0$ and a natural number $\ell$ such that we have
$\mathbf{d}_{a+[0,L]^k}(x, y) \leq C\left(\bar{\mathbf{d}}_1\right)^{\mathbb{Z}}_{L+\ell} (x,y)$ for every natural number $L$.
The purpose of this subsection is to show that the metric on the space of Brody curves satisfies them.

Recall that $\mathcal{B}^N$ is the space of Brody curves in $\mathbb{C}P^N$ with the natural action 
$T\colon  \mathbb{C}\times \mathcal{B}^N \to \mathcal{B}^N$ and that 
we defined the metric $\mathbf{d}$ on it by 
\[  \mathbf{d}(f, g) = \max_{z\in [0,1]^2} \dfs\left(f(z), g(z)\right). \]

\begin{lemma} \label{lemma: L^infty and L^infty L^1 metrics in the case of Brody curves}
Let $a=\frac{1}{2} + \frac{\sqrt{-1}}{2}$.
For $f, g \in \mathcal{B}^N$ and any natural number $L$ we have
\[  \mathbf{d}_{a+[0, L]^2}(f, g) \leq 4\left(\bar{\mathbf{d}}_1\right)^{\mathbb{Z}}_{L+1}(f, g). \]
\end{lemma}

\begin{proof}
For every $u\in \left[0, \frac{1}{2}\right]^2$ we have 
\[ \max_{z\in \left[\frac{1}{2}, 1\right]^2} \dfs\left(f(z), g(z)\right) \leq \mathbf{d}\left(T^u f, T^u g\right). \]
Hence 
\[ \max_{z\in \left[\frac{1}{2}, 1\right]^2} \dfs\left(f(z), g(z)\right) 
    \leq  4\int_{\left[0, \frac{1}{2}\right]^2} \mathbf{d}\left(T^u f, T^u g\right)\, d\mathbf{m}(u). \]
Similarly we have 
\begin{align*}
   \max_{\left[1, \frac{3}{2}\right]\times \left[\frac{1}{2}, 1\right]} \dfs\left(f(z), g(z)\right) 
   & \leq 4\int_{\left[\frac{1}{2}, 1\right]\times \left[0, \frac{1}{2}\right]} \mathbf{d}\left(T^u f, T^u g\right)\, d\mathbf{m}(u) \\
   \max_{\left[\frac{1}{2}, 1\right] \times \left[1, \frac{3}{2}\right]} \dfs\left(f(z), g(z)\right) 
   & \leq 4\int_{\left[0, \frac{1}{2}\right] \times \left[\frac{1}{2}, 1\right]} \mathbf{d}\left(T^u f, T^u g\right)\, d\mathbf{m}(u) \\
    \max_{\left[1, \frac{3}{2}\right]^2} \dfs\left(f(z), g(z)\right) 
   & \leq  4\int_{\left[\frac{1}{2}, 1\right]^2} \mathbf{d}\left(T^u f, T^u g\right)\, d\mathbf{m}(u).
\end{align*}
Thus 
\[ \mathbf{d}(T^a f, T^a g) = \max_{\left[\frac{1}{2}, \frac{3}{2}\right]^2} \dfs\left(f(z), g(z)\right)
     \leq 4 \int_{[0,1]^2}\mathbf{d}\left(T^u f, T^u g\right)\, d\mathbf{m}(u) 
     = 4 \bar{\mathbf{d}}_1(f, g). \]
For any natural number $L$ 
\begin{align*}
   \mathbf{d}_{a+[0,L]^2}(f, g)  & = \mathbf{d}_L \left(T^a f, T^a g\right)  \\
   & = \max_{z\in [0, L+1]^2} d_{\mathrm{FS}}\left((T^a f)(z), (T^a g)(z)\right) \\
  &  = \max_{u\in [L+1]} \mathbf{d}\left(T^{a+u}f, T^{a+u}g\right) \\
   & \leq  \max_{u\in [L+1]} 4 \bar{\mathbf{d}}_1\left(T^u f, T^u g\right) \\
   &= 4 \left(\bar{\mathbf{d}}_1\right)^{\mathbb{Z}}_{L+1}(f, g).  
\end{align*}  
\end{proof}

\begin{proposition} \label{prop: Brody curves and tame growth of covering numbers}
  Let $A$ be a bounded subset of $\mathbb{C}$. Then $(\mathcal{B}^N, \mathbf{d}_A)$ has the tame
  growth of covering numbers. In particular, $(\mathcal{B}^N, \bar{\mathbf{d}}_1)$ has the tame growth of covering numbers.
\end{proposition}

\begin{proof}
We have $\bar{\mathbf{d}}_1(f, g) \leq \mathbf{d}_1(f, g)$. Hence the latter statement follows from the former.
The following proof was inspired by the argument of Gromov \cite[pp. 393-394]{Gromov}.

First we note that the diameter and injectivity radius of the Fubini--Study metric of $\mathbb{C}P^N$ are both equal to 
$\sqrt{\pi}/2 = 0.8862\dots$ 
in our normalization\footnote{We normalize the metric so that the standard $\mathbb{C}P^1 \subset \mathbb{C}P^N$ has unit area.}.

We fix a small positive number $r$, say $r = \frac{1}{10}$.
Set $S = [0, r]^2$. We also set $S_n = [-nr, (n+1)r]^2$ for $n\geq 0$. 
We have $S=S_0\subset S_1 \subset S_2 \subset \dots$. 

Let $\delta$ be a positive number much smaller than $r$ and satisfying $r/\delta  \in \mathbb{Z}$.
We will fix $\delta$ in Claim \ref{claim: Cauchy estimate and metric} below.
(Indeed $\delta = 10^{-10}$ is enough if $r= \frac{1}{10}$.)
Set $\Lambda = \mathbb{Z}\delta + \mathbb{Z}\delta\sqrt{-1}$.
This is a lattice of the plane $\mathbb{C}$.
The proof is divided into several claims. 

\begin{claim} \label{claim: Cauchy estimate and metric}
We can choose $\delta>0$ so small that for any $f, g\in \mathcal{B}^N$ 
\[  \max_{z\in S} \dfs\left(f(z), g(z)\right) \leq C_6 \max_{z\in S\cap \Lambda}\dfs\left(f(z), g(z)\right) + 
     \frac{1}{2} \max_{z\in S_1} \dfs\left(f(z), g(z)\right). \]
Here $C_6$ is an universal positive constant. Notice that $S\cap \Lambda$ is a finite set.
\end{claim}

\begin{proof}
If $\dfs\left(f(0), g(0)\right) >\frac{1}{10}$ then 
\[ \max_S \dfs\left(f(z), g(z)\right) \leq \diam \, \mathbb{C}P^N < 1 
    \leq 10 \max_{S\cap \Lambda} \dfs\left(f(z), g(z)\right),  \]
and hence the claim holds for $C_6\geq 10$.
Therefore we assume $\dfs\left(f(0), g(0)\right) \leq \frac{1}{10}$.

By symmetry we can also assume $f(0) = [1:0:\dots:0]$ without loss of generality.
Then both $f(z)$ and $g(z)$ $(z\in S_1)$ belong to the $(1/2)$-neighborhood of $[1:0:\dots:0]$
in $\mathbb{C}P^N$.
Let $f(z) = [1:f_1(z):\dots:f_n(z)]$ and $g(z) = [1:g_1(z):\dots:g_n(z)]$.
All $f_i(z)$ and $g_i(z)$ are holomorphic functions over $S_1$.
Set $F(z) = \left(f_1(z), \dots, f_N(z)\right)$ and $G(z) = \left(g_1(z), \dots, g_N(z)\right)$.
We have 
\[ C^{-1}|F(z)-G(z)| \leq \dfs\left(f(z), g(z)\right) \leq C |F(z)-G(z)|\, \quad 
    (z\in S_1). \]
By Cauchy’s estimate on the derivative,
\[ \norm{F^\prime - G^\prime}_{L^\infty(S)} \leq C^\prime \norm{F-G}_{L^\infty(S_1)}. \]
For any point $z\in S$ there exists $z^\prime \in S\cap \Lambda$ with $|z-z^\prime|\leq \delta$.
Therefore 
\begin{align*}
  \norm{F-G}_{L^\infty(S)} & \leq \norm{F-G}_{\ell^\infty(S\cap \Lambda)} + \delta \norm{F^\prime-G^\prime}_{L^\infty(S)} \\
  & \leq  \norm{F-G}_{\ell^\infty(S\cap \Lambda)} + C^\prime \delta \norm{F-G}_{L^\infty(S_1)}.
\end{align*}
This implies 
\[ \max_S \dfs\left(f(z), g(z)\right) \leq C^2 \max_{S\cap \Lambda} \dfs\left(f(z), g(z)\right) 
    + C^2 C^\prime \delta \max_{S_1} \dfs\left(f(z), g(z)\right). \]
We choose $\delta>0$ with $C^2 C^\prime \delta \leq \frac{1}{2}$.
\end{proof}

\begin{claim} \label{claim: Lipschitz discretization}
For $f, g\in \mathcal{B}^N$ and every natural number $m$ we have 
\[  \max_{z\in S_m} \dfs\left(f(z), g(z)\right)
     \leq  C_6 \sum_{n=m}^\infty 2^{m-n} \max_{z\in S_n\cap \Lambda} \dfs\left(f(z), g(z)\right). \]
\end{claim}

\begin{proof}
$S_1$ consists of 9 squares congruent to $S$.
We apply the argument of Claim \ref{claim: Cauchy estimate and metric} to each square and get 
\[  \max_{S_1} \dfs\left(f(z), g(z)\right) \leq C_6 \max_{S_1\cap \Lambda} \dfs\left(f(z), g(z)\right) + 
     \frac{1}{2} \max_{S_2} \dfs\left(f(z), g(z)\right). \]
By repeating this argument, for every $n\geq 0$ we have
\[  \max_{S_n} \dfs\left(f(z), g(z)\right) \leq C_6 \max_{S_n\cap \Lambda} \dfs\left(f(z), g(z)\right) + 
     \frac{1}{2} \max_{S_{n+1}} \dfs\left(f(z), g(z)\right). \]
Combining these estimates, we get the claim.
\end{proof}

Let $(\mathbb{C}P^N)^{\Lambda}$ be the infinite product of the copies of $\mathbb{C}P^N$ indexed by $\Lambda$.
We define a metric $D$ on $\left(\mathbb{C}P^N\right)^{\Lambda}$ by 
\[ D(u, v) = \sum_{n=0}^\infty 2^{-n} \max_{\lambda\in S_n \cap \Lambda} \dfs\left(u_\lambda, v_\lambda\right). \] 

Let $A$ be a bounded subset of $\mathbb{C}$. We choose a natural number $m$ so that $S_m$ contains $A+[0,1]^2$.
Then $\mathbf{d}_A(f, g) \leq \max_{S_m} \dfs\left(f(z), g(z)\right)$.
By Claim \ref{claim: Lipschitz discretization} 
\begin{equation}  \label{eq: comparison between d_A and D under the discretezation}
   \mathbf{d}_A(f, g) \leq 2^m C_6 \sum_{n=m}^\infty 
   2^{-n} \max_{z\in S_n\cap \Lambda} \dfs\left(f(z), g(z)\right)
     \leq 2^m C_6 \cdot D\left(f|_{\Lambda}, g|_{\Lambda}\right). 
\end{equation}     
Here $f|_\Lambda$ and $g|_{\Lambda}$ are the restrictions of $f$ and $g$ to $\Lambda$.
They are elements of $(\mathbb{C}P^N)^{\Lambda}$.

The estimate (\ref{eq: comparison between d_A and D under the discretezation}) shows that the map 
\[ Q \colon  \mathcal{B}^N \to (\mathbb{C}P^N)^{\Lambda}, \quad f\mapsto f|_{\Lambda} \]
is injective and that the inverse map $Q^{-1}\colon Q\left(\mathcal{B}^N\right) \to \mathcal{B}^N$ is 
Lipschitz with respect to the metrics $D$ and $\mathbf{d}_A$. 
Therefore it is enough to show the next claim for proving that $(\mathcal{B}^N, \mathbf{d}_A)$ has the 
tame growth of covering numbers.

\begin{claim}
 The space $\left((\mathbb{C}P^N)^{\Lambda}, D\right)$ has the tame growth of covering numbers.
\end{claim}      
 
\begin{proof}
First we recall a terminology: Let $(\mathcal{X}, d)$ be a compact metric space.
For $\varepsilon>0$, a subset $A$ of $\mathcal{X}$ is said to be an \textbf{$\varepsilon$-spanning set} if 
for every $x\in \mathcal{X}$ there exists $a\in A$ with $d(x, a) < \varepsilon$.
If $A$ is an $\varepsilon$-spanning set then $\#\left(X, d, 3\varepsilon\right) \leq |A|$.

Given $0<\varepsilon < 1$.
We can take a $(\varepsilon/10)$-spanning set $A$ of $\mathbb{C}P^N$ with $|A|\leq (C/\varepsilon)^{2N}$.
We choose a natural number $n_0$ so that $\sum_{n\geq n_0} 2^{-n}\diam (\mathbb{C}P^N) < \frac{\varepsilon}{10}$
and $n_0 \leq \log (1/\varepsilon) + C^\prime$.
We have 
\[  |\Lambda \cap S_{n_0}| \leq C^{\prime\prime}  n_0^2 
     \leq C^{\prime\prime} \left(\log (1/\varepsilon) + C^\prime\right)^2. \]
     
Fix a point $p\in \mathbb{C}P^N$.
We define a finite subset $B\subset \left(\mathbb{C}P^N\right)^{\Lambda}$ as the set of $u\in (\mathbb{C}P^N)^{\Lambda}$ satisfying 
$u_\lambda\in A$ for $\lambda\in \Lambda\cap S_{n_0}$ and $u_\lambda = p$ outside of $\Lambda\cap S_{n_0}$.
Then $B$ is a $(\varepsilon/3)$-spanning set of $(\mathbb{C}P^N)^{\Lambda}$. 
We have 
\[  \#\left((\mathbb{C}P^N)^{\Lambda}, D, \varepsilon\right) \leq  |B|  \leq  |A|^{|\Lambda\cap S_{n_0}|} 
    \leq  \left(\frac{C}{\varepsilon}\right)^{2NC^{\prime\prime} \left(\log (1/\varepsilon) + C^\prime\right)^2}. \]
Then for any positive number $\tau$ we have 
\[ \varepsilon^\tau \log \#\left((\mathbb{C}P^N)^{\Lambda}, D, \varepsilon\right)
    \leq \varepsilon^\tau \cdot  2NC^{\prime\prime} \left(\log (1/\varepsilon) + C^\prime\right)^2 \log (C/\varepsilon) 
    \to  0 \quad (\text{as $\varepsilon \to 0$}). \]
\end{proof}
\end{proof}

\subsection{Proof of Theorem \ref{theorem: construction of extremal measure}}
\label{subsection: proof of Theorem construction of extremal measure}

We need to recall some notations of \S \ref{subsection: nondegenerate Brody curves}.
Let $T\mathbb{C}P^N$ be the tangent bundle of $\mathbb{C}P^N$ with the Fubini--Study metric.
Let $f\colon \mathbb{C}\to \mathbb{C}P^N$ be a Brody curve.
Let $E:=f^* T\mathbb{C}P^N$ be the pull-back of $T\mathbb{C}P^N$ by $f$.
This is a holomorphic vector bundle over the plane $\mathbb{C}$.
Its Hermitian metric is given by the pull-back of the Fubini--Study metric.
Let $H_f$ be the space of holomorphic sections $u$ of $E$ satisfying $\norm{u}_{L^\infty(\mathbb{C})} < \infty$.
The pair $(H_f, \norm{\cdot}_{L^\infty(\mathbb{C})})$ is a (possibly infinite dimensional) complex Banach space.
We set $B_r(H_f) = \{u\in H_f \mid  \norm{u}_{L^\infty(\mathbb{C})} \leq r\}$ for $r>0$.

The next two propositions were proved in \cite[Proposition 3.1, Proposition 3.2]{Matsuo--Tsukamoto}.

\begin{proposition} \label{prop: construction of deformation}
Let $f\colon \mathbb{C}\to \mathbb{C}P^N$ be a nondegenerate holomorphic map with 
$\norm{df}_{L^\infty(\mathbb{C})} < 1$. 
Then there exists $r_1>0$ and a map $\mathcal{F} \colon B_{r_1}(H_f)\to \mathcal{B}^N$
such that 
  \begin{enumerate}
    \item  $\mathcal{F}(0) = f$, 
    \item  for any $u, v \in B_{r_1}(H_f)$ and $z\in \mathbb{C}$ we have 
    \[ \left|\dfs\left(\mathcal{F}(u)(z), \mathcal{F}(v)(z)\right) - |u(z)-v(z)|\right| \leq
        \frac{1}{8} \norm{u-v}_{L^\infty(\mathbb{C})}. \]
  \end{enumerate}
\end{proposition}

\begin{proposition}  \label{prop: study of deformation parameters}
Let $f\colon \mathbb{C}\to \mathbb{C}P^N$ be a nondegenerate Brody curve.
For any square $S\subset \mathbb{C}$ of side length $L>2$ there exists a finite dimensional 
complex linear subspace $V\subset H_f$ satisfying the following two conditions.
   \begin{enumerate}
      \item  The complex dimension of $V$ is bounded from below by  
      \[   \dim_{\mathbb{C}} V \geq  (N+1) \int_S |df|^2(z)\, dxdy - C_f L, \]
       where $C_f$ is a positive constant depending only on $f$ and independent of $S$ and $L$.
      \item  For all $u \in V$ we have $\norm{u}_{L^\infty(\mathbb{C})} \leq 2 \norm{u}_{L^\infty(S)}$.
   \end{enumerate}
\end{proposition}

\begin{proof}[Sketch of the proofs of Propositions \ref{prop: construction of deformation} and 
\ref{prop: study of deformation parameters}]
Proposition \ref{prop: construction of deformation} is a consequence of the deformation theory 
described in \S \ref{subsection: nondegenerate Brody curves}.
In the notations of \S \ref{subsection: nondegenerate Brody curves}, the map 
$\mathcal{F}\colon B_{r_1}(H_f) \to \mathcal{B}^N$ is given by 
$\mathcal{F}(u)(z) = \exp_{f(z)}\left(u(z) + \alpha(u)(z)\right)$ for $u\in B_{r_1}(H_f)$.

Proposition \ref{prop: study of deformation parameters} follows from the Riemann--Roch theorem and 
analytic machinery developed for nondegenerate Brody curves (in particular, a right inverse $F$ of the 
Dolbeault operator $\bar{\partial}\colon \Omega^0(E) \to \Omega^{0,1}(E)$ 
given in (\ref{eq: right inverse of the Dolbeault operator})).
\end{proof}

\begin{theorem}  \label{theorem: construction of measure from nondegenerate curve}
Let $f\colon \mathbb{C}\to \mathbb{C}P^N$ be a nondegenerate holomorphic map with 
$\norm{df}_{L^\infty(\mathbb{C})} < 1$.
Then there exists $\mu \in \mathscr{M}^T\left(\mathcal{B}^N\right)$ satisfying 
\[  \rdim\left(\mathcal{B}^N, T, \mathbf{d}, \mu\right) = \int_{\mathcal{B}^N} \psi \, d\mu 
     = 2(N+1) \rho(f). \]
\end{theorem}

\begin{proof}
We apply Proposition \ref{prop: construction of deformation} to $f$.
We take a positive number $r$ such that 
\begin{equation} \label{eq: choice of r in the construction of extremal measure}
   r\leq r_1, \quad  \frac{9r}{8} \leq  \delta_4.
\end{equation}
Here $r_1$ and $\delta_4$ are the positive numbers 
introduced in Proposition \ref{prop: construction of deformation} and 
Lemma \ref{lemma: energy integral is not sensitive} respectively.
By Proposition \ref{prop: construction of deformation} 
there exists a map $\mathcal{F}\colon B_r(H_f)\to \mathcal{B}^N$
such that $\mathcal{F}(0) = f$ and 
$\left|\dfs\left(\mathcal{F}(u)(z), \mathcal{F}(v)(z)\right) - |u(z)-v(z)|\right| \leq
        \frac{1}{8} \norm{u-v}_{L^\infty(\mathbb{C})}$ 
for $u, v\in B_r(H_f)$ and $z\in \mathbb{C}$.
Then for any $u \in B_r(H_f)$
\[  \mathbf{d}_{\mathbb{C}}\left(f, \mathcal{F}(u)\right) \leq \frac{9}{8} \norm{u}_{L^\infty(\mathbb{C})}
    \leq  \frac{9r}{8} \leq \delta_4. \]
By Lemma \ref{lemma: energy integral is not sensitive} we have 
\begin{equation}  \label{eq: deformation has almost the same energy}
    \left|\int_S |df|^2 \, dxdy - \int_S \left|d\left(\mathcal{F}(u)\right)\right|^2 \, dxdy\right|
    < C_5 L 
\end{equation}
for any square $S\subset \mathbb{C}$ of side length $L\geq 1$.

From Corollary \ref{corollary: reformulation of energy density} 
we have $\rho(f) = \lim_{L\to \infty} \left(\frac{1}{L^2} \sup_{a\in \mathbb{C}} \int_{a+[0,L]^2}|df|^2\, dxdy\right)$.
Then there are $a_n\in \mathbb{C}$ and $L_n>0$ 
$(n\geq 1)$ such that 
$L_n\to \infty$ as $n\to \infty$ and 
\[ \rho(f) = \lim_{n\to \infty} \frac{1}{L_n^2} \int_{a_n+[0,L_n]^2}|df|^2\, dxdy. \]
We can assume $L_n>2$ for all $n$.
Set $S_n = a_n + [0, L_n]^2$.
We apply Proposition \ref{prop: study of deformation parameters} to $f$ and $S_n$.
Then there is a finite dimensional complex linear subspace $V_n\subset H_f$ such that 
  \begin{itemize}
    \item  $\dim_{\mathbb{C}} V_n \geq (N+1) \int_{S_n} |df|^2\, dxdy - C_f L_n$, 
    \item  for every $u\in V_n$ we have $\norm{u}_{L^\infty(\mathbb{C})} \leq 2 \norm{u}_{L^\infty(S_n)}$.   
  \end{itemize}
Denote $B_r(V_n) = \{u\in V_n\mid \norm{u}_{L^\infty(\mathbb{C})} \leq  r \}$.
Let $m_n$ be the Lebesgue measure on $B_r(V_n)$ normalized so that $m_n\left(B_r(V_n)\right) = 1$.
For any Borel subset $A\subset B_r(V_n)$ we have 
\begin{equation}  \label{eq: power law of Lebesgue measure}
  m_n(A) \leq  
  \left(\frac{\diam\left(A, \norm{\cdot}_{L^\infty(\mathbb{C})}\right)}{r}\right)^{2\dim_{\mathbb{C}} V_n}
\end{equation}
because $A$ is contained in a ball of radius $\diam\left(A, \norm{\cdot}_{L^\infty(\mathbb{C})}\right)
:= \sup_{u, v\in A} \norm{u-v}_{L^\infty(\mathbb{C})}$.

We define a (not necessarily invariant) Borel probability measure $\nu_n$ on $\mathcal{B}^N$ by
$\nu_n = \left(T^{a_n}\circ \mathcal{F}\right)_* m_n$.
This is the push-forward of $m_n$ by the map 
$T^{a_n}\circ \mathcal{F}\colon B_r(V_n) \to \mathcal{B}^N$.

Let $u, v\in B_r(V_n)$.
For any $z\in \mathbb{C}$
\[ |u(z)-v(z)| \leq  \dfs\left(\mathcal{F}(u)(z), \mathcal{F}(v)(z)\right) 
   +  \frac{1}{8} \norm{u-v}_{L^\infty(\mathbb{C})}. \]
Taking the supremum over $z\in S_n = a_n+[0,L_n]^2$, we get
\[  \norm{u-v}_{L^\infty(S_n)} 
    \leq  \mathbf{d}_{L_n}\left(T^{a_n} \mathcal{F}(u), T^{a_n}\mathcal{F}(v)\right)
    + \frac{1}{8} \norm{u-v}_{L^\infty(\mathbb{C})}. \]
Since $\norm{u-v}_{L^\infty(\mathbb{C})} \leq  2 \norm{u-v}_{L^\infty(S_n)}$, we have 
\[  \norm{u-v}_{L^\infty(\mathbb{C})} 
    \leq \frac{8}{3} \mathbf{d}_{L_n}\left(T^{a_n} \mathcal{F}(u), T^{a_n}\mathcal{F}(v)\right). \]
Hence for any subset $A\subset \mathcal{B}^N$ we have 
\[ \diam\left(B_r(V_n) \cap \left(\left(T^{a_n}\circ \mathcal{F}\right)^{-1}A\right), \norm{\cdot}_{L^\infty(\mathbb{C})}\right)
     \leq  \frac{8}{3}\cdot \diam\left(A, \mathbf{d}_{L_n}\right). \]
By (\ref{eq: power law of Lebesgue measure})
\[ \nu_n(A) = m_n\left( \left(T^{a_n}\circ \mathcal{F}\right)^{-1}A\right) 
    \leq  \left(\frac{8}{3r}\diam\left(A, \mathbf{d}_{L_n}\right)\right)^{2\dim_{\mathbb{C}} V_n}. \]
Therefore if $\diam\left(A, \mathbf{d}_{L_n}\right) \leq \frac{3r}{8}$ then we have 
\begin{equation}  \label{eq: power law of nu_n in construction of extremal measure}
 \nu_n(A) \leq 
  \left(\frac{8}{3r}\diam\left(A, \mathbf{d}_{L_n}\right)\right)^{2(N+1)\int_{S_n}|df|^2\, dxdy - 2C_f L_n}. 
\end{equation}

Define a Borel probability measure $\mu_n$ on $\mathcal{B}^N$ by 
\[ \mu_n = \frac{1}{L_n^2} \int_{[0, L_n]^2} T^u_*\nu_n\, d\mathbf{m}(u). \]
By taking a subsequence (also denoted by $\mu_n$), we can assume that $\mu_n$
converges to some $\mu\in \mathscr{M}^T\left(\mathcal{B}^N\right)$ in the weak$^*$ topology.

Now we would like to use Proposition \ref{prop: general method to construct invariant measures in L^infty metric}.
By Lemma \ref{lemma: L^infty and L^infty L^1 metrics in the case of Brody curves} and 
Proposition \ref{prop: Brody curves and tame growth of covering numbers}, the metric structure of $\mathcal{B}^N$
satisfies the assumptions of Proposition \ref{prop: general method to construct invariant measures in L^infty metric}.
Therefore we can apply Proposition \ref{prop: general method to construct invariant measures in L^infty metric} to $\mu$ and get 
\begin{align*}
   \lrdim\left(\mathcal{B}^N, T, \mathbf{d}, \mu\right)  & \geq \varlimsup_{n\to \infty} \frac{1}{L_n^2}
   \left(2(N+1)\int_{S_n} |df|^2\, dxdy - C_f L_n\right) \\
   & = 2(N+1)\rho(f).
\end{align*}

Next we will show that 
   \begin{equation}  \label{eq: integral of psi over mu and energy density}
       \int_{\mathcal{B}^N} \psi \, d\mu = 2(N+1)\rho(f).
   \end{equation}
Once this is proved, we have 
\[ \urdim\left(\mathcal{B}^N, T, \mathbf{d}, \mu\right) \leq \int_{\mathcal{B}^N} \psi \, d\mu 
   = 2(N+1)\rho(f) \]
by Theorem \ref{theorem: Ruelle inequality for Brody curves}, and hence we can conclude 
\[  \rdim\left(\mathcal{B}^N, T, \mathbf{d}, \mu\right) = \int_{\mathcal{B}^N} \psi \, d\mu 
   = 2(N+1)\rho(f).  \]
Thus it is enough to show (\ref{eq: integral of psi over mu and energy density}).

From the definitions of $\mu$ and $\mu_n$ 
\begin{align*}
  \int_{\mathcal{B}^N} \psi \, d\mu & = \lim_{n\to \infty} \int_{\mathcal{B}^N} \psi \, d\mu_n \\
  &= \lim_{n\to \infty}  \int_{\mathcal{B}^N} 
   \left(\frac{1}{L_n^2} \int_{[0,L_n]^2} \psi(T^u g)\, d\mathbf{m}(u)\right) d\nu_n(g) \\
  & =  \lim_{n\to \infty}  \int_{\mathcal{B}^N} 
   \left(\frac{2(N+1)}{L_n^2}\int_{[0,L_n]^2}|dg|^2 \, dxdy\right) d\nu_n(g).
\end{align*}
By (\ref{eq: deformation has almost the same energy}), for $\nu_n$-almost every $g$ we have
\[ \left|\int_{[0,L_n]^2}|dg|^2 \, dxdy - \int_{a_n+[0,L_n]^2} |df|^2\, dxdy\right| < C_5 L_n. \]
Thus 
\[  \int_{\mathcal{B}^N} \psi \, d\mu  
    = \lim_{n\to \infty} \frac{2(N+1)}{L_n^2} \int_{a_n+[0,L_n]^2} |df|^2\, dxdy
    = 2(N+1)\rho(f). \]
\end{proof}

Now we can prove Theorem \ref{theorem: construction of extremal measure}.
We write the statement again.

\begin{theorem}[$=$ Theorem \ref{theorem: construction of extremal measure}]
For any real number $c$ with $0\leq c < 2(N+1)\rho(\mathbb{C}P^N)$ there exists 
$\mu \in \mathscr{M}^T\left(\mathcal{B}^N\right)$ such that 
\[  \rdim\left(\mathcal{B}^N, T, \mathbf{d}, \mu\right) = \int_{\mathcal{B}^N} \psi \, d\mu = c. \]
\end{theorem}

\begin{proof}
The statement is trivial for $c=0$. (Any delta measure at a constant curve satisfies the statement.)
So assume $0<c<2(N+1)\rho(\mathbb{C}P^N)$.
By Corollary \ref{cor: energy density and nondegenerate curve}, there exists a nondegenerate 
Brody curve $g\colon \mathbb{C}\to \mathbb{C}P^N$ satisfying 
$c < 2(N+1)\rho(g)$. Set 
\[ \lambda = \sqrt{\frac{c}{2(N+1)\rho(g)}}. \]
We have $0< \lambda < 1$. Define $f(z) = g(\lambda z)$.
This is a nondegenerate holomorphic curve with $\norm{df}_{L^\infty} = \lambda \norm{dg}_{L^\infty} < 1$.
Applying Theorem \ref{theorem: construction of measure from nondegenerate curve} to $f$, we can 
construct $\mu \in \mathscr{M}^T\left(\mathcal{B}^N\right)$ satisfying 
\[  \rdim\left(\mathcal{B}^N, T, \mathbf{d}, \mu\right) = \int_{\mathcal{B}^N} \psi \, d\mu 
     = 2(N+1) \rho(f). \]
We have 
\[  \rho(f) = \lambda^2 \rho(g) = \frac{c}{2(N+1)}. \]
\end{proof}

We have finished all the proofs of Theorems \ref{theorem: Ruelle inequality for Brody curves} and 
\ref{theorem: construction of extremal measure}.
Lastly, we briefly discuss an open problem.
This paper is a starting point of the ergodic study of entire holomorphic curves.
There definitely exist plenty of open problems and new research directions.
But here we mention just one problem directly related to the results of this paper.

Theorem \ref{theorem: construction of extremal measure} shows that we can construct a rich variety of invariant probability measures $\mu$
on $\mathcal{B}^N$ satisfying $\rdim\left(\mathcal{B}^N, T, \mathbf{d}, \mu\right) = \int_{\mathcal{B}^N}\psi \, d\mu$.
However it is not easy to detect their ergodic theoretic properties (e.g. ergodicity, mixing properties) from the above proof.
Hence we propose the following problem.

\begin{problem}  \label{problem: ergodicity of extremal measures}
Let $0< c < 2(N+1)\rho(\mathbb{C} P^N)$.
Is there an ergodic measure $\mu\in \mathscr{M}^T(\mathcal{B}^N)$ satisfying 
$\rdim\left(\mathcal{B}^N, T, \mathbf{d}, \mu\right) = \int_{\mathcal{B}^N}\psi \, d\mu = c$?
We notice that this does not (at least, directly) follow from the ergodic decomposition theorem.
Let $\mu \in \mathscr{M}^T(\mathcal{B}^N)$ be the measure given by Theorem \ref{theorem: construction of extremal measure}.
Namely it satisfies $\rdim\left(\mathcal{B}^N, T, \mathbf{d}, \mu\right) = \int_{\mathcal{B}^N}\psi \, d\mu =c$.
Let $\mathscr{M}^T_{\mathrm{erg}}(\mathcal{B}^N)$ be the set of ergodic measures on $\mathcal{B}^N$, and let
\[  \mu = \int_{\mathscr{M}^T_{\mathrm{erg}}(\mathcal{B}^N)} \nu \, d\lambda(\nu) \]
be the ergodic decomposition of $\mu$. 
Here $\lambda$ is a probability measure on $\mathscr{M}^T_{\mathrm{erg}}(\mathcal{B}^N)$.
Then we have 
\[ \urdim\left(\mathcal{B}^N, T, \mathbf{d}, \mu\right) 
    \leq \int_{\mathscr{M}^T_{\mathrm{erg}}(\mathcal{B}^N)} \urdim\left(\mathcal{B}^N, T, \mathbf{d}, \nu\right)\, d\lambda(\nu). \]
Notice that here we consider upper rate distortion dimension. 
For lower rate distortion dimension, we have a similar inequality in the reverse direction\footnote{Namely we have 
\[ \lrdim\left(\mathcal{B}^N, T, \mathbf{d}, \mu\right) 
    \geq  \int_{\mathscr{M}^T_{\mathrm{erg}}(\mathcal{B}^N)} \lrdim\left(\mathcal{B}^N, T, \mathbf{d}, \nu\right)\, d\lambda(\nu). \]}.
We also have 
\[  \int_{\mathcal{B}^N} \psi\, d\mu = \int_{\mathscr{M}^T_{\mathrm{erg}}(\mathcal{B}^N)} 
     \left(\int_{\mathcal{B}^N}\psi\, d\nu\right)d\lambda(\nu). \]
By Theorem \ref{theorem: Ruelle inequality for Brody curves} (a Ruelle inequality for Brody curves), 
\[  \urdim\left(\mathcal{B}^N, T, \mathbf{d}, \nu\right) \leq  \int_{\mathcal{B}^N} \psi \, d\nu \]
for all $\nu\in \mathscr{M}^T_{\mathrm{erg}}(\mathcal{B}^N)$.
It follows from $\rdim\left(\mathcal{B}^N, T, \mathbf{d}, \mu\right) = \int_{\mathcal{B}^N}\psi \, d\mu$ that
\begin{align*}
   \int_{\mathscr{M}^T_{\mathrm{erg}}(\mathcal{B}^N)} \left(\int_{\mathcal{B}^N}\psi\, d\nu\right)d\lambda(\nu)   & \leq 
    \int_{\mathscr{M}^T_{\mathrm{erg}}(\mathcal{B}^N)} \urdim\left(\mathcal{B}^N, T, \mathbf{d}, \nu\right)\, d\lambda(\nu) \\
    & \leq  \int_{\mathscr{M}^T_{\mathrm{erg}}(\mathcal{B}^N)} \left(\int_{\mathcal{B}^N}\psi\, d\nu\right)d\lambda(\nu). 
\end{align*}    
The left-most side and the right-most side are the same.
Hence we conclude that 
\begin{equation} \label{upper rate distortion dimension of ergodic components}
    \urdim\left(\mathcal{B}^N, T, \mathbf{d}, \nu\right) = \int_{\mathcal{B}^N} \psi\, d\nu 
\end{equation}    
for $\lambda$-almost every $\nu\in \mathscr{M}^T_{\mathrm{erg}}(\mathcal{B}^N)$.
Therefore we have obtained ergodic measures $\nu$ 
satisfying the equation (\ref{upper rate distortion dimension of ergodic components}).
However we do not have control of the lower rate distortion dimension of these measures $\nu$.
So we cannot conclude $\rdim\left(\mathcal{B}^N, T, \mathbf{d}, \nu\right) = \int_{\mathcal{B}^N}\psi \, d\nu$.
We even do not know whether there exists a single ergodic measure $\nu$ on $\mathcal{B}^N$ satisfying 
\[  \rdim\left(\mathcal{B}^N, T, \mathbf{d}, \nu\right) = \int_{\mathcal{B}^N} \psi \, d\nu >0. \]
\end{problem}

If Problem \ref{problem: ergodicity of extremal measures} is affirmatively solved, then we can also 
ask the mixing properties of these measures $\mu$.

\end{document}